\documentclass{AIMS}
\usepackage{amsmath}

\usepackage{afterpage}
  \usepackage{paralist}
\usepackage[pdftex]{graphicx}
\usepackage[normal]{subfigure}
\usepackage{amssymb}
\usepackage{subfig,float}
\usepackage{color}
 \usepackage[colorlinks=true]{hyperref}
 \hypersetup{urlcolor=blue, citecolor=red}
 
\def\currentvolume{X}
 \def\currentissue{X}
  \def\currentyear{200X}
   \def\currentmonth{XX}
    \def\ppages{X--XX}
     \def\DOI{10.3934/xx.xx.xx.xx}

\newcommand{\al}[1]{\textcolor{black}{#1}}
\newcommand{\alt}[1]{\textcolor{black}{#1}}

\newcommand{\ds}[0]{\displaystyle}

\newcommand{\eps}[0]{\varepsilon}

\newcommand{\bsub}{\begin{subequations}}
\newcommand{\esub}{\end{subequations}$\!$}

\numberwithin{equation}{section}

\title[Blow up sets of fourth order PDEs]{An Asymptotic Study of Blow up Multiplicity in fourth order parabolic partial differential equations.}


\author[Alan E. Lindsay]{}
\email{a.lindsay@nd.edu}

\begin{document}

\begin{abstract}
Blow-up in second and fourth order semi-linear parabolic partial differential equations (PDEs) is considered in bounded regions of one, two and three spatial dimensions with uniform initial data. A phenomenon whereby singularities form at multiple points simultaneously is exhibited and explained by means of a singular perturbation theory. In the second order case we predict that points furthest from the boundary are selected by the dynamics of the PDE for singularity. In the fourth order case, singularities can form simultaneously at multiple locations, even in one spatial dimension. In two spatial dimensions, the singular perturbation theory reveals that the set of possible singularity points depends subtly on the geometry of the domain and the equation parameters. In three spatial dimensions, preliminary numerical simulations indicate that the multiplicity of singularities can be even more complex. For the aforementioned scenarios, the analysis highlights the dichotomy of \al{behaviors} exhibited between the second and fourth order cases. 
\end{abstract}

\maketitle

\centerline{\scshape Alan E. Lindsay}
\medskip
{\footnotesize
\al{ \centerline{Department of Applied and Computational Mathematics and Statistics,}
   \centerline{University of Notre Dame,}
    \centerline{Notre Dame, IN, USA, 46556.}}
} 

 \centerline{(Communicated by)}

\baselineskip=12pt

\noindent 
\baselineskip=16pt 

\section{Introduction}\label{sec:intro}

Partial differential equations (PDEs) of general form
\bsub\label{intro_1}
\begin{equation}\label{intro_1a}
 \left\{\begin{array}{lc} u_t = \Delta u + f(u), &\quad (x,t)\in\al{\Omega_T}; \\[5pt]
 u = 0, \quad& (x,t)\in\al{\partial\Omega_T};\\[5pt]
 u = \psi(x), &\quad x\in\Omega_{\al{0}},\end{array}\right. 
  \end{equation}
  and their fourth order equivalents
  \begin{equation}\label{intro_1b}
 \left\{\begin{array}{lc} u_t = -\Delta^2 u + f(u), &\quad (x,t)\in\al{\Omega_T}; \\[5pt]
 u =\partial_n u = 0, \quad& (x,t)\in\partial\al{\Omega_{T}};\\[5pt]
 u = \psi(x), &\quad x\in\Omega_{\al{0}},\end{array}\right. 
  \end{equation}
arise in the study of countless physical and natural phenomena. In the above formulations, $\Omega$ is a bounded region of $\mathbb{R}^n$ and
 \begin{equation}\label{timedomain}
\al{ \Omega_{T} = \Omega \times (0,T), \qquad  \partial\Omega_{T} = \partial\Omega \times (0,T).}
 \end{equation} 
 \esub
 In the absence of spatial terms, it is well known (cf. \cite{BB98}) that the ordinary differential equation (ODE) $u_t = f(u) \geq 0$ does not necessarily have \al{a global solution}. Indeed, depending on the form of $f(u)$ and the initial value, a solution may only exist on a finite time interval $(0,\al{T_{0}})$. The solution is said to \emph{blow-up}, in the case where \al{$|u(t)|\to\infty$} as $t \to \al{T_{0}}^{-}$. Alternatively, in the case of \emph{quenching} or \emph{rupture}, $u$ remains finite while $u_t$ diverges as $t\to \al{T_{0}}^{-}$. In either case, the phrase \emph{finite time blow-up} is taken to mean a divergence in the solution or some derivative of the solution in a particular norm at a finite time $\al{T_{0}}$. 

An interesting and long studied problem is to describe the corresponding finite time blow up of the ODE problem $u_t=f(u)$, in the spatial setting described by \alt{problems} \eqref{intro_1}. Classical studies into this area revolve around five central questions; 1. Does a singularity occur? 2. When do the singularities occur? 3. Where does the singularity occur? 4. How do singularities occur? and 5. What happens after a singularity occurs? In the present work we focus on the third of the aforementioned questions: when singularities occur, what is their location and multiplicity? 

According to their \al{relevance} in applications, the above questions have been addressed extensively for \alt{problem} \eqref{intro_1a} with exponential nonlinearities $\al{f(u) = e^u}$ (\cite{Frank38,ZBLM85,AFMC85,BE,Lacey84}), power nonlinearities $\al{f(u) = u^p}$, $p>1$ (\cite{FK92,AFMC85,Kap63,Fujita66,Fujita68,VelGAPoHer93,GaVa95,Levine1990,Levine2000}) and inverse power nonlinearities $\al{f(u) = u^{-p}}$ (\cite{GPW,GUO2,Pel07}). The question of existence of global solutions to the initial boundary value problem
\begin{equation}\label{intro_up}
 \left\{\begin{array}{lc} u_t = \Delta u + u^p, &\quad (x,t)\in\Omega\times(0,\infty); \\[5pt]
 u = 0, \quad& (x,t)\in\partial\Omega \times(0,\infty);\\[5pt]
 u = \psi(x), &\quad x\in\Omega, \quad t = 0,\end{array}\right. 
\end{equation}
is now well known (cf. \cite{BandleLevine89,Fujita66,Fujita68,Levine2000}) in terms of the theory of critical exponents. If $1<p \leq p_c(n)
\equiv 1+2/n$, then $u=0$ is the only global solution of 
\eqref{intro_up}. If $p>p_c(n)$, then global solutions of \eqref{intro_up} can exist, provided the initial data is sufficiently small (cf. \cite{Levine2000}). For further information regarding the large literature in studies of \eqref{intro_1a}, the interested reader is directed to \cite{BB98,GALVA02,BE,Levine1990,Levine2000} and the references therein.

The corresponding fourth order problem \eqref{intro_1b}, which had previously attracted somewhat less attention, has more recently enjoyed significant interest. In \cite{GALAK2002}, the critical exponent $p_c(n) = 1 + 4/n$ was established for \alt{problem} \eqref{intro_1b} with the nonlinearity \al{$f(u) = |u|^p$} and $\Omega= \mathbb{R}^n$. The existence and stability of blow-up profiles to \eqref{intro_1b} with power and exponential nonlinearities was investigated in \cite{GJFW,GALAK,BGW} where it was found that, in contrast to the second order problem, there exists \al{stable} self-similar singularity profiles. In the context of Micro-Electro Mechanical Systems (MEMS), where for $f(u) = 1/(1-u)^2$ \alt{problem} \eqref{intro_1b} models the deflection of a beam/plate under Coulomb forcing, quenching solutions have been studied in \cite{GUO,AL2012,ALS2013}. Related higher order problems in the context of thin film dynamics have been studied extensively (\cite{BW1,BW2,BW3,BW4,BGW2001}) where de-wetting processes occur via ruptures in degenerate fourth order parabolic equations.

Fourth order problems exhibit many interesting and surprising solution features when compared to their second order counterparts. As a simple example, consider the reduced setting where $\Omega = \mathbb{R}^n$ and $f(u)=0$, and the evolution \alt{problem}
\begin{equation}\label{intro_poly1}
\left\{ \begin{array}{lc} 
u_t  + \Delta^2 u = 0, & \qquad (x,t) \in \mathbb{R}^n \times(0,T);\\[5pt]
u = \psi(x), & \qquad x\in\mathbb{R}^n, \quad t=0.
\end{array}
\right.
\end{equation}
 If $\psi(x) \in C^0 \cap L^{\infty}(\mathbb{R}^n)$, then the unique global solution (cf. \cite{Gazzola1,Gazzola2}) of \eqref{intro_poly1} is given by
\bsub\label{intro_poly2}
\begin{equation}\label{intro_poly2a}
u(x,t) = C t^{-n/4} \int_{\mathbb{R}^n} \psi(x-y) \al{k_n}\left( \frac{|y|}{t^{1/4}}\right)\, dy
\end{equation}
where $C= C_n$ are normalization constants and $\al{k_n(z)}$ are the fundamental solutions
\begin{equation}\label{intro_poly2b}
\al{k_n(z)} = z^{1-n} \int_{0}^{\infty} e^{-s^4}  (zs)^{n/2} J_{(n-2)/2}(zs)\, ds.
\end{equation}
\esub
In contrast to the \al{G}aussian kernels of the corresponding second order heat equation, the Bessel functions $J_{\nu}$ in the integrand of \eqref{intro_poly2b} generate highly oscillatory solution behavior. Consequently, many second order features such as the positivity preserving property ($\psi>0$ implies $u>0$) and the maximum principle do not extend to higher order \alt{problems}.

\begin{figure}[htbp]
\centering
\subfigure[$L=5$]{\includegraphics[width=0.45\textwidth]{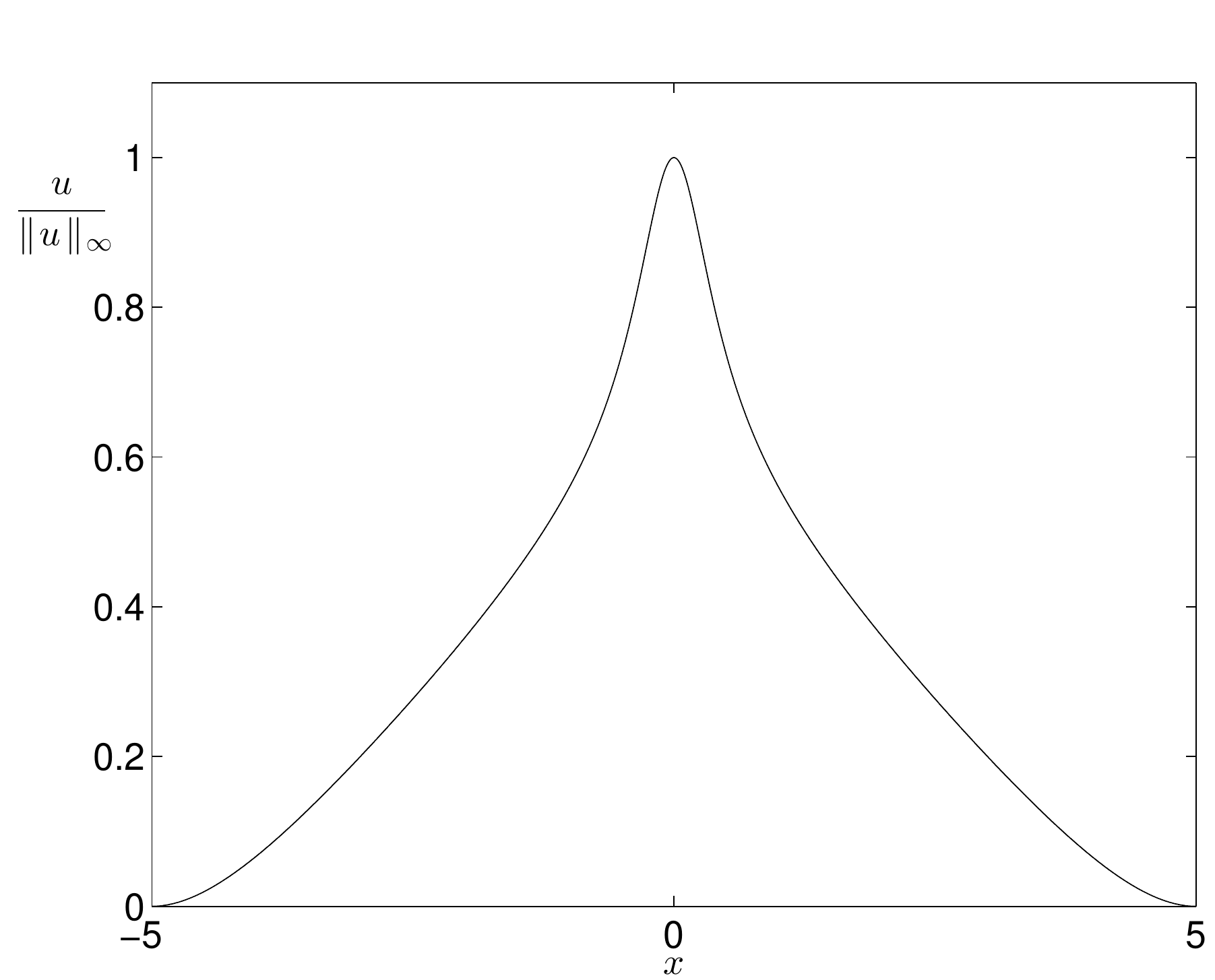}\label{fig:intro1a}}
\qquad
\subfigure[$L=7$]{\includegraphics[width=0.45\textwidth]{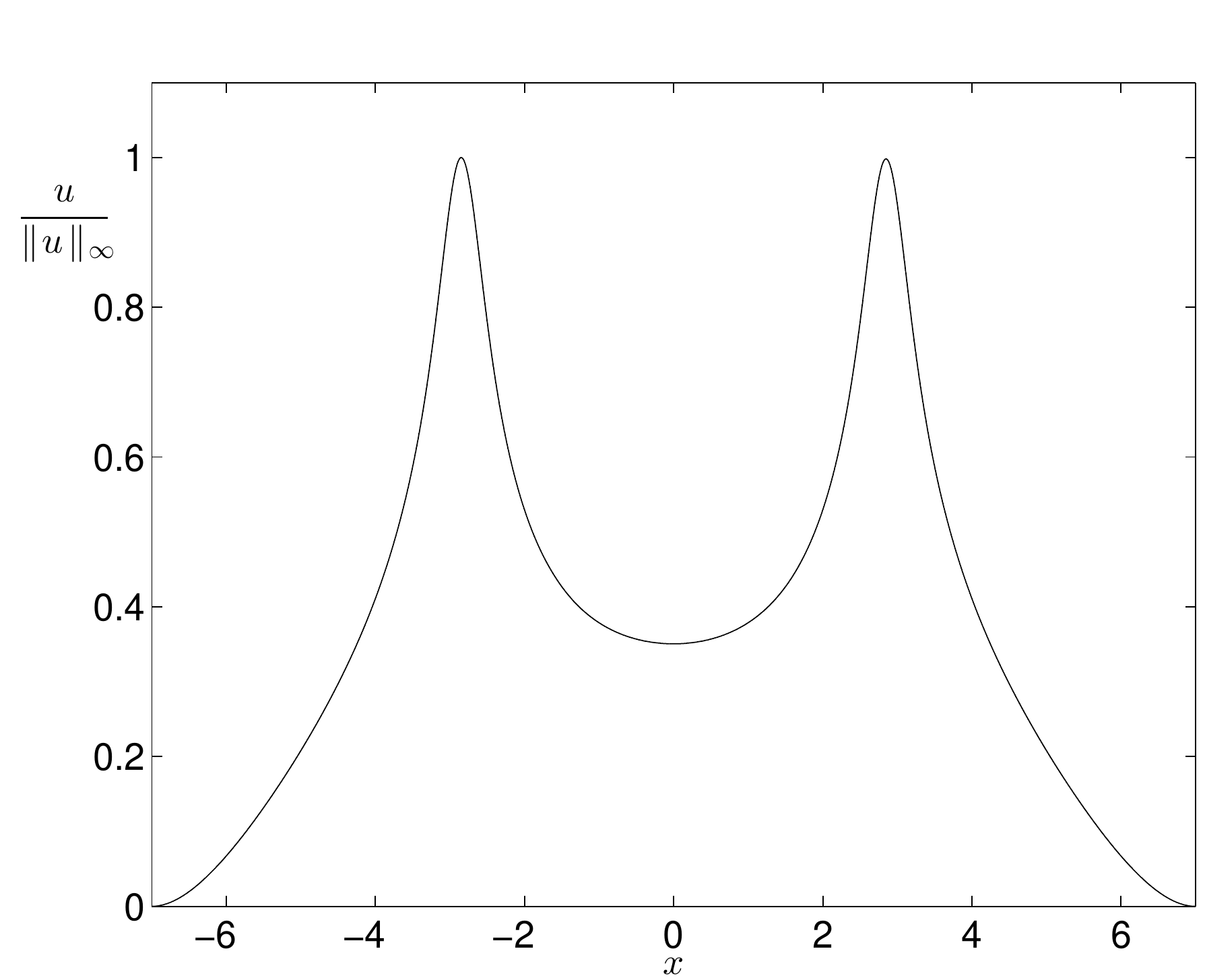}\label{fig:intro1b}}
\parbox{5in}{\caption{Numerical solutions of \alt{problem} \eqref{intro_1b} with $f(u)=e^u$ on the interval $\Omega=[-L,L]$ for $L=5$ and $L=7$ integrated from uniform zero initial data to $\| u\|_{\infty} = 10$. The simulations suggest the existence of a critical value $L = L_c$, over which the multiplicity of singularities changes from one to two. \label{fig:intro1}}}
\end{figure}

As a result of these oscillatory features, the blow up dynamics for the fourth order \alt{problem} \eqref{intro_1b} are quite different to that of the well-studied second order problem \eqref{intro_1a}. To demonstrate the contrasting \al{behaviors} of these two problems, consider the example case of the 1D strip $\Omega = [-L,L]$, with nonlinearity $f(u) = e^u$ and uniform zero initial data $\psi(x) = 0$. In the second order case, it is well known (cf. \cite{AFMC85}) that for $L$ large enough, no equilibriums solutions are present and, from maximum principle considerations, that the solution blows up uniquely at the origin. In Fig.~\ref{fig:intro1}, numerical solutions of the equivalent fourth order problem \eqref{intro_1b} are displayed for values $L=5$ and $L=7$ and are observed to be very different from the previously mentioned second order \al{behavior}. For $L=5$, the blow-up point is observed to occur uniquely at the origin while for $L=7$, we observe two singularities forming simultaneously at distinct points.
This multiple singularity phenomenon of \eqref{intro_1b} was recently observed in \cite{AL2012} for the MEMS case $f(u) = 1/(1-u)^2$ and radially symmetric solutions in one and two dimensions. In the 2D radially symmetric case, the singularities form simultaneously along a ring of points for $L$ sufficiently large. In this particular application, the singularities indicate the contact points between two elastic surfaces and so their multiplicity and location is of practical importance. 

\begin{figure}[htbp]
\centering
\subfigure[$L=1.5$]{\includegraphics[width=0.38\textwidth]{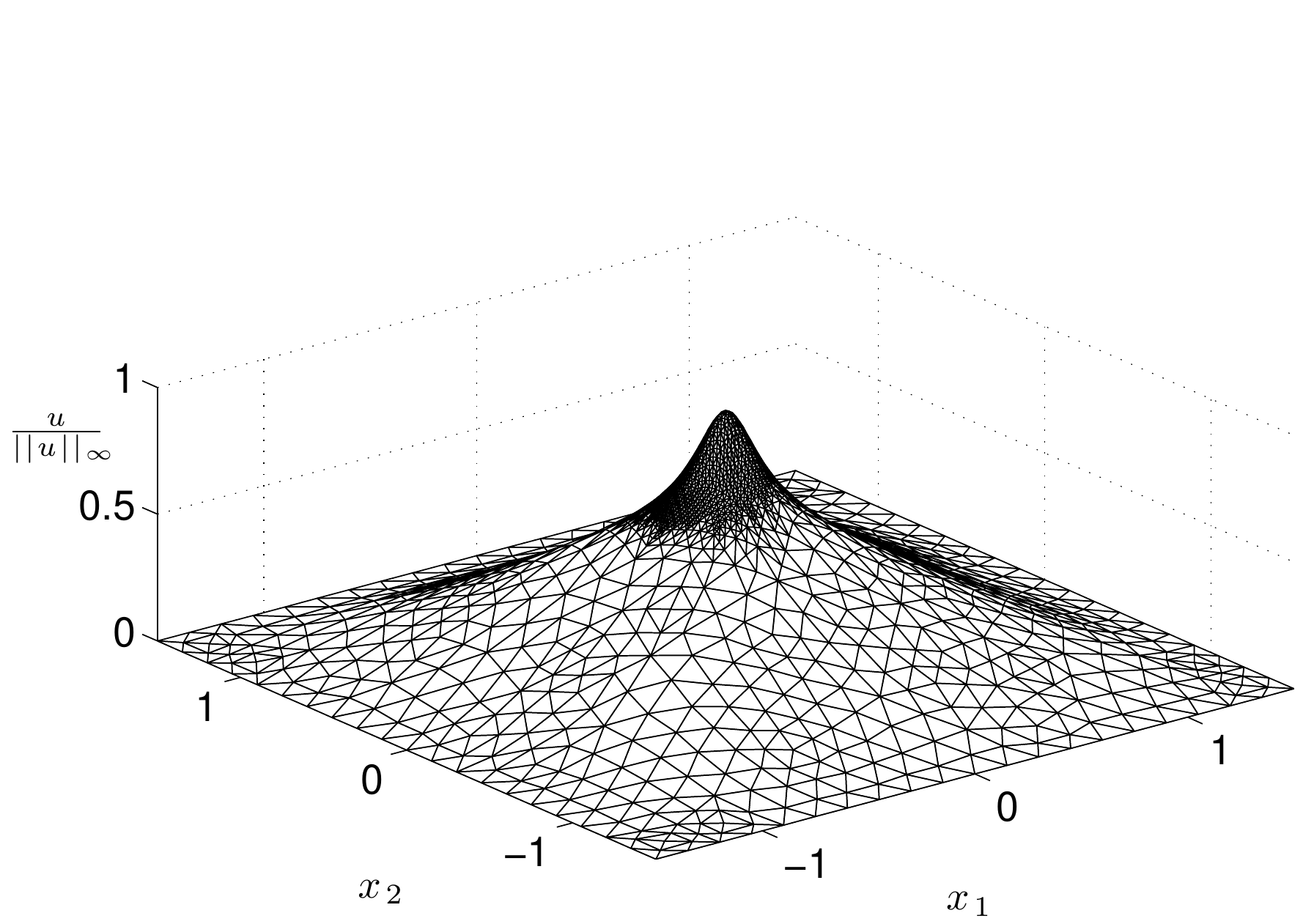}\label{fig:intro2a}}
\qquad
\subfigure[$L=1.8$]{\includegraphics[width=0.42\textwidth]{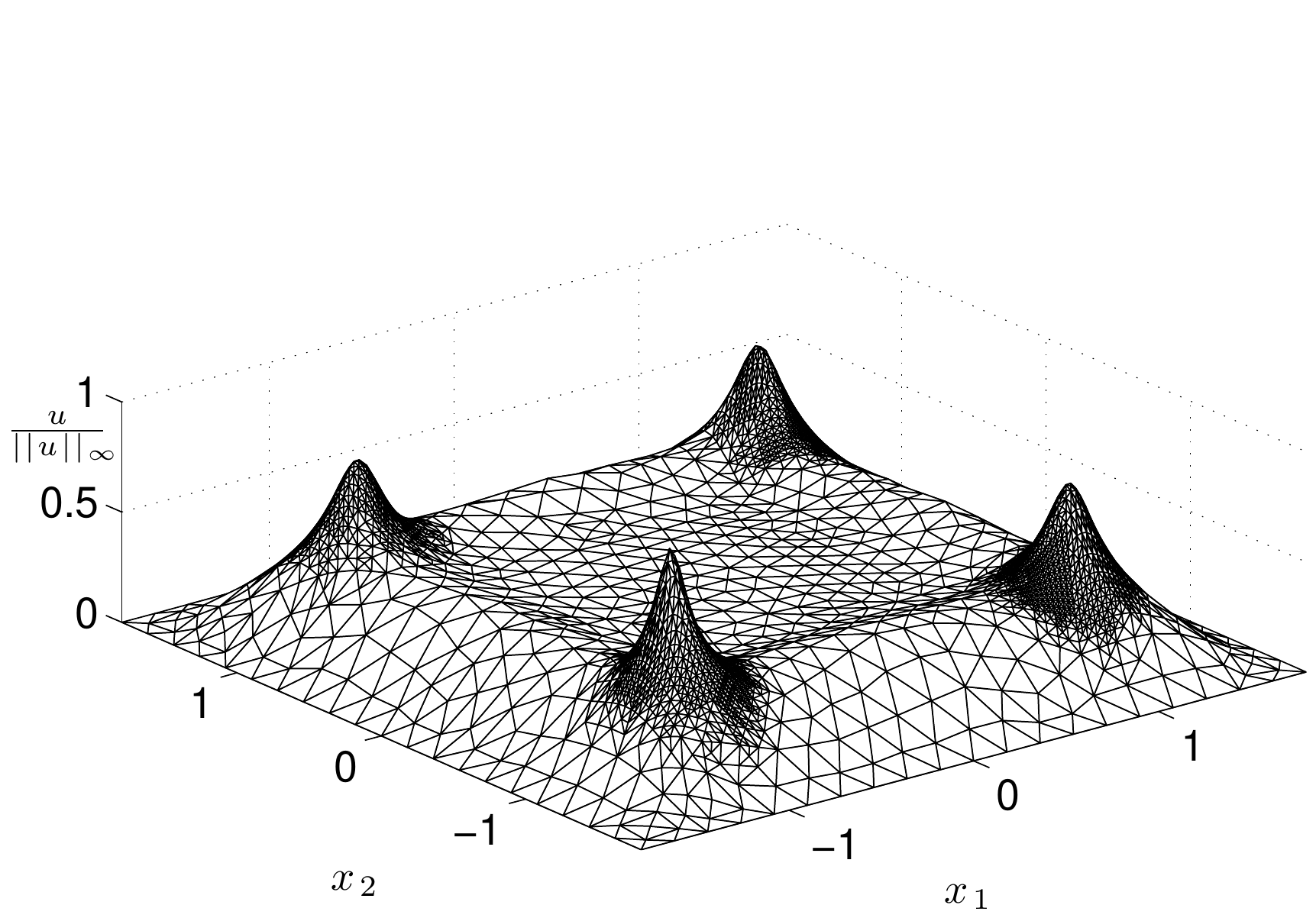}\label{fig:intro2b}}
\parbox{5.6in}{\caption{Numerical solutions of \alt{problem} \eqref{intro_1b} with $f(u)=e^u$ on the square region $\Omega=[-L,L]^2$ for $L=1.5$ and $L=1.8$ integrated from uniform zero initial data to $\al{\|u\|_{\infty}} = 10$. The simulations suggest the existence of a critical value $L = L_c$, over which the multiplicity of singularities changes from one to four. The dark regions near the peaks indicate high grid refinement in the vicinity of the forming singularities. \label{fig:intro2}}}
\end{figure}

The multiplicity of quenching singularities for the inverse square nonlinearity has also been investigated (cf. \cite{ALS2013}) for general 2D geometries. As before, an illustration of the dichotomy between the second \eqref{intro_1a} and fourth \eqref{intro_1b} order cases is provided through the $f(u) = e^u$ case and the square region $\Omega = [-L,L]^2$. For $L$ sufficiently large, \alt{problems} \eqref{intro_1} blow-up in finite time. In the second order case, solutions of \eqref{intro_1a} blow-up uniquely at the origin. As indicated by the numerical simulations of \eqref{intro_1b} shown in Fig.~\ref{fig:intro2}, the multiplicity of singularities is remarkably different for $L=1.5$ and $L=1.8$. In the former case, the singularity occurs uniquely at the origin, while in the later case, the singularity occurs simultaneously at four distinct points.

In the recent works of \cite{AL2012,ALS2013}, the location and multiplicity of singularities \alt{in problem \eqref{intro_2b}} was studied for the nonlinearity $f(u) = 1/(1-u)^2$. In the present work, the contribution is threefold: 

First, we consolidate and generalize the results of \cite{AL2012,ALS2013} to the fourth order problem \eqref{intro_1b} with general positive convex nonlinearities in one and two dimensions. The result of this analysis is a geometric framework which predicts the multiple singularity phenomenon in \eqref{intro_1b} for a general class of nonlinearities $f(u)$. The analysis demonstrates that this phenomenon is due to a combination of boundary effects and the dynamics of fourth order problems which do not admit a maximum principle. In addition, a geometric framework is established for predicting the possible blow up set of \eqref{intro_1b}, for general regions $\Omega$.

Second, we apply the same geometric framework developed herein to the second order problem \eqref{intro_1a} with general nonlinearities. This analysis elucidates the underlying reasons why the multiple singularity phenomena does not occur generically in \eqref{intro_1a}. In addition, it provides a simple geometric framework for predicting the singularity location in such problems. Consequently, we establish an understanding of how the location of singularities in \alt{problems} \eqref{intro_1} is determined by the geometry of $\Omega$ and the dynamics of the PDE in one and two dimensions.

Third and finally, we present numerical simulations of \eqref{intro_1b} in three dimensions to gain insight into the corresponding multiple singularities phenomenon in higher spatial dimensions. The preliminary simulations suggest that the multiplicity of singularities in \eqref{intro_2b} may be greater than in 2D. For example, in the cube $\Omega = [-L,L]^3$, numerical simulation show the number of blow up points changes from eight to one as $L$ increases through $L_c$. A natural conjecture based on observations is the following:  There exists an $L_c(n)$ such that \eqref{intro_1b} for $\Omega=[-L,L]^n$ and $\psi=0$, exhibits blow-up at $2^n$ distinct points whenever $L>L_c(n)$.

 To analyze the multiple singularity phenomenon, and thus obtain predictions of the multiplicity and location of singularities in the system \eqref{intro_1}, we employ formal asymptotic methods in the limit of large domain size. If the length scale of the domain $\Omega$ is $L$, then the problem of interest, as motivated by the previous examples, is to understand the singularity set for large $L$, or equivalently small $\eps = L^{-1}$. When the rescaling $\Omega \to \eps^{-1} \Omega$ is applied to \eqref{intro_1}, the second order singularly perturbed initial boundary value problem
\bsub\label{intro_2}
\begin{equation}\label{intro_2a}
 \left\{\begin{array}{lc} u_t = \eps^2\Delta u + f(u), & \quad \al{(x,t)\in\Omega_{\al{T_{\eps}}}}; \\[5pt]
 u = 0, \quad& \al{(x,t)\in\partial\Omega_{\al{T_{\eps}}}};\\[5pt]
 u = 0, &\quad \al{(x,t)\in\Omega_{0}},\end{array}\right.
  \end{equation}
  and its fourth order equivalent
  \begin{equation}\label{intro_2b}
 \left\{\begin{array}{lc} u_t = -\eps^4 \Delta^2 u + f(u), &\quad \al{(x,t)\in\Omega_{\al{T_{\eps}}}}; \\[5pt]
 u =\partial_n u = 0, \quad& \al{(x,t)\in\partial\Omega_{\al{T_{\eps}}}};\\[5pt]
 u = 0, &\quad \al{(x,t)\in\Omega_{0}},\end{array}\right. 
  \end{equation}
are arrived at. For the purposes of the present work, general positive, convex source terms $f$ satisfying
\begin{equation}\label{intro_2c}
f \in C^{1}, \quad f(t)>0 \qquad t\geq 0, \qquad \al{f(0)=1},
\end{equation}
 \esub 
 are considered. \al{In this particular formulation, \alt{problem} \eqref{intro_2a} is a generic model for a slowly diffusing concentration field with local reaction kinetics $f(u)$. The fourth order \alt{problem} \eqref{intro_2b} is a ubiquitous model for deflection of a plate with small flexural rigidity \cite{PB,AL2012}, undergoing forcing $f(u)$.}
 
 For demonstration of the theory developed herein, classical examples of form $f(t) = e^t$, $f(t) = (1+ t)^p$ for $p>1$ will be used. \al{The premise of our approach is to construct explicit short time solutions to \eqref{intro_2} in the limit $\eps\to0$, the critical points of which act as surrogates for the singularity location(s) of \eqref{intro_2} for one dimensional strips and general bounded regions of $\mathbb{R}^2$.}
  
Note that in (\ref{intro_2a}-\ref{intro_2b}), the initial conditions have been chosen to be uniformly zero. In general, the particular form of the initial condition will play a large role in determining the singularity set of the system. \al{In the present work, however,} we restrict our attention to the uniform initial condition case to specifically study how the geometry of $\Omega$ and the dynamics of the underlying PDE alone combine to determine the singularity set. 

\al{Heuristically, the dynamics of \alt{problems} \eqref{intro_2} can be decomposed into two temporal regimes. The first is a short time regime whereby the solution is shaped by the initial data, the geometry of the region and the dynamics of the governing equation. The second is the blow up regime $(t-T_{\eps})\ll1$, in which the solution is changing very rapidly in a highly localized vicinity of the singularity points. Whichever point(s) enter the basin of attraction of the blow-up regime first, will eventually be the blow-up locations of the system. This mechanism promotes single point blow-up which is the standard behavior observed for \alt{problems} \eqref{intro_2} with general initial data - consequently simultaneous multiple point blow-up is not generically stable to asymmetric perturbations in the initial condition.  This instability is discussed further in the disc geometry example of \S\ref{sec:ex_disc} and the square example of \S\ref{sec:ex_square}.  } 
 
\section{One dimensional Theory}\label{sec:1d}

In this section, we develop a leading order, short time description of solutions to \alt{problems} \eqref{intro_2} in the one dimension strip $\Omega=[-1,1]$, by means of matched asymptotic expansions. A similar expansion, including higher order corrections, was obtained in \cite{AL2012} for \alt{problem} \eqref{intro_2b} with $f(u) = 1/(1-u)^2$. The higher order correction terms in the asymptotic expansion allow for a more accurate quantitative description of the solution. However, as will become apparent, the multiple singularities shown in Fig.~\ref{fig:intro1} are fully explained by considerations at leading order. In addition, the leading order analysis herein, is valid for the more general class of non-linearities described by \eqref{intro_2c}.

Once a uniformly valid asymptotic expansion is established, its validity is confirmed with comparison to numerical simulations and found to be good, even for moderately large times. In the second order case \eqref{intro_2a}, we find that the global maximum of the uniformly valid solution is $x=0$ for all $\eps$. This indicates, as shown in \cite{AFMC85}, that this point is selected by the dynamics of the full PDE for singularity. In contrast, when this asymptotic theory is applied to the fourth order case \eqref{intro_2b}, the location and multiplicity of the global maxima are found to depend on the value of $\eps$.

\subsection{Laplacian Case}\label{sec1d:lap}

In this section, a small amplitude asymptotic solution to the PDE
\begin{equation}\label{eqn:1d:1}
u_t = \eps^2 u_{xx} + f(u), \quad -1<x<1, \quad 0<t<\al{T_{\eps}}; \qquad u(\pm 1,t) = 0; \qquad u(x,0) = 0
\end{equation}
is developed. In the outer region, away from $x=\pm1$, the solution is spatially uniform and satisfies
\begin{equation}\label{eqn:1d:2}
\frac{d u_0}{dt} = f(u_0), \quad 0< t < \al{T_{0}}; \qquad u_0(0) = 0.
\end{equation}
Boundary layers are required in the vicinity of $x=\pm1$ to enforce the conditions $u(\pm 1)=0$. The formulation for the boundary at $x=1$ is established through the stretching variables
\begin{equation}\label{eqn:1d:3}
u(x,t) = u_0(t) \al{[ v(\eta) + \mathcal{O}(\phi)]}, \qquad \eta = \frac{1-x}{\phi(t;\eps)}, \qquad \phi(t;\eps) = \eps \, u_0(t)^{1/2}
\end{equation}
Substituting variables \eqref{eqn:1d:3} into \eqref{eqn:1d:1} \al{yields the equation
\begin{equation}\label{eqn:1d:extra1}
f(u_0) \left[ v - \frac{\eta}{2} \frac{dv}{d\eta}  \right] = \frac{d^2v}{d\eta^2} + f(u_0v).
\end{equation}
This equation further simplifies from noticing that $u_0 = \mathcal{O}(t)$ as $t\to0$ and that $f(t) = 1 + \mathcal{O}(t)$ from the assumptions \eqref{intro_2c} on $f(t)$. Therefore, a short time and small amplitude approximation of \eqref{eqn:1d:extra1} satisfies}
\bsub\label{eqn:1d:4}
\begin{equation}\label{eqn:1d:4a}
\frac{d^2v}{d\eta^2} + \frac{\eta}{2} \frac{dv}{d\eta} - v = -1, \quad \eta >0; \qquad v(0) = 0; \qquad v(\eta) \to 1, \quad \eta\to\infty.
\end{equation}
The solution of the linear problem \eqref{eqn:1d:4a} and its far field asymptotic \al{behavior} are given by
\begin{equation}\label{eqn:1d:4b}
\begin{array}{rcl}
v ( \eta ) &=& 1 - e^{-\frac{\eta^2}{4}} \left[ -\frac{\eta}{\sqrt{\pi}} +  \Big( 1+ \frac{\eta^2}{2} \Big) \, e^{\frac{\eta^2}{4}}\mbox{\al{erfc}} \, \big(\frac{\eta}{2}\big) \right], \\[10pt]
v ( \eta ) & \sim & 1- \frac{8}{\sqrt{\pi} \eta^3} \,  e^{-\frac{\eta^2}{4}} \left[ 1 + \mathcal{O}\Big( \frac{1}{\eta^2} \Big) \right], \qquad \eta \to \infty,
\end{array}
\end{equation}
\esub
where \al{erfc$(z)$} is the complementary error function. \al{The solution of \eqref{eqn:1d:4}, displayed in Fig.~\ref{fig:two profiles_a}, represents a short-time similarity solution of the boundary layer equation.} By superposing the contributions from each boundary and the uniform region, followed by a subtraction of overlapping terms, the uniformly valid short time asymptotic solution
\begin{equation}\label{eqn:1d:5}
u(x,t) \sim u_0(t)  +  u_0(t) \left[ v\left(\frac{1-x}{\phi(t;\eps)}\right) + v\left(\frac{1+ x}{\phi(t;\eps)}\right) -2 \right],  \qquad \phi(t;\eps) = \eps \, u_0(t)^{1/2},
\end{equation}
is established. The asymptotic \al{behavior} \eqref{eqn:1d:4b} of the profile $v(\eta)$ is increasing monotonically to the limiting value which indicates that the solution will have larger value in regions farther away from the boundary. Therefore, the regime of interest in the global approximation \eqref{eqn:1d:5} is when $1\pm x \gg \phi$. After applying the far field \al{behavior} \eqref{eqn:1d:4b} relevant to this regime, the solution with exponential corrections, valid away from $x=\pm1$,
\begin{equation}\label{eqn:1d:6}
u(x,t) \sim u_0(t) \left[ 1 -  \frac{8 \phi^3}{\sqrt{\pi}} \left[   \frac{1}{(1-x)^3} \exp\left[ -\frac{\al{(1-x)^2}}{4\phi^2} \right] + \frac{1}{(1+x)^3} \exp\left[ -\frac{\al{(1+x)^2}}{4\phi^2} \right]    \right] \right],
\end{equation}
is obtained. \al{As $\lim_{t\to T_0}\phi(t;\eps)=\infty$, this expression is defined for $0<t<\al{T_{0}}$, however, we can only expect good quantitative validity for $t\ll1$. Indeed, from \eqref{eqn:1d:3}, we have that $u = u_0(t) [ v(\eta) + \mathcal{O}(u_0)]$, and calculate that the relative $L^{\infty}$ error estimate for the asymptotic approximation is
\begin{equation}\label{eqn:1d:error}
\tau_{\mbox{rel}} = \frac{\| u - u_0 v \|_{\infty}}{\|u\|_{\infty}} = \mathcal{O}(u_0) =  \mathcal{O}(t), \quad \mbox{as} \quad t\to0.
\end{equation}}
We now posit that the global maximum \al{$x=0$} of \eqref{eqn:1d:6} provides a predictor of the location of the singularity of the full problem \eqref{eqn:1d:1}. Heuristically, this maximum of the small amplitude solution will enter the basin of attraction of a stable similarity solution regime of \eqref{eqn:1d:1} before any others.
\subsubsection{Example: Power Nonlinearity}\label{ex:1d:pow}

\begin{figure}[H]
\centering
\subfigure[$t = 0.4$]{\includegraphics[width = 0.45\textwidth]{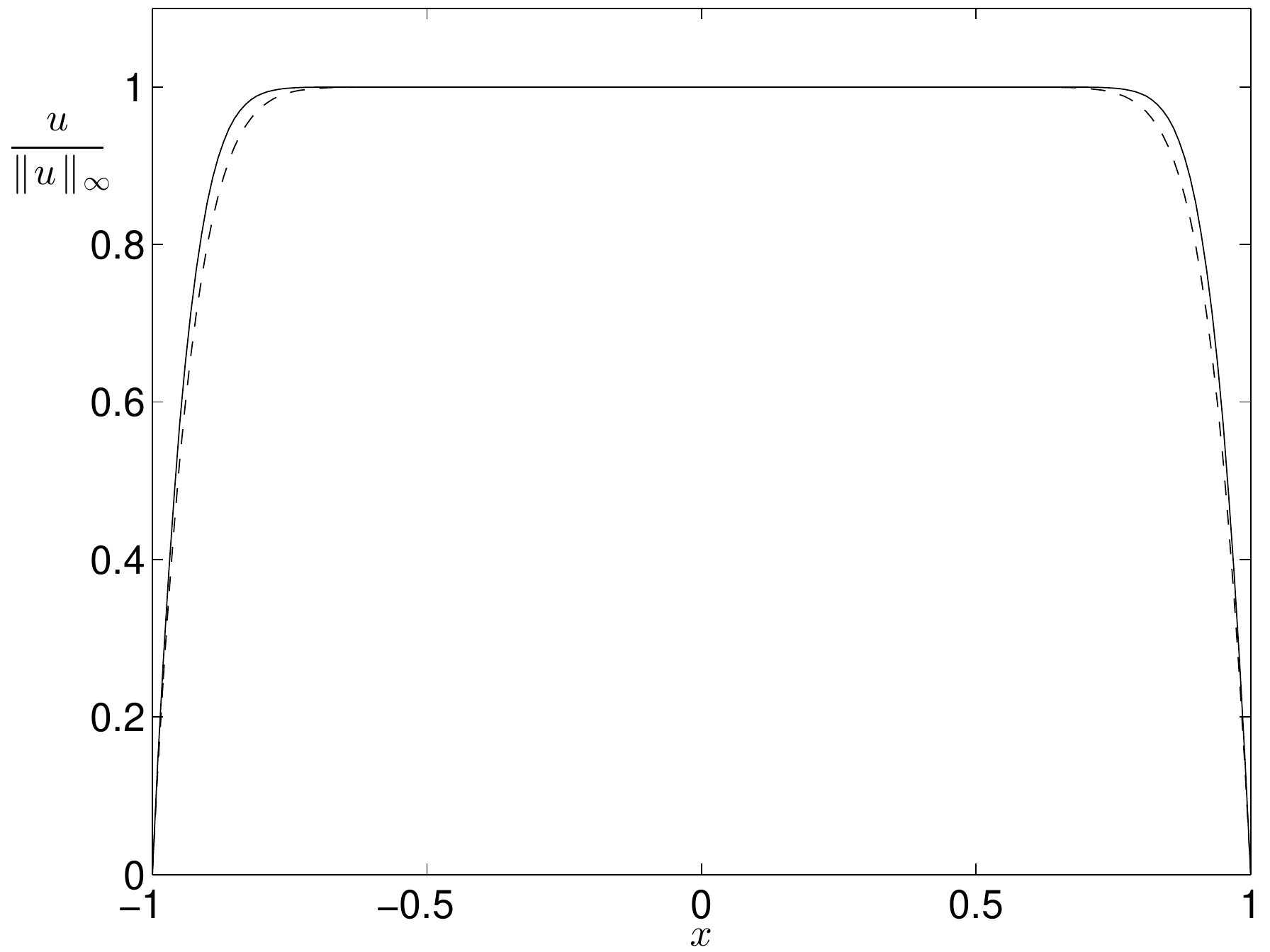}\label{fig:ex_pow_1Da}}
\qquad
\subfigure[Relative Error]{\includegraphics[width = 0.45\textwidth]{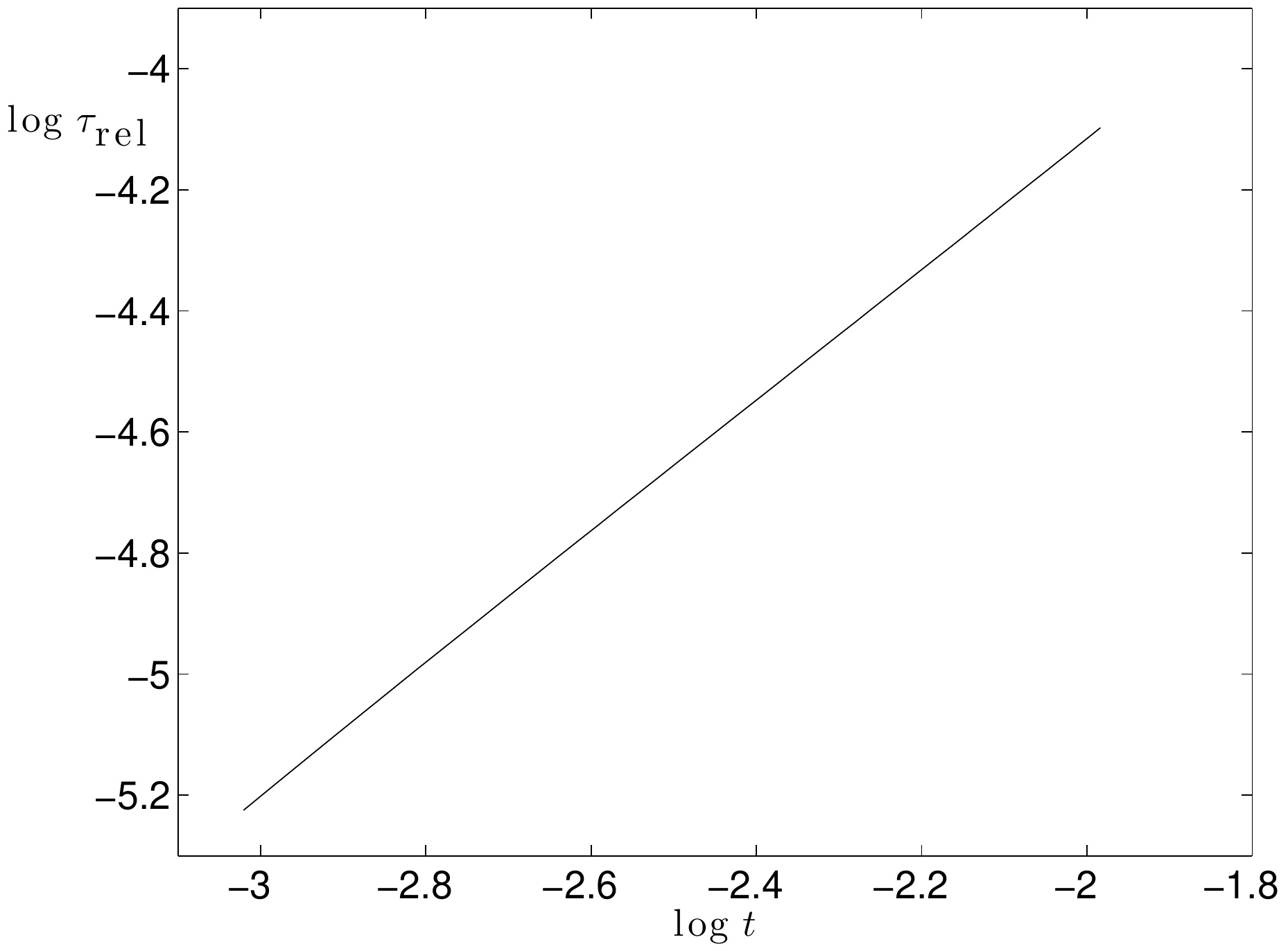}\label{fig:ex_pow_1Db}}
\parbox{5in}{\caption{Comparison of full numerical (solid) solution of \eqref{ex1:1d:pow1} with $\eps=0.1$, together with the asymptotic (dashed) solution \eqref{eqn:1d:6}. The normalized solution \al{with $\|u\|_{\infty} \simeq 0.446$} is displayed and good agreement is observed, even up to the moderately large value $t=0.4$. \al{Right panel shows a log-log plot of the relative $L^{\infty}$ error \eqref{eqn:1d:error} against time and exhibits the predicted $\mathcal{O}(t)$ behavior as $t\to0$}. \label{fig:ex_pow_1D}}}
\end{figure}

To demonstrate the efficacy of the short time asymptotic solution \eqref{eqn:1d:6}, consider the specific choice $f(u) = (1+u)^2$ and the problem
\begin{equation}\label{ex1:1d:pow1}
\begin{array}{c}
u_t = \eps^2 u_{xx} + (1+u)^2, \quad -1<x<1, \quad 0<t<\al{T_{\eps}}; \\[5pt]
u(\pm 1,t) = 0; \qquad u(x,0) = 0.
\end{array}
\end{equation}
For this case, the purely reaction problem $u_{0t} = (1+u_0)^2$, $u_0(0)=0$ has the solution $u_0(t) = t/(1-t)$ which blows up at $\al{T_{0}}=1$. Note that since $u_0(t)$ is a supersolution of \eqref{ex1:1d:pow1}, standard comparison principles imply $u_0(t) >u(x,t)$ so that $1=\al{T_{0}}<\al{T_{\eps}}$. For more details on upper and lower bounds for $\al{T_{\eps}}$, see \cite{BE}. In Fig.~\ref{fig:ex_pow_1D}, a comparison between the numerical solution of \eqref{ex1:1d:pow1} and the asymptotic solution \eqref{eqn:1d:6} is displayed for $\eps = 0.1$. In Fig.~\ref{fig:ex_pow_1Da} excellent agreement is observed for even moderately large values of $t$. To obtain accurate numerical solutions of \eqref{ex1:1d:pow1} very close to singularity, r-adaptive moving mesh methods are employed together with computational time stepping. For more details, see \cite{BJFW}.

\subsection{Bi-Laplacian Case}\label{sec:1d:biharm}

\begin{figure}[H]
\centering
\subfigure[Laplacian Case.]{\includegraphics[width=0.45\textwidth]{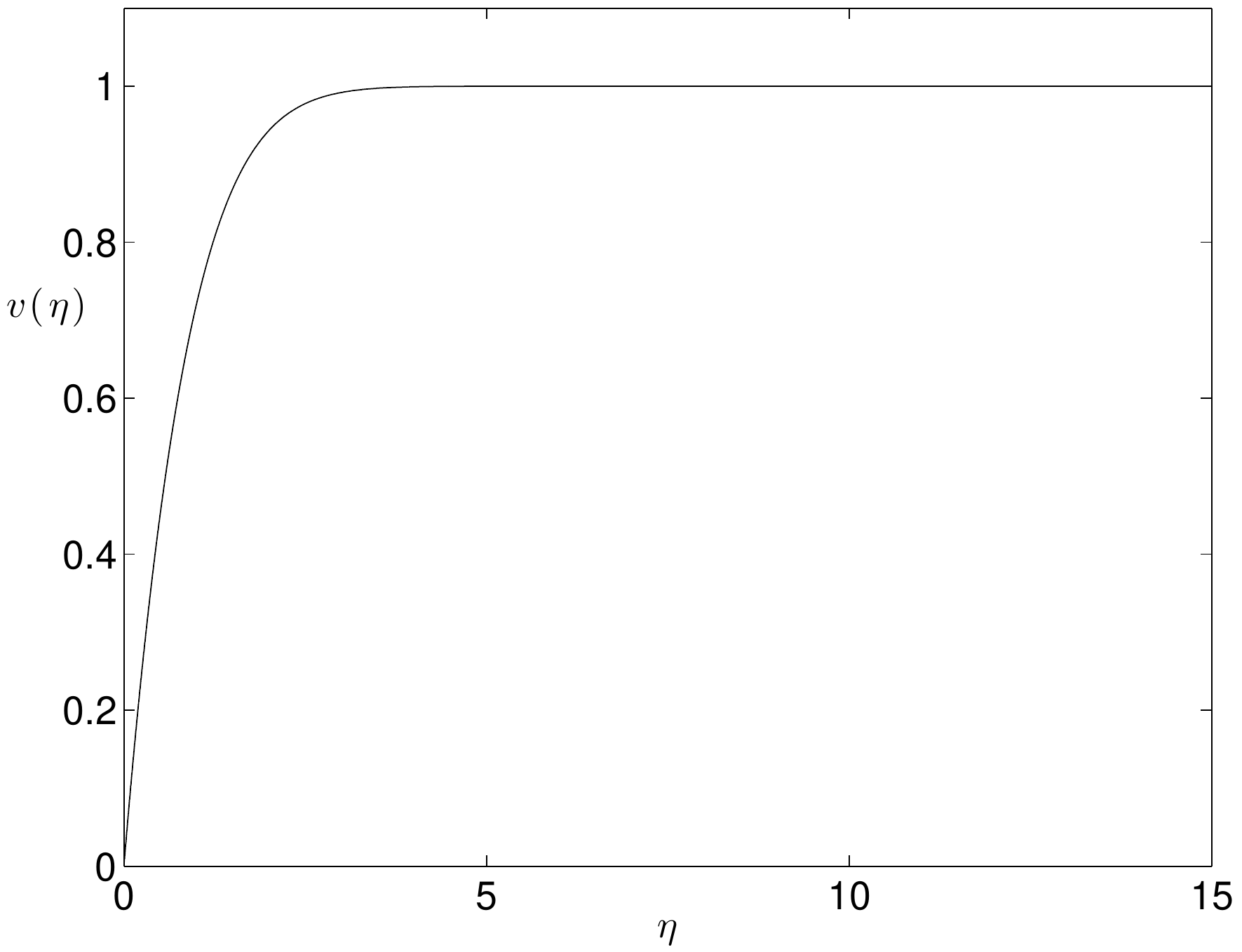}\label{fig:two profiles_a}}
\qquad
\subfigure[Bi-Laplacian Case.]{\includegraphics[width=0.475\textwidth]{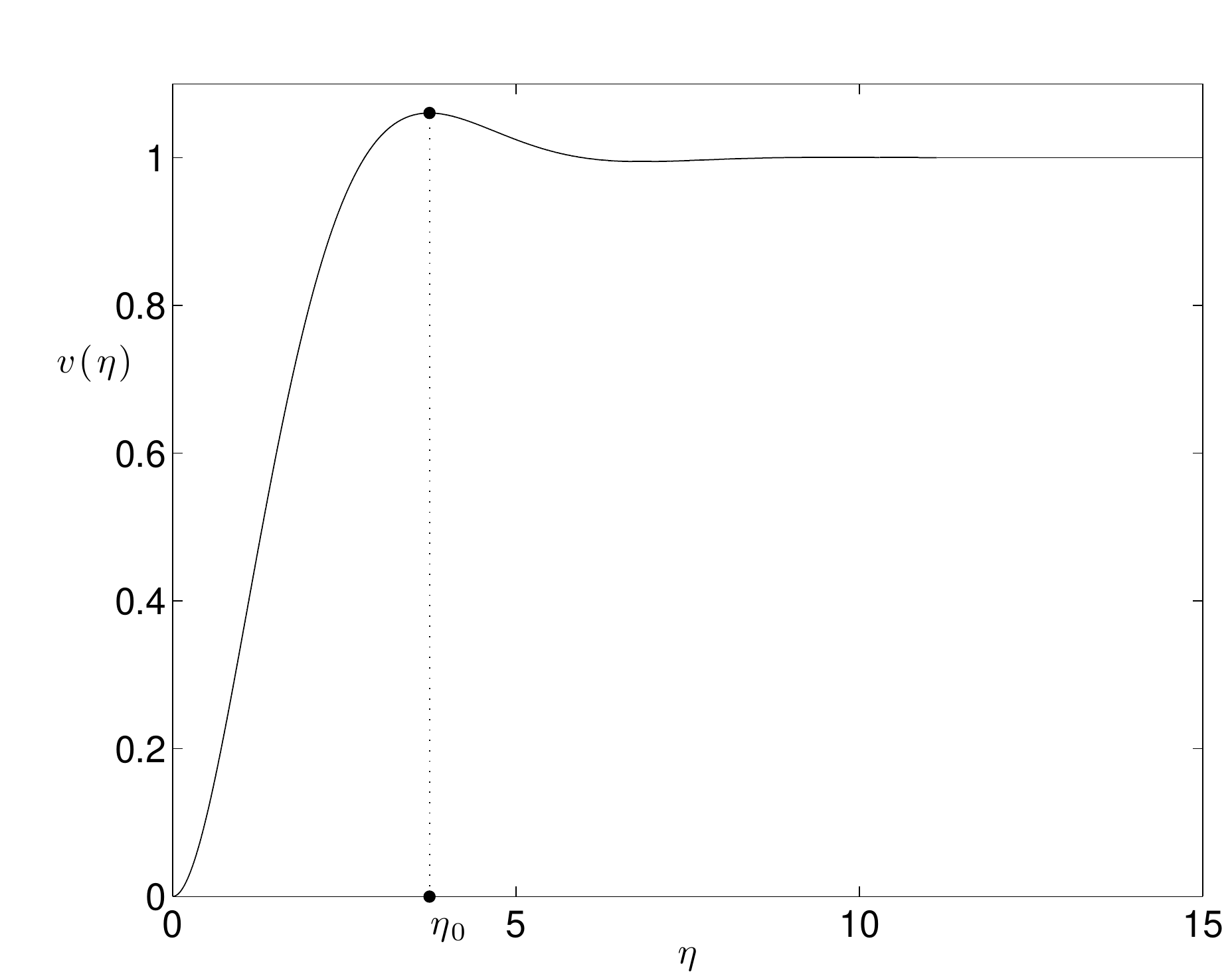}\label{fig:two profiles_b}}
\parbox{5in}{\caption{Panel (a) displays the solution of the second order problem \eqref{eqn:1d:4} while panel (b) displays the solution to the fourth order problem \eqref{eqn:1d:12}. The important distinction between the two cases is that the second order profile is monotone while the fourth order attains a global maximum at $\eta = \eta_0$.  \label{fig:two profiles}}}
\end{figure}

With the fourth order problem
\begin{equation}\label{eqn:1d:7}
\begin{array}{c}
u_t = -\eps^4 u_{xxxx} + f(u), \quad -1<x<1, \quad 0<t<\al{T_{\eps}}; \\[5pt]
u(\pm 1) = u_{x}(\pm 1) = 0; \qquad u(x,0) = 0, \quad -1<x<1,
\end{array}
\end{equation}
the spirit of the analysis is the same. The outer region is the solution to the ODE problem $u_{0t} = f(u_0)$ while the boundary conditions are enforced in layers at $x=\pm1$. To establish a boundary layer in the vicinity of $x=1$, the stretching variables
\begin{equation}\label{eqn:1d:8}
u(x,t) = u_0(t) \al{[v(\eta) +\mathcal{O}(u_0)]}, \qquad \eta = \frac{1-x}{\phi(t;\eps)}, \qquad \phi(t;\eps) = \eps \, u_0(t)^{1/4}
\end{equation}
are introduced. The leading order term in the expansion as $\phi\to0$, is the profile $v(\eta)$ satisfying
\begin{equation}\label{eqn:1d:9}
-\frac{d^4v}{d\eta^4} + \frac{\eta}{4} \frac{dv}{d\eta} - v = -1, \quad \eta >0; \qquad v(0) =v'(0)=  0; \qquad v(\eta) \to 1, \quad \eta\to\infty.
\end{equation}
While a closed form solution to \eqref{eqn:1d:9} is available, it is rather unsightly and not very useful. The key feature is the far field \al{behavior}, which can be obtained from a WKB analysis. By applying the large $\eta$ anzatz \al{$v(\eta) = 1 +\exp [\psi]$}, $\psi$ is found to satisfy
\begin{equation}\label{eqn:1d:10}
-\left(\frac{d\psi}{d\eta}\right)^4  + \frac{\eta}{4} \frac{d\psi}{d\eta}  \sim 0,
\end{equation}
at leading order. This ODE admits three non-trivial solutions
\begin{equation}\label{eqn:1d:11}
\psi_{j} =  3 \, \eta^{4/3}\, 2^{-8/3} \, \exp\left[ \frac{2\pi i j}{3} \right], \qquad  j = 0,1,2,
\end{equation}
however, \al{only $\exp[\psi_{1}]$ and $\exp[\psi_{2}]$ decay as $\eta\to\infty$}. Therefore the complete specification of \eqref{eqn:1d:9} with the obtained far field \al{behavior} is
\bsub\label{eqn:1d:12}
\begin{align}
\label{eqn:1d:12a} -\frac{d^4v}{d\eta^4} + \frac{\eta}{4} \frac{dv}{d\eta} - v &= -1, \quad \eta >0; \qquad v(0) =v'(0)=  0; \\[5pt]
\label{eqn:1d:12b} v(\eta) &= 1 + A \sin\left[ \sqrt{3} \omega \, \eta^{4/3} + \theta\right] e^{-\omega \eta^{4/3}} \Big[ 1 + {\bf o}(1) \Big], \qquad \eta\to\infty,
\end{align}
\esub
where $A$, $\theta$ are arbitrary constants and $\omega = 3\cdot 2^{-11/3}$. It is crucial to notice that the far field \al{behavior} of the equivalent second order problem \eqref{eqn:1d:4b} is strictly increasing towards the limiting value, while the solution of the fourth order problem \eqref{eqn:1d:12b} exhibits oscillations about it (cf. Fig.~\ref{fig:two profiles}). This loss of monotonicity in the fourth order case has a dramatic effect on the location of the singularity selected by the dynamics of the full PDE.

The uniformly valid short time asymptotic solution to \eqref{eqn:1d:7} is now given by
\begin{equation}\label{eqn:1d:13}
u(x,t) \sim u_0(t)  +  u_0(t) \left[ v\left(\frac{1-x}{\phi(t;\eps)}\right) + v\left(\frac{1+ x}{\phi(t;\eps)}\right) -2 \right],  \qquad \phi(t;\eps) = \eps \, u_0(t)^{1/4}
\end{equation}
where $v(\eta)$ is the profile determined in \eqref{eqn:1d:12}. The aim now is to understand which value(s) $x_c(t)\in(-1,1)$ will be global maxima of expression \eqref{eqn:1d:13}. In the previously considered (cf. \S\ref{sec1d:lap}) second order case, the monotonicity of the boundary layer solution meant that exponentially small corrections in the outer region were responsible for determining critical points. While in the fourth order case the exponentially small corrections to $u(x,t)$ in the outer region are still present, they are subsidiary to the contributions to the global maxima at $1\pm x = \eta_0\phi$ where $\eta_0$ is the global maximum of the profile $v(\eta)$ solving \eqref{eqn:1d:12}. These maximum points are time dependent and move from the boundary to the interior of the strip. Heuristically, the multiple singularity phenomenon can now be explained by the transit of these peaks in the stretching boundary layer: if the two maxima have not \al{met} by \al{the} time of singularity, there will be two singularity points. If they have met by the time of singularity, the origin is the singularity point.

In addition, if we posit that the global maxima of \eqref{eqn:1d:13} act as surrogates for the singularity points selected by the dynamics of the full PDE, then a crude approximation for such points is
\begin{equation}\label{eqn:1d:14}
x_c(T_{\eps}) \sim \left\{ \begin{array}{cl} \pm(1- \eta_0 \phi_c  ), & \eta_0 \phi_c \leq 1; \\[5pt] 0, & \eta_0 \phi_c >1\end{array} \right. \qquad \phi_c(\eps) = \phi(\al{T_{\eps}};\eps),
\end{equation}
where $\al{T_{\eps}}$ is the finite time singularity of the full PDE \eqref{eqn:1d:7}. It is crucial to note here that for $\phi_c$ to be well-defined, the fact that $\al{T_{\eps}}<\al{T_{0}}$ for sufficiently small $\eps$  has been used in \eqref{eqn:1d:14}. In other words, the PDE problem \eqref{intro_2b} can blow up faster than the ODE problem $u_{0t} = f(u_0)$. 

In the second order case, the lower bound $\al{T_{0}}\leq \al{T_{\eps}}$ is a simple consequence of the maximum principle. In the fourth order problem \eqref{intro_1b}, this approach to bounding $\al{T_{\eps}}$ from below is not valid. The challenge of obtaining lower bounds to the blow-up time of \eqref{intro_1b} in $\Omega = \mathbb{R}^N$ via a classical comparison argument with $u_0(t)$ was discussed in \cite{GALAK2002, Galak2001}. The authors' solution was to instead compare with the so-called \emph{order-preserving majorizing equation} for which a comparison principle can be established.  \al{It is also known that the blow-up time of \eqref{intro_2b} satisfies 
\begin{equation}\label{FriedmanResult}
T_{\eps} \leq T_0 + C\eps^{4/3},
\end{equation}
for constant $C>0$ \cite{FriedmanOswald88}. It would be interesting to extend these results and ideas to rigorously establish a lower bound on $T_{\eps}$ and therefore show that for \alt{problem} \eqref{intro_2b} on bounded regions, $\al{T_{\eps}}<\al{T_{0}}$ for certain ranges of $\eps$.} 


\subsubsection{Example: Exponential Nonlinearity}\label{sec:ex_biharm_exp}

\begin{figure}[H]
\centering
\subfigure[$t=0.5$]{\includegraphics[width=0.46\textwidth]{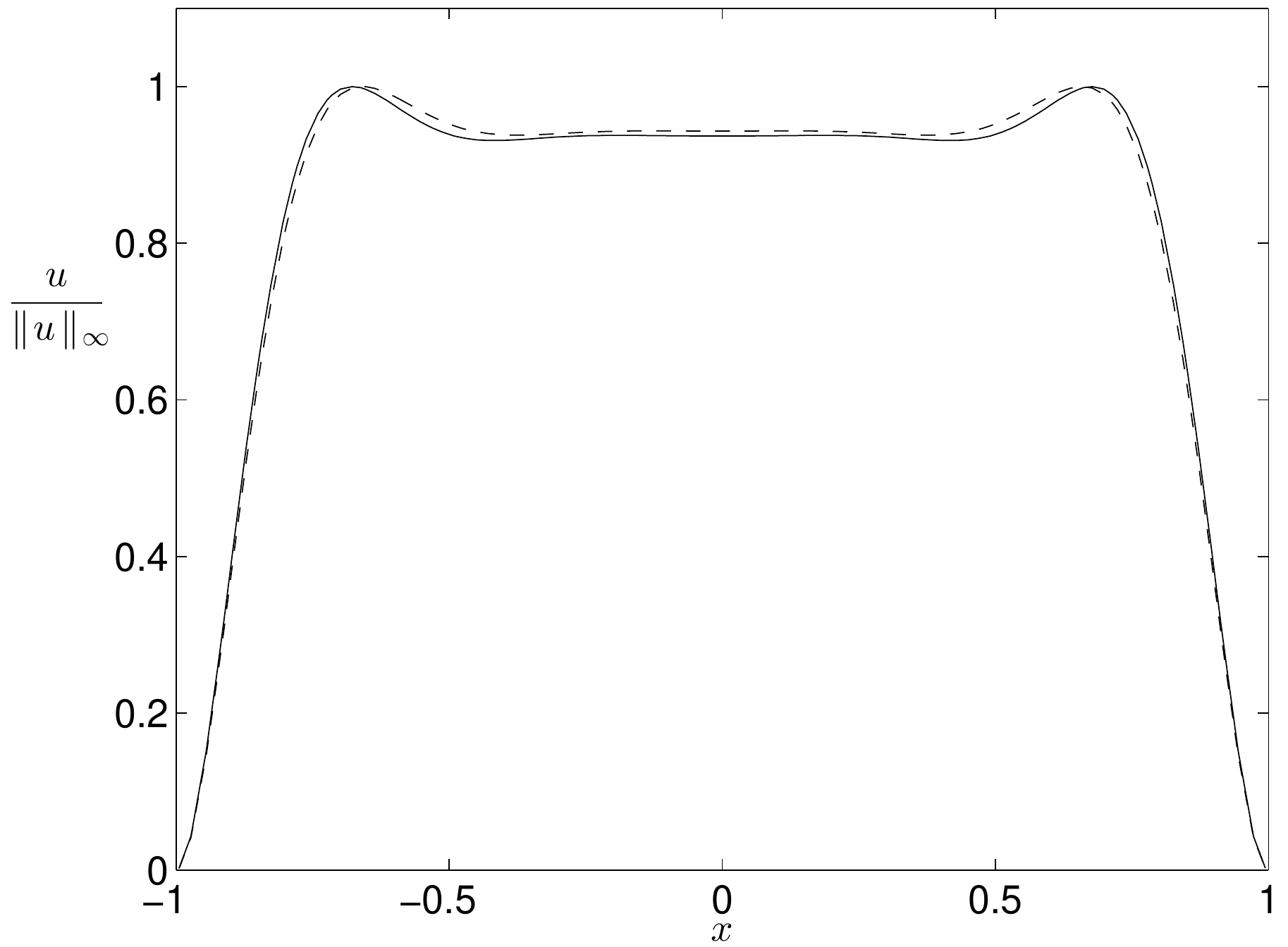}\label{fig:blowupexp:a} }
\hspace{0.05\textwidth}
\subfigure[Relative Error]{\includegraphics[width=0.45\textwidth]{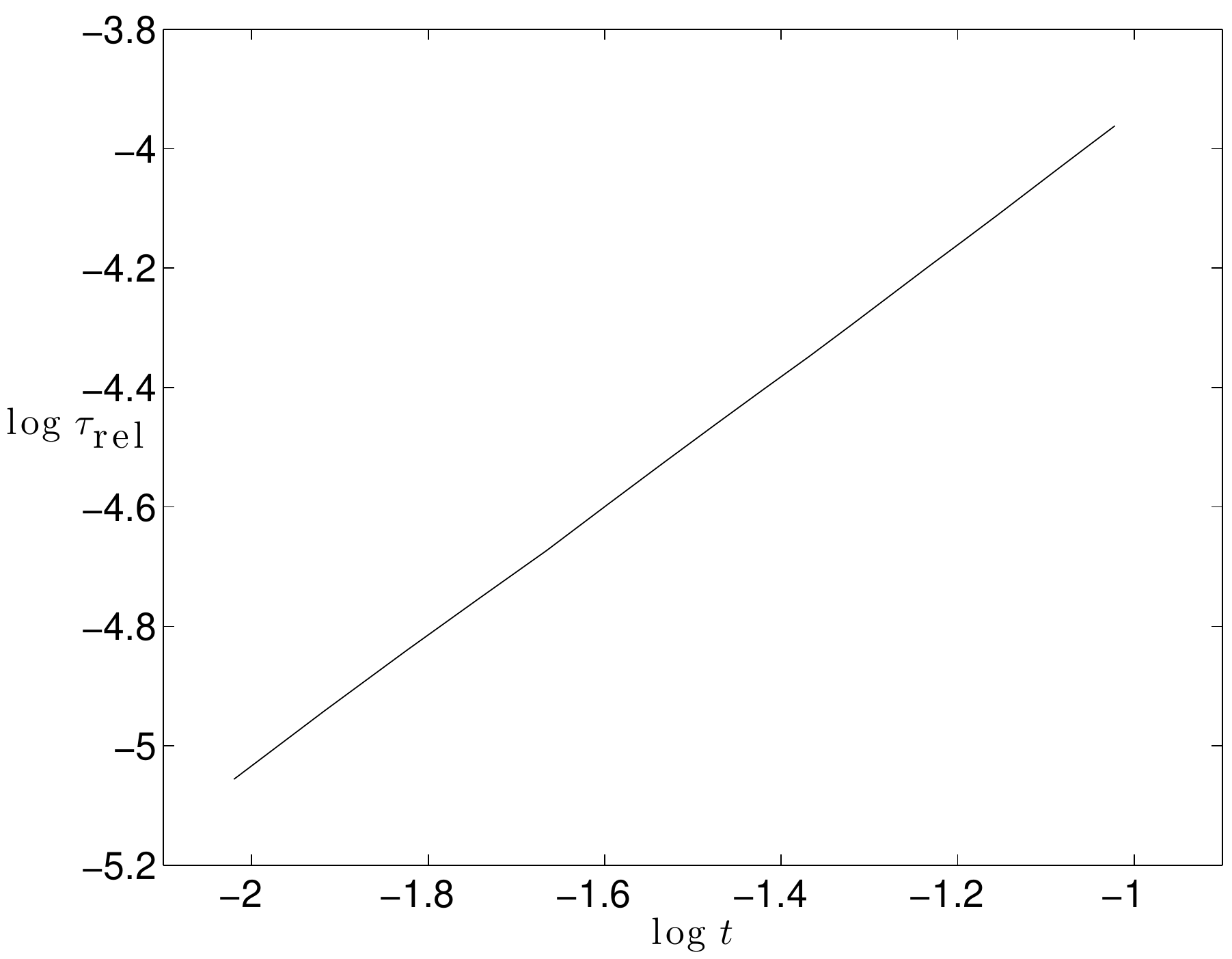} \label{fig:blowupexp:c} }
\parbox{5in}{\caption{Comparison of numerical (solid line) simulations of \eqref{eqn:ex:biharm:exp:1} and asymptotic (dashed line) predictions for $\eps=0.1$. In panels (a), we see that the asymptotic solution \eqref{eqn:1d:13} provides good agreement to the numerical solution, even for moderately large values of $t$. \al{The numerical solution shown here is for $\|u\|_{\infty}\simeq0.74$ and the blow up time is numerically estimated to be $T_{0.1} \simeq 0.9779$. Panel (b) displays a log-log plot of the relative $L^{\infty}$ error \eqref{eqn:1d:error} against time confirming that $\tau_{\mbox{rel}} = \mathcal{O}(t)$ as $t\to0$.} \label{fig:blowupexp}}}
\end{figure}

\begin{figure}[H]
\centering
\subfigure[Peak trajectory]{\includegraphics[width=0.45\textwidth]{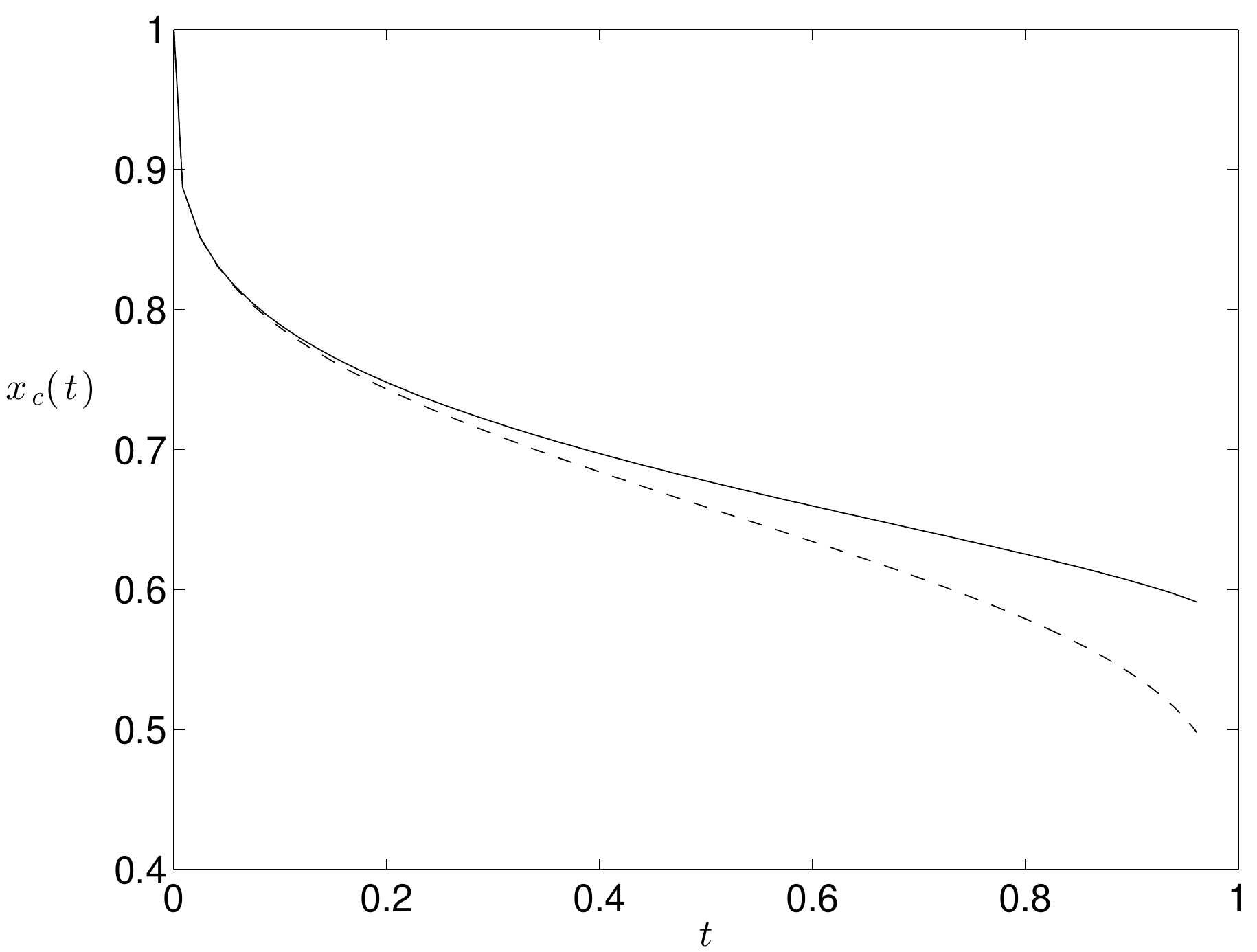}\label{fig:blowupexp:d}}
\qquad
\subfigure[Blow up points]{\includegraphics[width=0.45\textwidth]{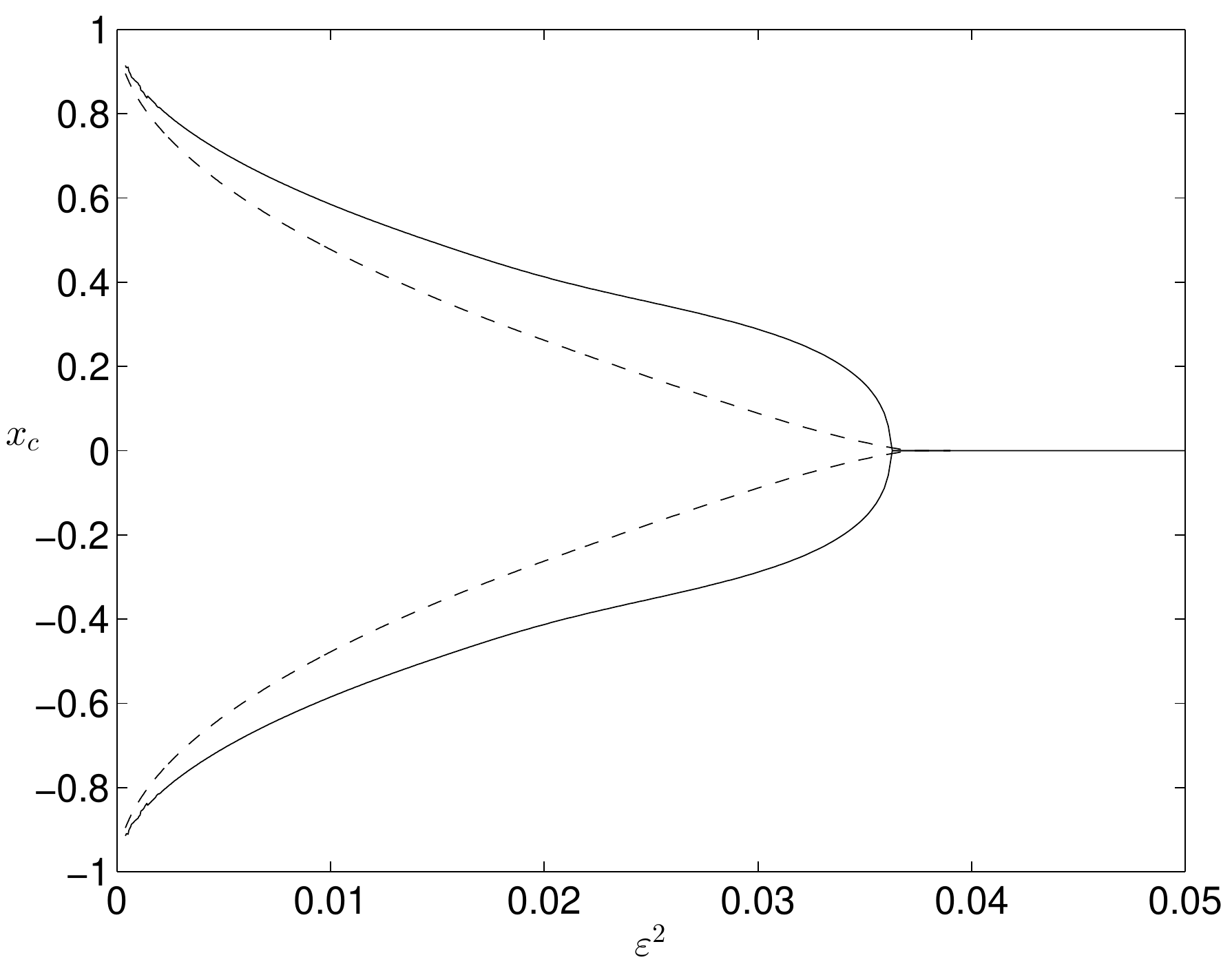}\label{fig:blowupexp:b}}
\parbox{5in}{\caption{\al{Panel (a) shows the peak trajectory from numerical simulations (solid line) and the asymptotic (dashed line) formulae \eqref{eqn:1d:13} for $\eps=0.1$. In panel (b), a comparison between the singularity points as predicted by the asymptotic formula \eqref{eqn:1d:14} and full numerical simulations, is displayed.}  \label{fig_biharm_exp}}}

\end{figure}

For a quantitative demonstration of the theory of the previous section, let us consider the case of the exponential nonlinearity $f(u)= e^u$. For this choice, the ODE problem $u_{0t} = e^{u_0}$, $u_0(0) = 0$ has the exact solution $u_0(t) = -\log(1-t)$. To simulate the PDE
\begin{equation}\label{eqn:ex:biharm:exp:1}
\begin{array}{c}
u_t = -\eps^4 u_{xxxx} + e^u, \quad -1<x<1, \quad 0<t<\al{T_{\eps}}; \\[5pt]
u(\pm 1) = u_{x}(\pm 1) = 0; \qquad u(x,0) = 0, \quad -1<x<1,
\end{array}
\end{equation}
very close to the singularity time $\al{T_{\eps}}$, the r-adaptive techniques described in \cite{BJFW,RXW,AL2012} are employed. As shown is Fig.~\ref{fig:blowupexp:a}, the leading order asymptotic prediction of \eqref{eqn:1d:13} is accurate, even for moderately large values of time. \al{The relative error of the approximation satisfies
\begin{equation}
\tau_{\mbox{rel}} = \frac{\|u-u_0 v\|_{\infty}}{\|u\|_{\infty}} = \mathcal{O}(u_0) = \mathcal{O}(t), \quad \mbox{as}\quad t\to 0,
\end{equation}
which is observed in Fig.~\ref{fig:blowupexp:c}}.
 \al{It is interesting to note here that $T_{\eps} < T_0 = 1$ for the value $\eps=0.1$ used here.} In Fig.~\ref{fig:blowupexp:b}, the asymptotic approximation is seen to provide a rather crude quantitative prediction of the singularity point \al{$x_c(T_\eps)$. In Fig.~\ref{fig:blowupexp:d}, the peak trajectory predicted by the asymptotic solution \eqref{eqn:1d:13} is plotted along with the numerically obtained peak trajectory for a fixed value of $\eps$. For $t\ll1$ the accuracy of the prediction is good, however, at later values of $t$ the asymptotic prediction overestimates the peak velocity.}  The accuracy of this prediction can be improved considerably by developing the asymptotic expansions beyond leading order. However, the main goal herein is to demonstrate the underlying principle behind the phenomenon, a purpose for which the leading order theory suffices.

\section{Two Dimensional Theory}\label{sec:2d}

The \al{generalization} of the previous analysis to bounded two-dimensional star-shaped regions is facilitated by implementing an arc-length tangent coordinate system $(\rho, s)$, where $\rho>0$ measures the distance from $x\in \Omega$ to $\partial\Omega$, whereas on $\partial\Omega$ the coordinate $s$ denotes arc-length along the boundary. In this co-ordinate system, the Laplacian operator admits the representation
\begin{equation}\label{eqn:2d:1}
\Delta = \frac{\partial^2}{\partial\rho^2} - \frac{\kappa}{1-\kappa\rho} \frac{\partial}{\partial\rho}
  + \frac{1}{1-\kappa\rho}\frac{\partial}{\partial s} \left( \frac{1}{1-\kappa\rho} \frac{\partial}{\partial s}
  \right)
\end{equation}
where $\kappa(s)$ indicates the curvature of $\partial\Omega$ as a function of arc-length $s$ along it. In this section, we analyze the second order \alt{problem}
\bsub\label{eqn:2d:2}
\begin{equation}\label{eqn:2d:2a}
 \begin{array}{lc} u_t = \eps^2\left[ \ds\frac{\partial^2}{\partial\rho^2} - \ds\frac{\kappa}{1-\kappa\rho} \ds\frac{\partial}{\partial\rho}
  + \ds\frac{1}{1-\kappa\rho}\frac{\partial}{\partial s} \left( \ds\frac{1}{1-\kappa\rho} \frac{\partial}{\partial s}
  \right) \right] u + f(u), & (x,t)\in\Omega_{\al{T_{\eps}}}; \\[5pt]
 u = 0, \quad& \al{(x,t)\in\partial\Omega_{T_{\eps}}};\\[5pt]
 u = 0, & \al{(x,t)\in\Omega_{0}},\end{array}
  \end{equation}
and its fourth order equivalent
  \begin{equation}\label{eqn:2d:2b}
  \begin{array}{c}
 \begin{array}{lc} u_t = -\eps^4\left[ \ds\frac{\partial^2}{\partial\rho^2} - \ds\frac{\kappa}{1-\kappa\rho} \ds\frac{\partial}{\partial\rho}
  + \ds\frac{1}{1-\kappa\rho}\frac{\partial}{\partial s} \left( \ds\frac{1}{1-\kappa\rho} \frac{\partial}{\partial s}
  \right) \right]^2 u + f(u), &  \hspace{-8pt}(x,t)\in\Omega_{T_{\eps}}; \\[5pt]
 u = 0, & \hspace{-9pt} (x,t)\in\partial\Omega_{T_{\eps}};\\[5pt]
 u = 0, &\hspace{-9pt}  (x,t)\in\Omega_{0},\end{array}
  \end{array}
\end{equation}
\esub
\al{where $\Omega_T$ and $\partial\Omega_T$ are defined in \eqref{timedomain}}. As previously, we posit that a uniform solution $u_0(t)$ satisfying
\begin{equation}\label{eqn:2d:3}
\frac{d u_0}{dt} = f(u_0), \quad 0< t < \al{T_{0}}; \qquad u_0(0) = 0,
\end{equation}
is valid away from $\partial\Omega$, while in the vicinity of the boundary, a local solution implements the boundary conditions. As in the 1D case, this local solution is monotone in the second order case and non-monotone in the fourth order problem. This distinction manifests itself in very different dynamical properties of the full PDEs \eqref{intro_2}. The accuracy of our analytical predictions are investigated on a sequence of test regions with finite element simulations of the PDEs \eqref{intro_2} carried out in MATLAB. Computational time stepping and local mesh refinement are applied in the vicinity of blow-up to ensure good resolution.  

\subsection{Laplacian Case}\label{sec:2d:lap}

The local solution in the vicinity of $\partial\Omega$ is established through the stretching variables
\begin{equation}\label{eqn:2d:4}
u(x,t) = u_0(t) \, v(\eta,s,t), \qquad \eta = \frac{\rho}{\phi}, \qquad \phi(t;\eps) = \eps\, u_0^{1/2}.
\end{equation}
and expansions in terms of $\phi\ll1$. After substituting \eqref{eqn:2d:4} and the expansion $v = v_0 + \phi v_1 + \cdots$ into \eqref{eqn:2d:2a}, we arrive at the problems
\bsub\label{eqn:2d:5}
\begin{align}
\label{eqn:2d:5a} \frac{d^2 v_{0}}{d\eta^2} &+ \frac{\eta}{2} \frac{dv_{0}}{d\eta} - v_{0} = -1, \quad \eta >0; \qquad v_{0}(0) = 0; \qquad v_{0}\to 1, \quad \eta\to\infty;\\[5pt]
\label{eqn:2d:5b} \frac{d^2 v_{1}}{d\eta^2} &+ \frac{\eta}{2} \frac{dv_{1}}{d\eta} - \frac{3}{2}v_{1} = \kappa \frac{dv_{0}}{d\eta}, \quad \eta >0; \qquad v_{1}(0) = 0; \qquad v_{1} \to 0, \quad \eta\to\infty.
\end{align}
\esub
The correction equation \eqref{eqn:2d:5b} also admits the decomposition $v_1(\eta,s) = \kappa(s) \bar{v}_1(\eta)$ where $\bar{v}_1(\eta)$ solves
\begin{equation}\label{eqn:2d:6}
\frac{d^2 \bar{v}_{1}}{d\eta^2} + \frac{\eta}{2} \frac{d\bar{v}_{1}}{d\eta} - \frac{3}{2}\bar{v}_{1} = \frac{dv_{0}}{d\eta}, \quad \eta >0; \qquad \bar{v}_{1}(0) = 0; \qquad \bar{v}_{1} \to 0, \quad \eta\to\infty.
\end{equation}
The leading order problem \eqref{eqn:2d:5a} depends only on the perpendicular distance from the boundary while the correction term \eqref{eqn:2d:5b} incorporates a dependence on the curvature of $\partial\Omega$. In fact equation \eqref{eqn:2d:5a} is precisely the one dimensional boundary profile established in \eqref{eqn:1d:4} and in particular, its far field asymptotic \al{behavior} is given by \eqref{eqn:1d:4b}.

We now turn to the problem of establishing a uniformly valid global solution at a general $x\in\Omega$. The previous analysis reveals that to leading order, the solution is principally determined by perpendicular distances from the boundary. Therefore to construct the solution at $x\in\Omega$, it is necessary to determine all boundary points $y\in\partial\Omega$ such that the straight line between $x$ and $y$ is contained in $\Omega$ and meets the boundary orthogonally at $y\in\partial\Omega$. With this in mind, suppose that for some $x\in\Omega$ there are $n$ boundary points $\{y_1,\ldots,y_n\}\in\partial\Omega$ such that the straight line $l(x,y_j)$ segment between $x$ and $y_j$ is contained in $\Omega$ and meets $\partial\Omega$ orthogonally at $y_j$, then a two term asymptotic solution $u(x,t)$ is given by
\begin{equation}\label{eqn:2d:7}
\begin{array}{c}
u(x,t) \sim u_0(t) +u_0(t)\ds\sum_{j=1}^{n} \left[ v_0 \left[ \ds\frac{|x-y_j|}{\phi(t;\eps)} \right]  + \phi\, \kappa(y_j) \bar{v}_1 \left[ \ds\frac{|x-y_j|}{\phi(t;\eps)} \right]-1  \right], \\[10pt]
\phi(t;\eps) = \eps \, u_0^{1/2}.
\end{array}
\end{equation}
In \eqref{eqn:2d:7}, $u_0(t)$ is the solution of the ODE \eqref{eqn:2d:3} and $\kappa(y_j)$ is the curvature of $\partial\Omega$ at $y_j\in\partial\Omega$. Equation \eqref{eqn:2d:7} represents a uniformly valid asymptotic expansion, derived in the limit $\phi\to 0$. Consequently, it is important to remark that the validity of this solution is restricted to regions $\Omega$ for which $\kappa = \mathcal{O}(1)$ as $\phi\to0$. In other words, the theory developed here is not valid for regions with rough or spiky boundaries $\partial\Omega$.
\begin{figure}[htbp]
\centering
\includegraphics[width=0.75\textwidth]{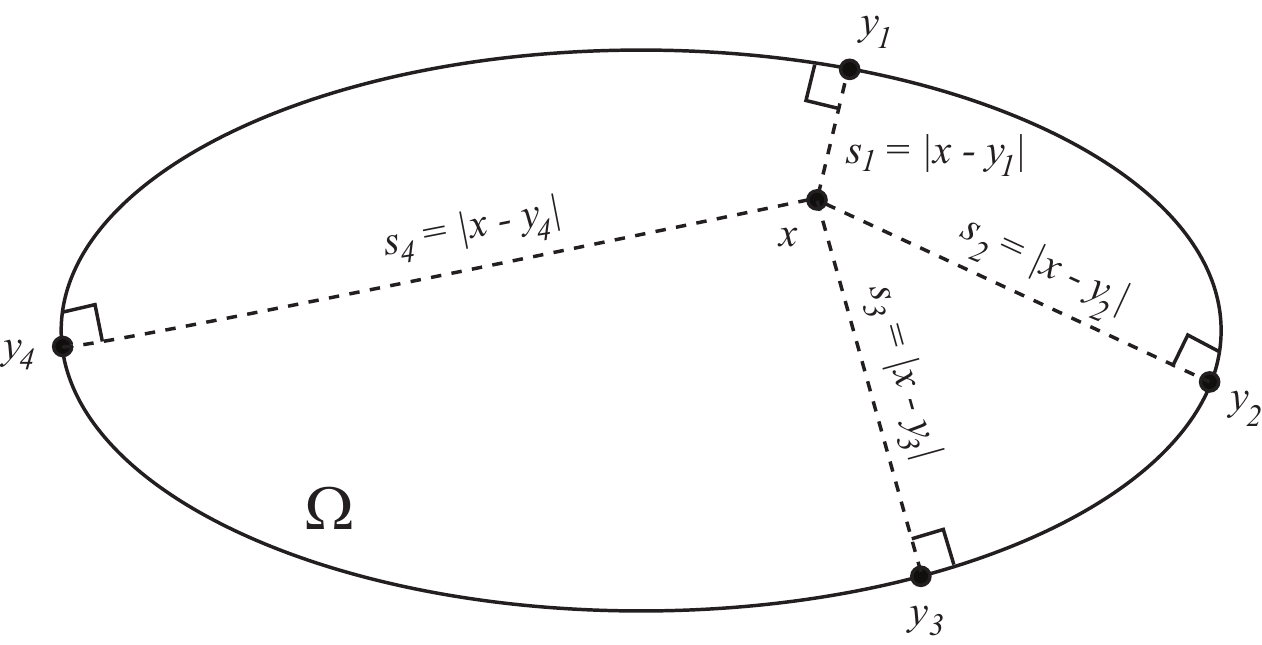}
\parbox{5in}{\caption{Schematic diagram illustrating the relationship between the points $\{y_j\}_{j=1}^n\in\partial\Omega$ and a point $x\in\Omega$ at which the uniform asymptotic solution is to be constructed. For this elliptical case and the chosen point $x\in\Omega$, there are $n=4$ boundary contributors. Figure reproduced from \cite{ALS2013}.\label{fig:rays} }}
\end{figure}
The objective now is to deduce the global maxima of the asymptotic formulation \eqref{eqn:2d:7} with the understanding that these critical points provide surrogates for those selected by the dynamics of the full PDE as the location of a finite time singularity. From the fact that the profile $v_0(\eta)$ is strictly increasing in $\eta$, the leading order theory predicts that the maximum will be away from $\partial\Omega$ and in the region where $|x-y_j|\gg\phi$. Therefore, \alt{in the case $\phi\ll1$, we are primarily interested in evaluating \eqref{eqn:2d:7} at points $x\in\Omega$ for which $v_0(|x-y_j|/\phi)$ and $v_1(|x-y_j|/\phi)$ can be replaced by their large argument \al{behavior} \eqref{eqn:1d:4b}. This reduces \eqref{eqn:2d:7}} to
\begin{equation}\label{eqn:2d:8}
u(x,t) \sim u_0(t) \left[1 - \frac{8 \phi^3}{\sqrt{\pi}}\sum_{j=1}^{n} \frac{1}{|x-y_j|^3} \exp\left[ -\frac{\al{|x-y_j|^2}}{4\phi^2} \right]  \right], \qquad \phi(t;\eps) = \eps \, u_0^{1/2},
\end{equation}
\alt{at points $x\in\Omega$ such that $|x-y_j|\gg\phi$. The form of \eqref{eqn:2d:8} can also be recovered by applying the anzatz
\[
u(x,t) = u_0(t) \left[1- \phi^3 A(x)\exp\left[-\frac{\Psi(x)}{4\phi^2} \right]  \right]
\]
to \eqref{intro_2a} and obtaining equations for $\Psi(x)$ and $A(x)$ in the limit $\phi\to0$}. The simplicity of \eqref{eqn:2d:8} allows several deductions to be made regarding the points $x\in\Omega$ at which $u(x,t)$ attains its maximum value and therefore the points which can be expected to reach singularity before others under the dynamics of the full PDE \eqref{eqn:1d:1}. As indicated by the form of \eqref{eqn:2d:8}, the distances $|x-y_j|$ play a large role in these determinations. As the correction terms to the uniform state are exponentially small, larger values of $|x-y_j|$ will have the smallest contribution to \eqref{eqn:2d:8} leaving the smaller values of $|x-y_j|$ to dominate the solution. Therefore, the maxima of \eqref{eqn:2d:8} will be located at points $x\in\Omega$ which are furthest from the boundary, \emph{i.e.},
\begin{equation}\label{eqn:2d:9}
x_c = \max_{x\in\Omega} \; \mathrm{d}(x,\partial\Omega).
\end{equation}
This result can be heuristically reconciled with the intuitive notion that, since the boundary condition $u=0$ on $\partial\Omega$ inhibits blow-up, and the solution is monotone away from the boundary, that points furthest from $\partial\Omega$ are those more likely to \al{develop a singularity as $t\to T_{\eps}^{-}$}.

In singularly perturbed elliptic problems, the distance function also plays a key role in the determining the location(s) of solution concentration. For example, it is well known (cf. \cite{Wei1,WeiDancer,WeiPino} and the references therein) that in the problem
\begin{equation}\label{eqn:2d:10}
 \left\{\begin{array}{c} \eps^2\Delta u - u + f(u) = 0, \quad x\in\Omega; \\[5pt]
 u>0, \quad x\in\Omega; \qquad u = 0,  \quad x\in\partial\Omega,\end{array} \right.
\end{equation}
for super-linear, sub-critical nonlinearity $f(u)$, spike solutions concentrate at \al{maxima of the distance function \eqref{eqn:2d:9}} as $\eps\to0$. In the following section, a similar analysis is applied to the fourth order problem \eqref{intro_2b}, however, the conclusions it provides are quite distinct.

\subsection{\alt{Bi-Laplacian} Case}\label{sec:2d:biharm}

To construct a short time asymptotic solution for \alt{problem} \eqref{eqn:2d:2b}, we again assume a uniform solution $u_0(t)$ satisfying \eqref{eqn:2d:3} in the interior coupled to a stretching boundary effect described by the variables
\begin{equation}\label{eqn:2d:biharm:1}
u(x,t) = u_0(t) \, v(\eta,s,t), \qquad \eta = \frac{\rho}{\phi}, \qquad \phi(t;\eps) = \eps\, u_0^{1/4}.
\end{equation}
After expanding \eqref{eqn:2d:biharm:1} with $v = v_0 + \phi v_1 + \cdots$ and collecting terms at the relevant order, $v_0$ and $v_1$ satisfy
\bsub\label{eqn:2d:biharm:2}
\begin{equation}\label{eqn:2d:biharm:2a}
\begin{array}{c}
 -\ds\frac{d^4v_{0}}{d\eta^4} + \ds\frac{\eta}{4} \frac{dv_{0}}{d\eta} - v_{0} = -1, \qquad \eta >0; \\[10pt]
  v_{0}(0) = v_{0\eta}(0) = 0; \qquad v_{0} \to 1, \qquad \eta\to\infty;
  \end{array}
\end{equation}
\begin{equation}\label{eqn:2d:biharm:2b}
\begin{array}{c}-\ds\frac{d^4v_{1}}{d\eta^4} + \ds\frac{\eta}{4} \frac{dv_{1}}{d\eta} - \frac{5}{4}v_{1} = -2\kappa \frac{d^3v_{0}}{d\eta^3}, \qquad \eta >0; \\[10pt]
 v_{1}(0) = v_{1\eta}(0) = 0; \qquad v_{1} \to 0, \qquad \eta\to\infty.
\end{array}
\end{equation}
Equation \eqref{eqn:2d:biharm:2b} for $v_1$ can be reduced by writing $v_1(\eta,s) = \kappa(s) \bar{v}_1(\eta)$ where
\begin{equation}\label{eqn:2d:biharm:2c}
\begin{array}{c}-\ds\frac{d^4\bar{v}_{1}}{d\eta^4} + \ds\frac{\eta}{4} \frac{d\bar{v}_{1}}{d\eta} - \frac{5}{4}\bar{v}_{1} = -2\frac{d^3v_{0}}{d\eta^3}, \qquad \eta >0; \\[10pt]
 \bar{v}_{1}(0) = \bar{v}_{1\eta}(0) = 0; \qquad \bar{v}_{1} \to 0, \qquad \eta\to\infty.
\end{array}
\end{equation}
\esub

As with the previous analysis in the Laplacian case, a uniformly valid asymptotic solution is constructed from the flat outer solution $u_0(t)$ and the boundary contributions established by (\ref{eqn:2d:biharm:1}-\ref{eqn:2d:biharm:2}). Therefore the short time asymptotic solution at $x\in\Omega$ is given by
\begin{equation}\label{eqn:2d:biharm:3}
\begin{array}{c}
u(x,t) \sim u_0(t) + u_0(t) \ds\sum_{j=1}^{n}\left[ v_0 \left[ \ds\frac{|x-y_j|}{\phi(t;\eps)} \right]  + \phi\, \kappa(y_j) \bar{v}_1 \left[ \ds\frac{|x-y_j|}{\phi(t;\eps)} \right]-1  \right], \\[10pt] 
\phi(t;\eps) = \eps \, u_0^{1/4},
\end{array}
\end{equation}
where $\{y_j\}_{j=1}^{n} \in \partial\Omega$ are the points for which the straight line $l(x,y_j)$ between $x$ and $y_j$ is contained in $\Omega$ and meets $\partial\Omega$ orthogonally at $y_j$. An illustration of how the points $\{y_j\}_{j=1}^n\in\partial\Omega$ depend on any particular $x\in\Omega$ is given in Fig.~\ref{fig:rays}. The goal now is to predict the location of blow up points by determining the global maxima at $t=\al{T_{\eps}}$ of \eqref{eqn:2d:biharm:3} for a range of bounded two-dimensional regions $\Omega$.

The explicit form of \eqref{eqn:2d:biharm:3} indicates that the quantities $s_j = |x-y_j|$ play a key role in determining whether or not a given point $x\in\Omega$ will be selected by the dynamics of the full PDE for singularity. As the profile $v_0(\eta)$ satisfying \eqref{eqn:2d:biharm:2a} has a global max at $\eta=\eta_0$, it follows from \eqref{eqn:2d:biharm:3} that points $x\in\Omega$ such that $|x-y_j| = \eta_0 \phi(t;\eps)$ are, to leading order, local maxima. Therefore, at time $t$ and for fixed $\eps$, we define the set $\omega(t)\subset\Omega$ where
\begin{equation}\label{eqn:2d:biharm:4}
\omega(t) = \left\{ x\in\Omega \; \vert \; \mathrm{d}(x,\partial\Omega) = \eta_0  \phi(t;\eps) \right\}.
\end{equation}
The set $\omega(\al{T_{\eps}})$ describes, to leading order in \eqref{eqn:2d:biharm:3}, points which are local maxima and are therefore more likely to be selected for singularity by the dynamics of the full PDE \eqref{eqn:2d:2b}. In general, singularities will not form simultaneously on all points of $\omega(\al{T_{\eps}})$. Instabilities along $\omega(\al{T_{\eps}})$ in combination with the higher order curvature effects present in \eqref{eqn:2d:biharm:3}, will increase the value of $u(x,t)$ at certain discrete points along $\omega(\al{T_{\eps}})$.

Amongst the points described by $\omega(\al{T_{\eps}})$, there may be special points $x\in\Omega$ which receive multiple boundary contributions, \emph{i.e.} there are $y_1,y_2\in\partial\Omega$ such that $\mathrm{d}(x,\partial\Omega) = \mathrm{d}(x,y_1) = \mathrm{d}(x,y_2)$. At such points, the leading order terms in the sum \eqref{eqn:2d:biharm:3} combine to increase the value of $u(x,t)$ therefore increasing the possibility that that point is selected for singularity. This motivates the definition of the \emph{skeleton} of the domain, $\mathcal{S}_{\Omega}$. The skeleton $\mathcal{S}_{\Omega}$ is the set of points $x\in\Omega$ for which there are at least two points $y_1,y_2\in\partial\Omega$ such that $\mathrm{d}(x,y_1)=\mathrm{d}(x,y_2)$ and such that the straight lines $l(x,y_1)$ and $l(x,y_2)$ are both contained in $\Omega$ and meet $\partial\Omega$ orthogonally. More compactly,
\begin{equation}\label{eqn:2d:biharm:5}
\mathcal{S}_{\Omega} = \left\{x \in \Omega \left\bracevert  \begin{array}{c} \exists y_1 \ne y_2 \in \partial \Omega, \quad \mathrm{d}(x,y_1) = \mathrm{d}(x,y_2),
\\[5pt]  l(x,y_1), \ l(x,y_2) \in\Omega, \\[5pt]  l(x,y_1)\perp \partial_{\tau} (y_1), \  l(x,y_2)\perp \partial_{\tau} (y_2). \end{array}\right.\right\},
\end{equation}
where $\partial_{\tau}(y)$ denotes the unit tangent vector to $\partial\Omega$ at $y\in\partial\Omega$.
\begin{figure}[htbp]
\centering
\subfigure[$\omega(t)$]{\includegraphics[width=0.35\textwidth,angle=-90]{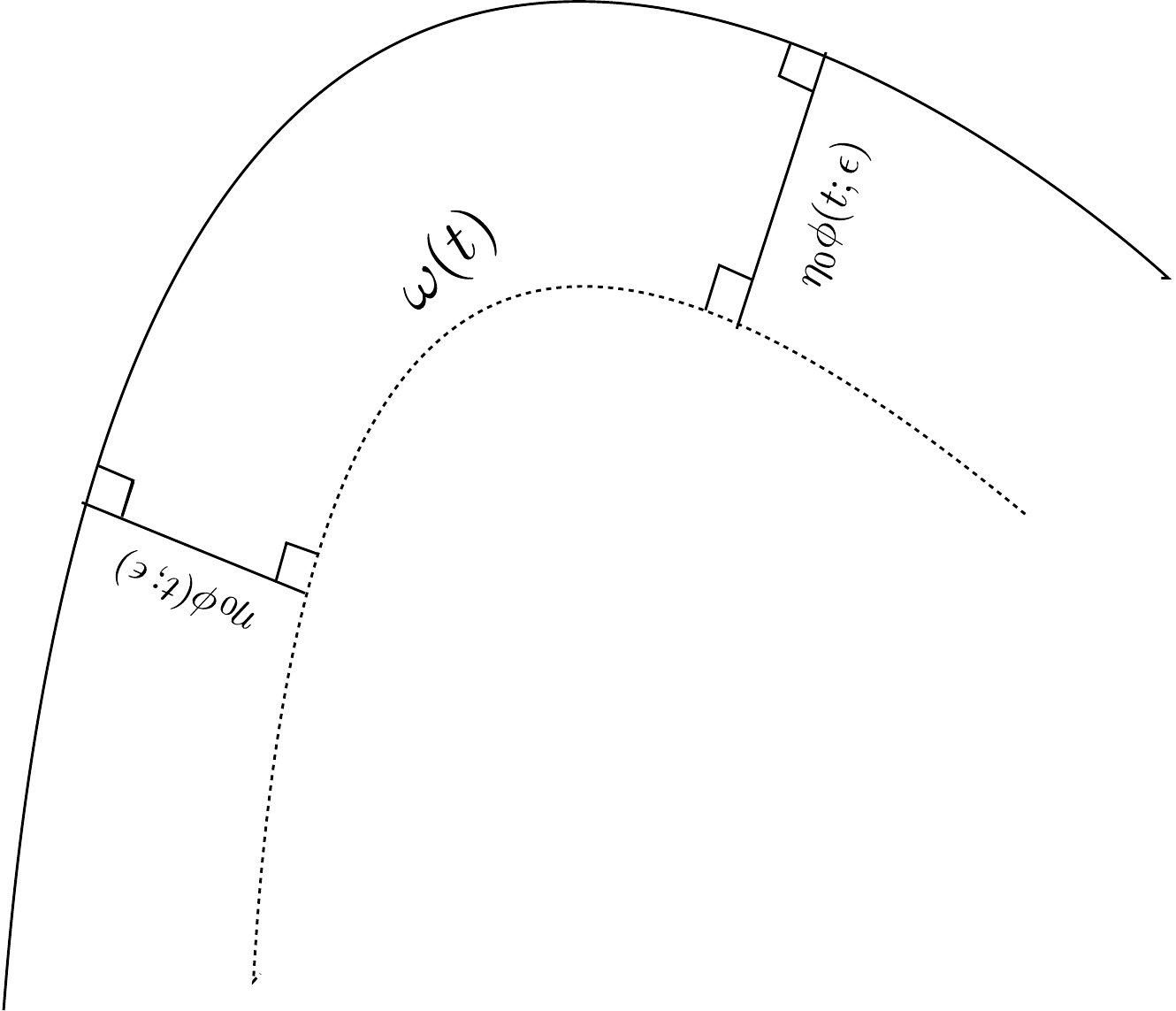}\label{fig:omega_skeleton_a}} \hspace{0.1\textwidth}
\subfigure[$x\in\mathcal{S}_{\Omega}$]{\includegraphics[width=0.35\textwidth,angle=-90]{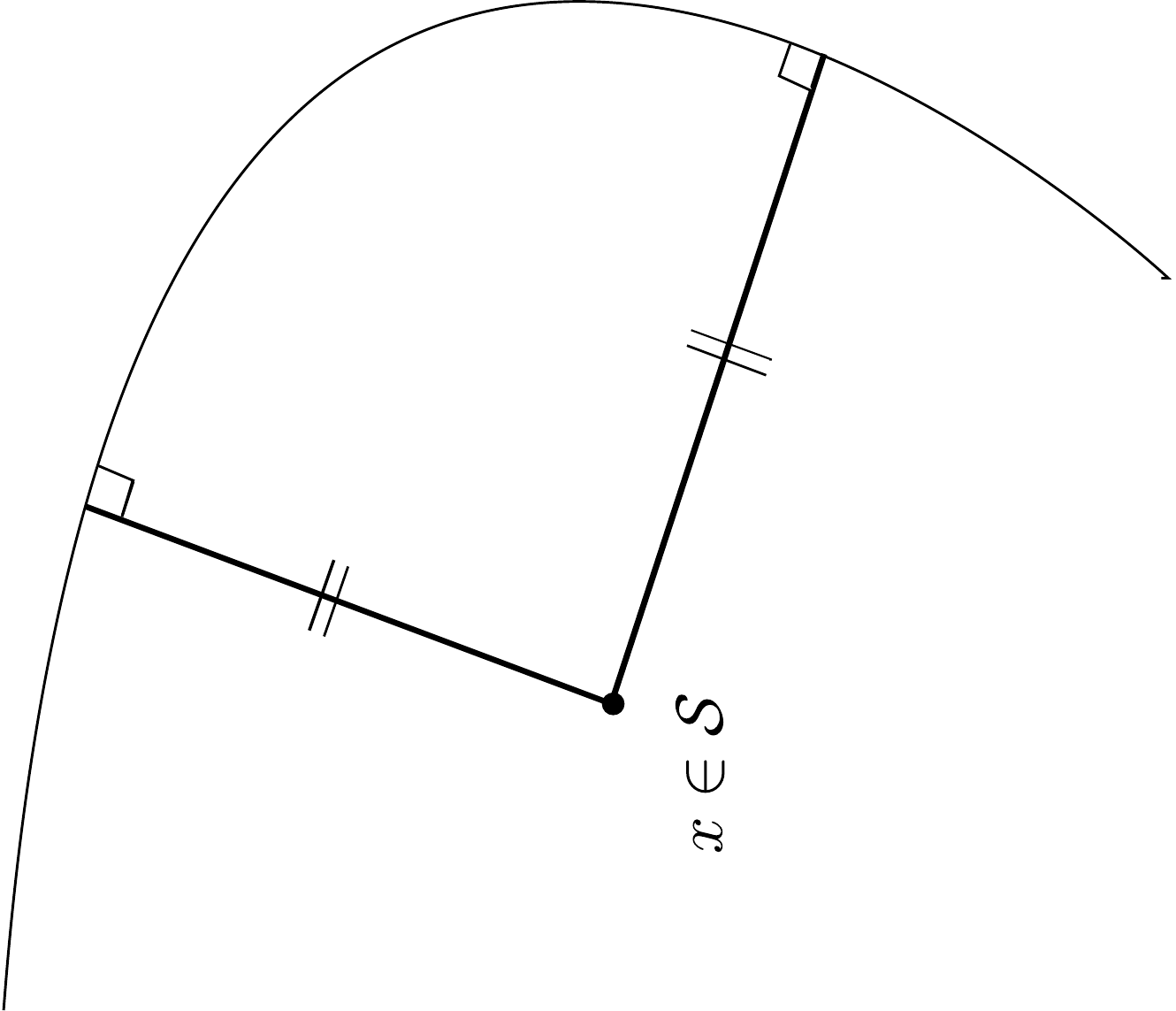}\label{fig:omega_skeleton_b}}
\parbox{5in}{\caption{The two structures $\omega(t)$ and $\mathcal{S}_{\Omega}$, defined in \eqref{eqn:2d:biharm:4} and \eqref{eqn:2d:biharm:5} respectively. Panel (a) The set $\omega(t)$ consists of the points in $\Omega$ that are at the distance $\eta_0 \, \phi(t; \eps)$ from $\partial \Omega$. Panel (b) The skeleton $\mathcal S_{\Omega}$ is the set of points in $\Omega$ that are equidistant to two or more points on $\partial \Omega$ such that the line segments between $x$ and those points are in $\Omega$ and meet $\partial\Omega$ orthogonally. \label{fig:omega_skeleton}}  }
\end{figure}
It is important to note that for a particular $x\in\mathcal{S}_{\Omega}$, there may be several sets of boundary points which mutually satisfy the conditions specified in \eqref{eqn:2d:biharm:5}. To make this notion concrete, consider the simple example of the elliptical region
\[ \Omega = \{ (x_1,x_2)\in\mathbb{R}^2 \ \left\bracevert  \ (x_1/a)^2 + (x_2/b)^2 = 1\} \right., \qquad 0<a<b.\]
For this particular region, $(0,0)\in\mathcal{S}_{\Omega}$ as both the sets of points $\{(-a,0), (a,0)\}$ and $\{(0,-b),(0,b)\}$ satisfy the conditions of \eqref{eqn:2d:biharm:5}. However, as will subsequently become apparent, it is necessary to distinguish between such sets of contributing points by their distances to $x\in\mathcal{S}_{\Omega}$. For this region, consider the function $s_{\Omega}: \mathcal{S}_{\Omega}\to\mathbb{R}$ such that
\begin{equation}\label{eqn:2d:biharm:6}
s_{\Omega}(x) = \min_{y_1\in\partial\Omega} \left\{ \mathrm{d}(x,y_1)  \left\bracevert  \begin{array}{c} \exists y_2 \ne y_1 \in \partial \Omega, \quad \mathrm{d}(x,y_1) = \mathrm{d}(x,y_2),
\\[5pt]  l(x,y_1), \ l(x,y_2) \in\Omega, \\[5pt]  l(x,y_1)\perp \partial_{\tau} (y_1), \  l(x,y_2)\perp \partial_{\tau} (y_2).  \end{array}\right.  \right\}.
\end{equation}
So for a skeleton point $x\in\mathcal{S}_{\Omega}$, which may have multiple sets of contributing boundary points satisfying \eqref{eqn:2d:biharm:5}, the function $s_{\Omega}(x)$ is the shortest distance from $x\in\mathcal{S}_{\Omega}$ to any of its contributing boundary points.

As $t$ increases, the set $\omega(t)$ propagates inwards from $\partial\Omega$ towards the \al{center} of the domain. At time $t=\al{T_{\eps}}$, the value of $u(x,t)$ predicted by \eqref{eqn:2d:biharm:3} will be quite different depending on whether $\omega(t)$ has intersected $\mathcal{S}_{\Omega}$ or not. To differentiate between these two distinct scenarios, the \emph{skeleton arrival time} $T_{\mathcal{S}}$ is introduced where
\begin{equation}\label{eqn:2d:biharm:7}
T_{\mathcal{S}} = \inf \{ t \,\left\bracevert \, \exists \; x\in\mathcal{S}_{\Omega}, \;  s_{\Omega}(x) = \eta_0 \phi(t; \eps) \; \right. \}.
\end{equation}
The value of $T_{\mathcal{S}}$ indicates the shortest time at which $\omega(t)$ intersects with the skeleton $\mathcal{S}_{\Omega}$. The dichotomy of possible singularity locations predicted by \eqref{eqn:2d:biharm:3} is now the following; If $\al{T_{\eps}}< T_{\mathcal{S}}$, then singularities are predicted to occur on $\omega(\al{T_{\eps}})$, typically at discrete points selected by the higher order curvature effects of \eqref{eqn:2d:biharm:3}. On the other hand, for $\al{T_{\eps}}\geq T_{\mathcal{S}}$, equation \eqref{eqn:2d:biharm:3} predicts that the singularity will occur on $\mathcal{S}_{\Omega}$ at the point(s) for which $s_{\Omega}(x) = \eta_0 \phi(\al{T_{\eps}};\eps)$.

\section{Application of the asymptotic theory}\label{sec:app}

In the following subsections, we apply the aforementioned asymptotic theory to predict the singularity set of \alt{problems} \eqref{intro_2} for a variety of spatial regions $\Omega$. As suggested by the analysis, the singularity set for the second order case \eqref{intro_2a} and the fourth order case \eqref{intro_2b} can be very different. 

\subsection{Example: Disc}\label{sec:ex_disc}

For the example of the disc geometry, $\Omega = \{(x_1,x_2)\in \mathbb{R}^2 \ \left\bracevert \ x_1^2 + x_2^2 \leq 1 \right. \}$, we compare predictions of the theory for the second order and fourth order cases. The theory for the second order case, developed in \S\ref{sec:2d:lap}, results in the simple prediction (cf. equation \eqref{eqn:2d:9}) that the singularity should occur at global maxima of the distance function $\mathrm{d}(x,\partial\Omega)$. For the unit disc case, the global maximum of this function is simply $\{(0,0)\}$. This prediction is in agreement with the results of \cite{AFMC85}.

To apply the theory developed in \S\ref{sec:2d:biharm} for the fourth order case \eqref{intro_2b}, we first note that the skeleton of the domain is $\mathcal{S}_{\Omega} = \{(0,0)\}$ and therefore $T_{\mathcal{S}}>0$. The predicted \al{behavior} is the following. If $\al{T_{\eps}}<T_{\mathcal{S}}$, singularities will form on $\omega(\al{T_{\eps}})$ which for this example is a ring or radius $1 - \eta_0\phi(\al{T_{\eps}},\eps)$ where $\phi(t;\eps) = \eps u_0^{1/4}$ and $u_{0t} = f(u_0)$. For $\al{T_{\eps}}>T_{\mathcal{S}}$, the singularity is predicted to occur at the origin. As the following example will clarify, this situation only \al{materializes} in the radially symmetric case - in general the solution develops an instability along $\omega(t)$, which results in a single point selected for singularity.

\begin{figure}[htbp]
\centering
\subfigure[Ring Blow-up, $\eps=0.1$.]{\includegraphics[width=0.45\textwidth]{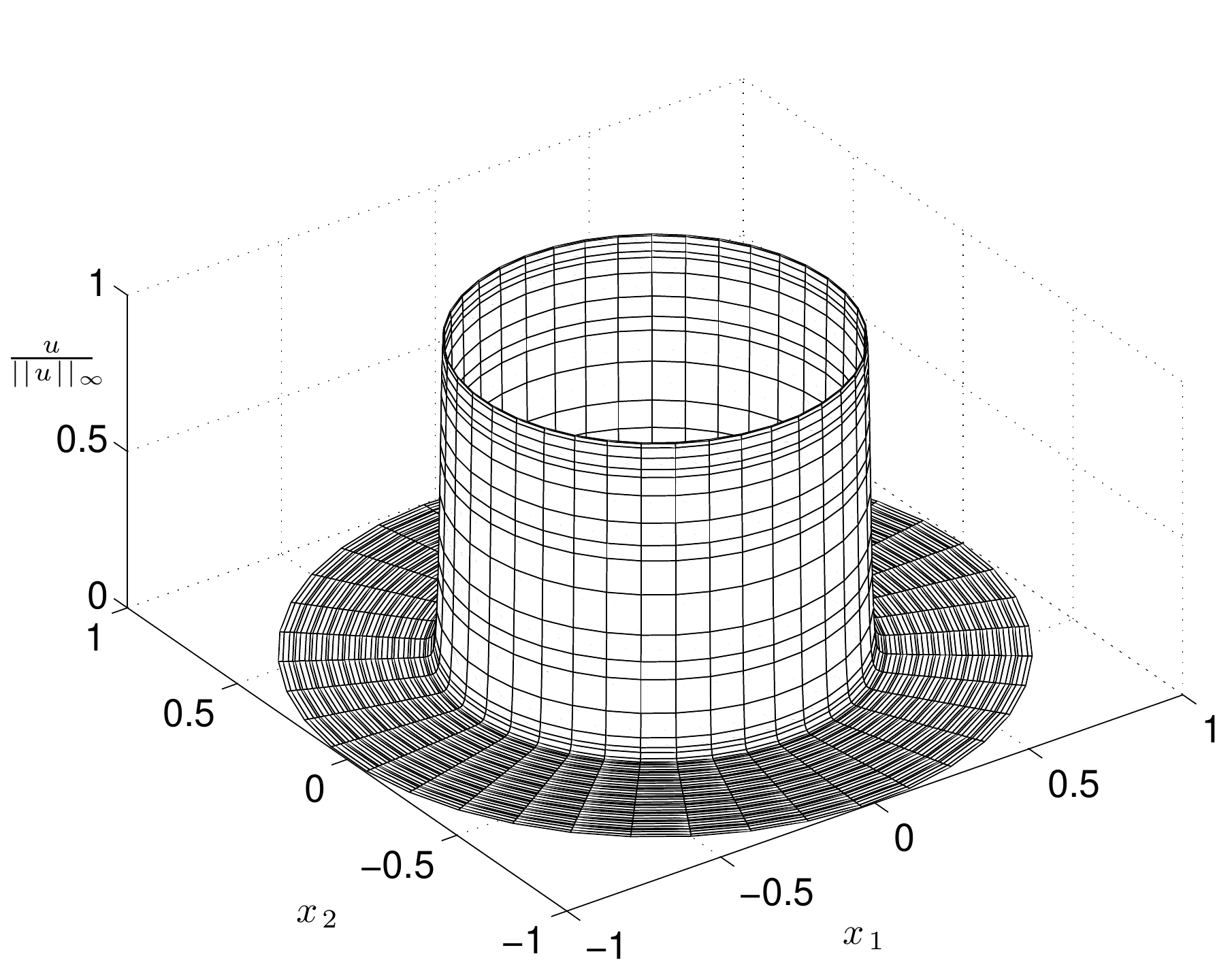}\label{fig:example_disc_a}}\qquad
\subfigure[Blow-up ring radius.]{\includegraphics[width=0.425\textwidth]{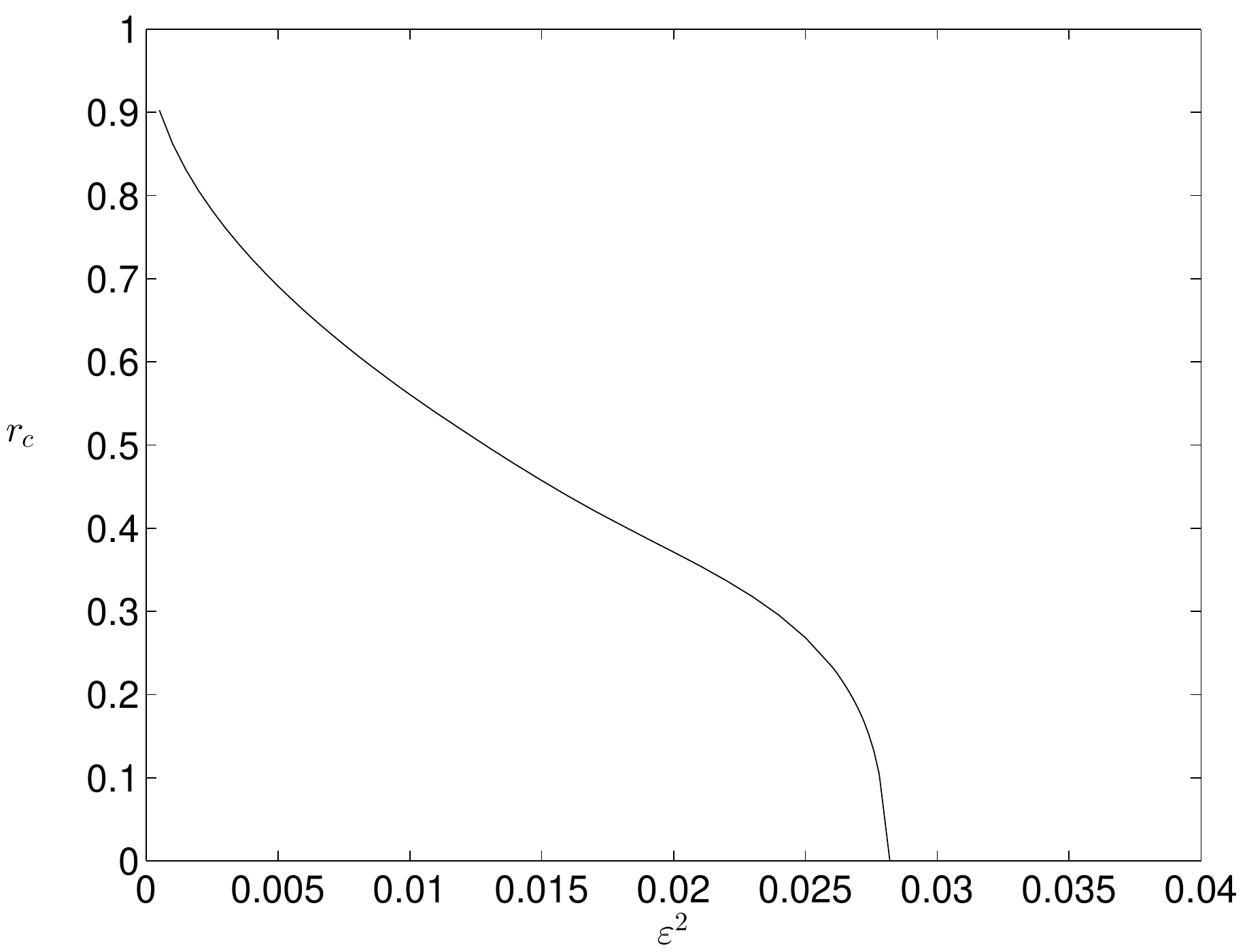}{\label{fig:example_disc_b}}}
\parbox{5in}{\caption{Numerical solutions of \alt{problem} \eqref{sec:ex_disc_radial}. Left panel, 
the solution profile for $\eps=0.1$ very close to blow-up $(\|u\|_{\infty}=10^{10})$. The singularity occurs simultaneously along an inner ring of points as predicted by the leading order asymptotic theory. Right panel, the radius of the blow up ring as a function of $\eps^2$.
\label{fig:example_disc}}}
\end{figure}

The consideration of radially symmetric solutions to $\eqref{intro_2b}$, along with the power nonlinearity $f(u) = (1+u)^2$, yields the reduced \alt{problem}
  \begin{equation}\label{sec:ex_disc_radial}
 \begin{array}{lc} u_t = -\eps^4 \Big[ u_{rrrr} + \ds\frac{2}{r}u_{rrr} -\ds\frac{1}{r^2} u_{rr} + \ds\frac{1}{r^3}u_{r} \Big] + (1+u)^2,& \ (r,t)\in (0,1)\times(0,\al{T_{\eps}}); \\[5pt]
 u = u_{r} = 0, \quad& r=1, \ t \in (0,\al{T_{\eps}});\\[5pt]
 u_{r} = u_{rrr} = 0, \quad& r=0, \ t \in (0,\al{T_{\eps}});\\[5pt]
 u = 0, & r\in(0,1), \quad t = 0,\end{array}
  \end{equation}
where $r=\ \sqrt{x_1^2+x_2^2}$. In Fig.~\ref{fig:example_disc_a}, a solution profile of \eqref{sec:ex_disc_radial} close to blow-up ($\|u\|_{\infty}=10^{10}$) is displayed for $\eps=0.1$. As predicted by the asymptotic theory, the solution is observed to blow-up simultaneously along an inner ring of points. The dependence of the radius of the blow-up ring and $\eps$ is illustrated in Fig.~\ref{fig:example_disc_b}.

\begin{figure}[htbp]
\centering
\subfigure[Instability along $\omega(t)$.]{\includegraphics[width=0.45\textwidth]{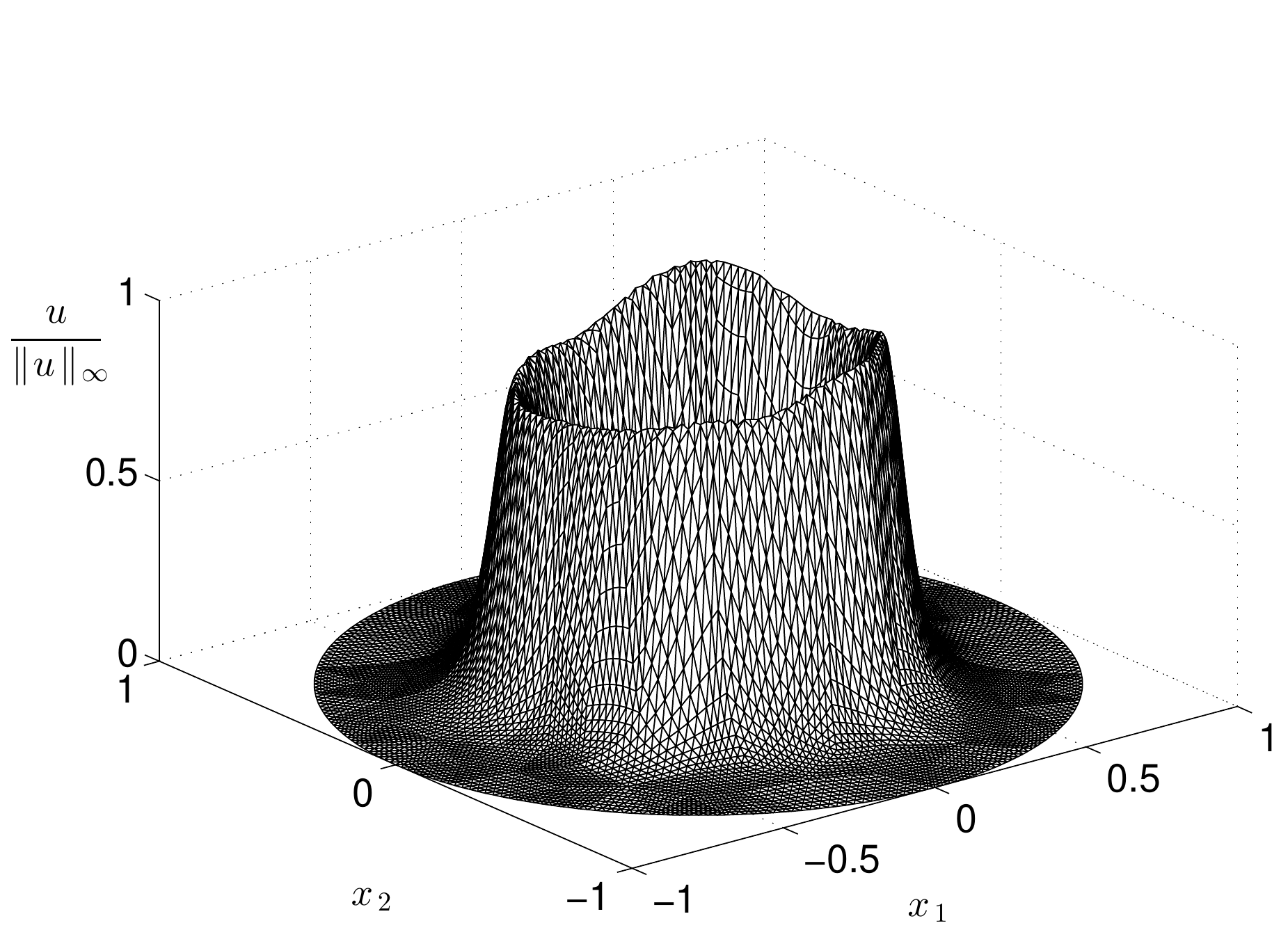}\label{fig:example_disc_noise_a}} \qquad
\subfigure[Single point blow-up.]{\includegraphics[width=0.45\textwidth]{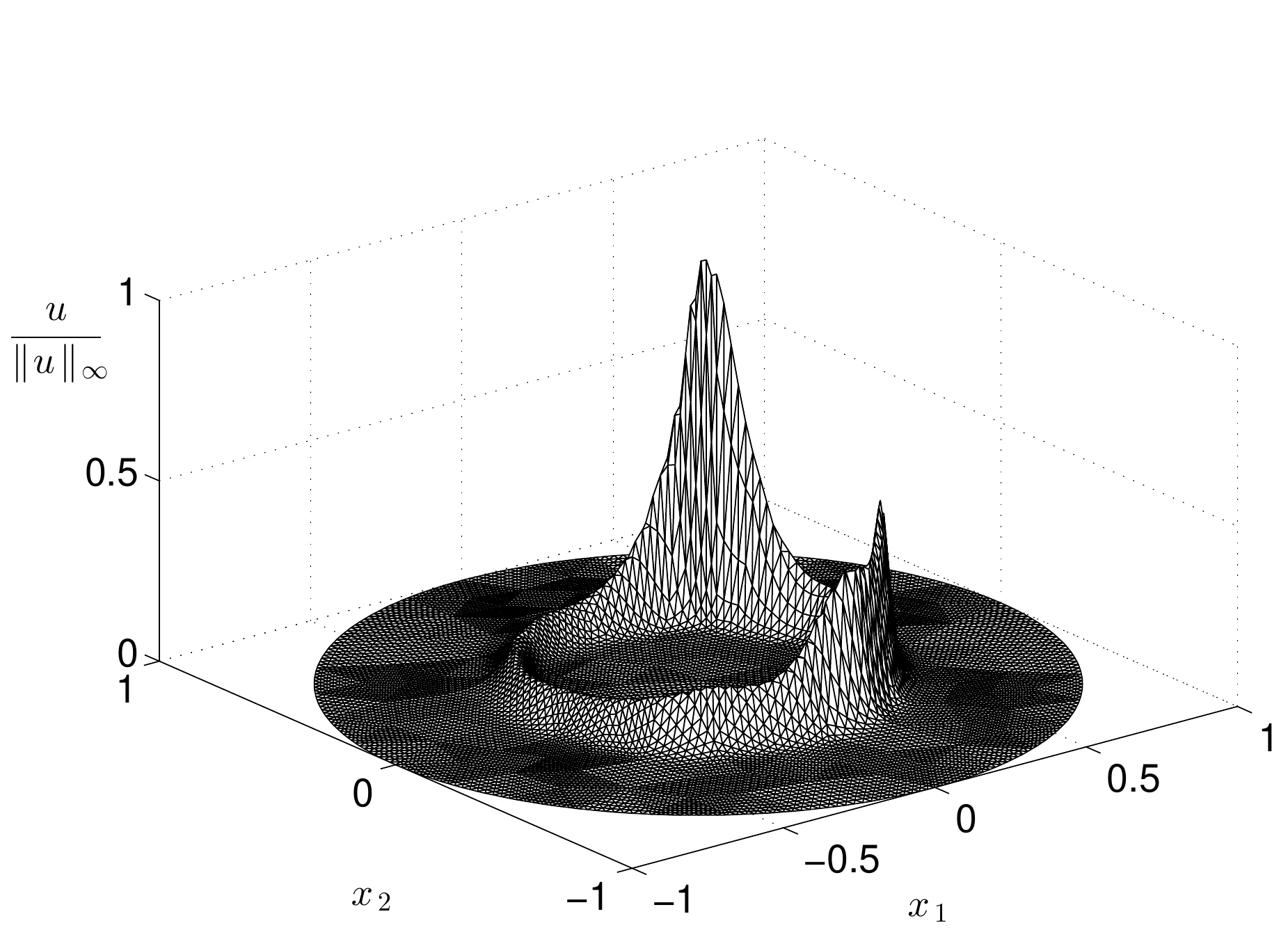}\label{fig:example_disc_noise_b}}
\parbox{5in}{\caption{Numerical solutions of \alt{problem} \eqref{intro_2b} with $f(u) = (1+u)^2$ and $\eps=0.1$ on the unit disc with noisy initial data of amplitude $5\times10^{-3}$. \al{The left and right panels display the solution profile for $\|u\|_{\infty} = 1\times10^{3}$ and $\|u\|_{\infty} = 1\times10^{5}$ respectively.} The initial noise is smoothed out and an instability develops along $\omega(t)$. The developing instability combines with the dynamics of the PDE to select a single point for blow-up.
\label{fig:example_disc_noise}}}
\end{figure}

However, this blow up ring solution is not stable in a full two dimensional setting. A similar instability was described in \cite{BW4} for ring rupture in thin-film equations. In Fig.~\ref{fig:example_disc_noise}, we display a full 2D (non-radially symmetric) solution of \eqref{intro_2b} initialized with small amplitude noise. As Fig.~\ref{fig:example_disc_noise_a} indicates, the initially noisy data is smoothed out by the dynamics of the PDE and the structure of $\omega(t)$ emerges. An instability develops along $\omega(t)$ in the angular direction, consequently the dynamics of the PDE selects a single point on $\omega(t)$ for blow-up, as shown in Fig.~\ref{fig:example_disc_noise_b}. A repetition of this numerical experiment for different \al{realizations} of initial data, results in blow-up points which are uniformly distributed along $\omega(t)$.

\subsection{Example: Square}\label{sec:ex_square}

We now compare the predictions of the theory for the second order case and the fourth order case for the square geometry $\Omega = [-1,1] ^2$. In the second order case \eqref{intro_2a}, the short time asymptotic theory predicts the critical point to be $x_c = \max_{x\in\Omega} \mathrm{d}(x,\partial\Omega)$, which gives $x_c = \{(0,0)\}$ in this scenario.

\begin{figure}[htbp]
\centering
\subfigure[Square Skeleton]{\includegraphics[width=0.36\textwidth]{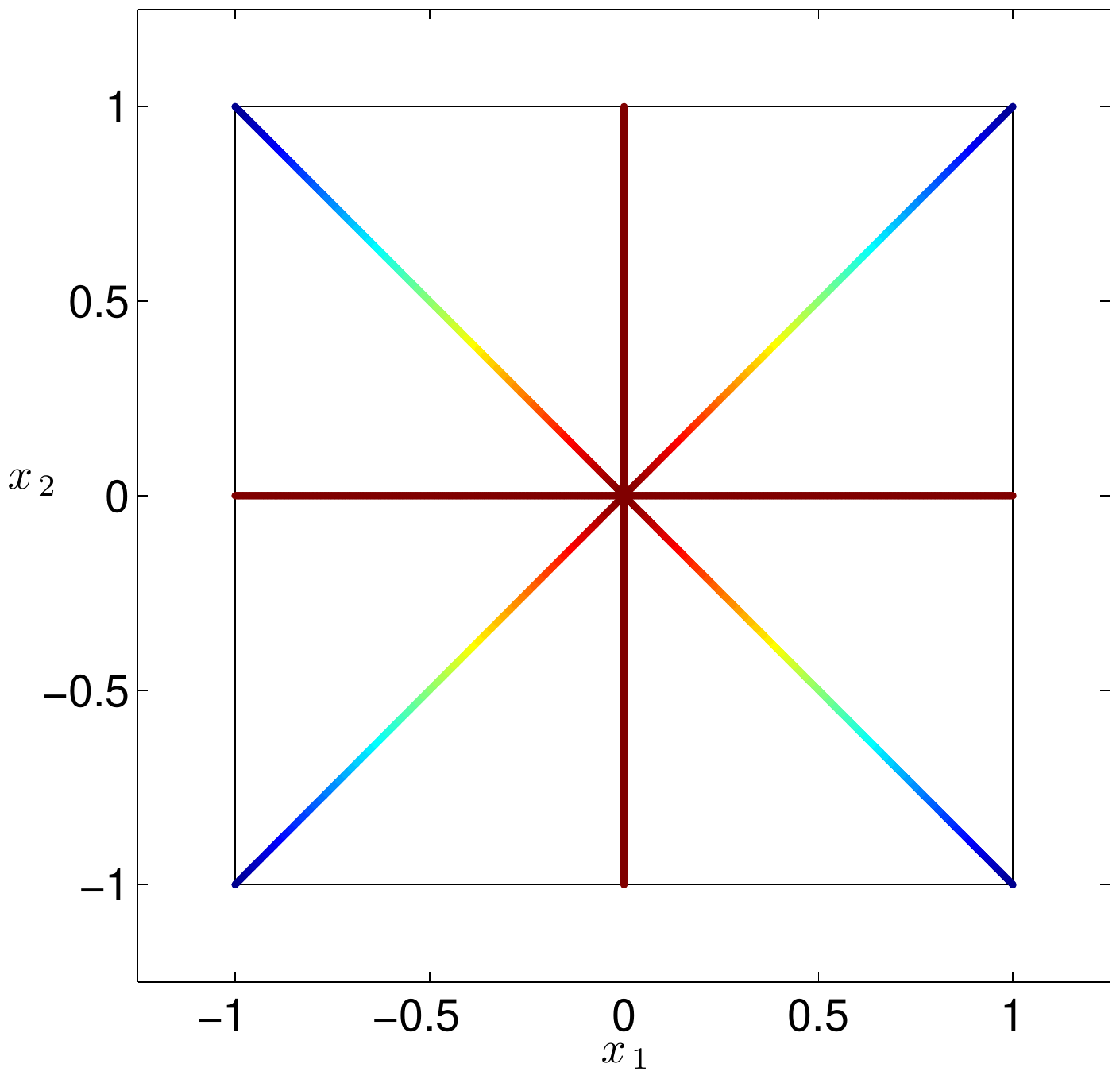}\label{fig:ex51_a}}
\qquad
\subfigure[Singularity Points]{\includegraphics[width=0.46\textwidth]{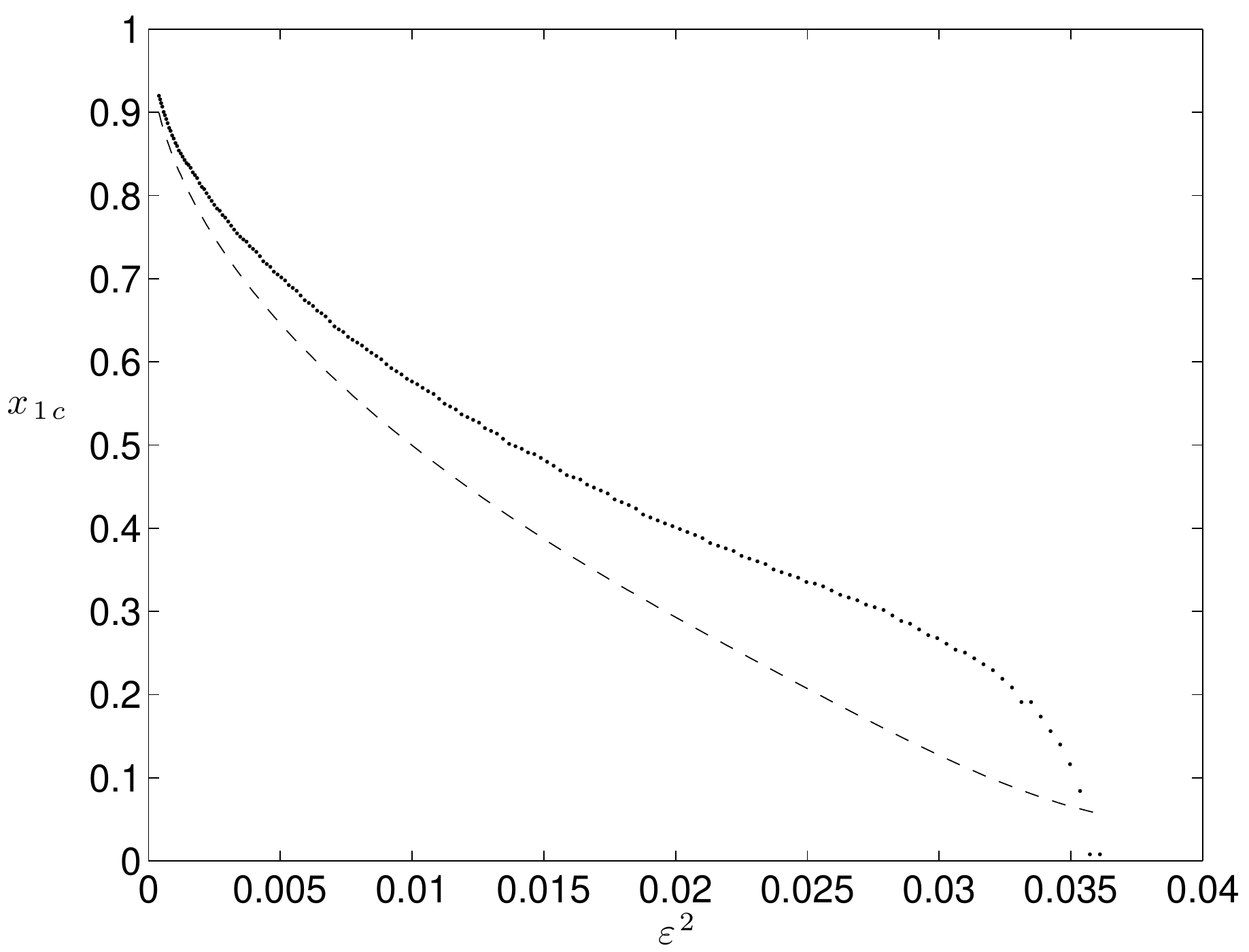}\label{fig:ex51_b}}\\[5pt]
\parbox{5in}{\caption{In panel (a), the skeleton of the square region is displayed with the shading indicating values of $s_{\Omega}(x)$. In panel (b), numerical (solid dots) and asymptotic (dashed line) predictions of the blow-up points in the first quadrant are shown for a range of $\eps^2$ values. The numerical blow up points are obtained from simulations of \eqref{intro_2b} with the nonlinearity $f(u)=e^u$ and $u=0$ initial data.\label{fig:ex51}}}
\end{figure}

To apply the theory for the fourth order problem, the first step is to construct the skeleton $\mathcal{S}_{\Omega}$. For this example, $\mathcal{S}_{\Omega}$ consists of the $x_1\,x_2$ axes and the diagonals of the square, \emph{i.e},
\begin{equation}\label{exsquare1}
\mathcal{S}_{\Omega} = \{ (x_1,x_2) \in \mathbb{R}^2 \; \left\vert \; x_1x_2 = 0 \; \mbox{or} \;  x_2 = \pm\, x_1\right. \} \cap \Omega.
 \end{equation}
An inspection of the skeleton (cf. Fig.~\ref{fig:ex51_a}) reveals that $T_{\mathcal{S}}=0$ and so the singularities are predicted to form on $\mathcal{S}_{\Omega}$. 

As the values of $s_{\Omega}(x)$ on the diagonal portion of $\mathcal{S}_{\Omega}$ are lower than those lying on the axes portion, the theory predicts that singularities will form simultaneously at four points on this diagonal for $\eps<\eps_c$ and at the origin for $\eps>\eps_c$. Moreover, the prediction of the $x_1$ locations of the singularities are given by
\begin{equation}\label{ex:square2}
x_{1c} = \left\{ \begin{array}{cl} \pm(1-\eta_0\phi_c(\eps)), & \eps\leq\eps_c, \\[5pt] 0, & \eps>\eps_c, \end{array}   \right.
\end{equation}
where the critical value $\eps_c$ is implicitly determined by $1 = \eta_0\phi_c(\eps_c)$. The $x_{2c}$ values are then $\pm x_{1c}$ to give a total of four distinct points in the case $\eps<\eps_c$.

\begin{figure}[htbp]
\centering
\subfigure[$\eps = 0.2$, One Singularity.]{\includegraphics[width=0.4\textwidth]{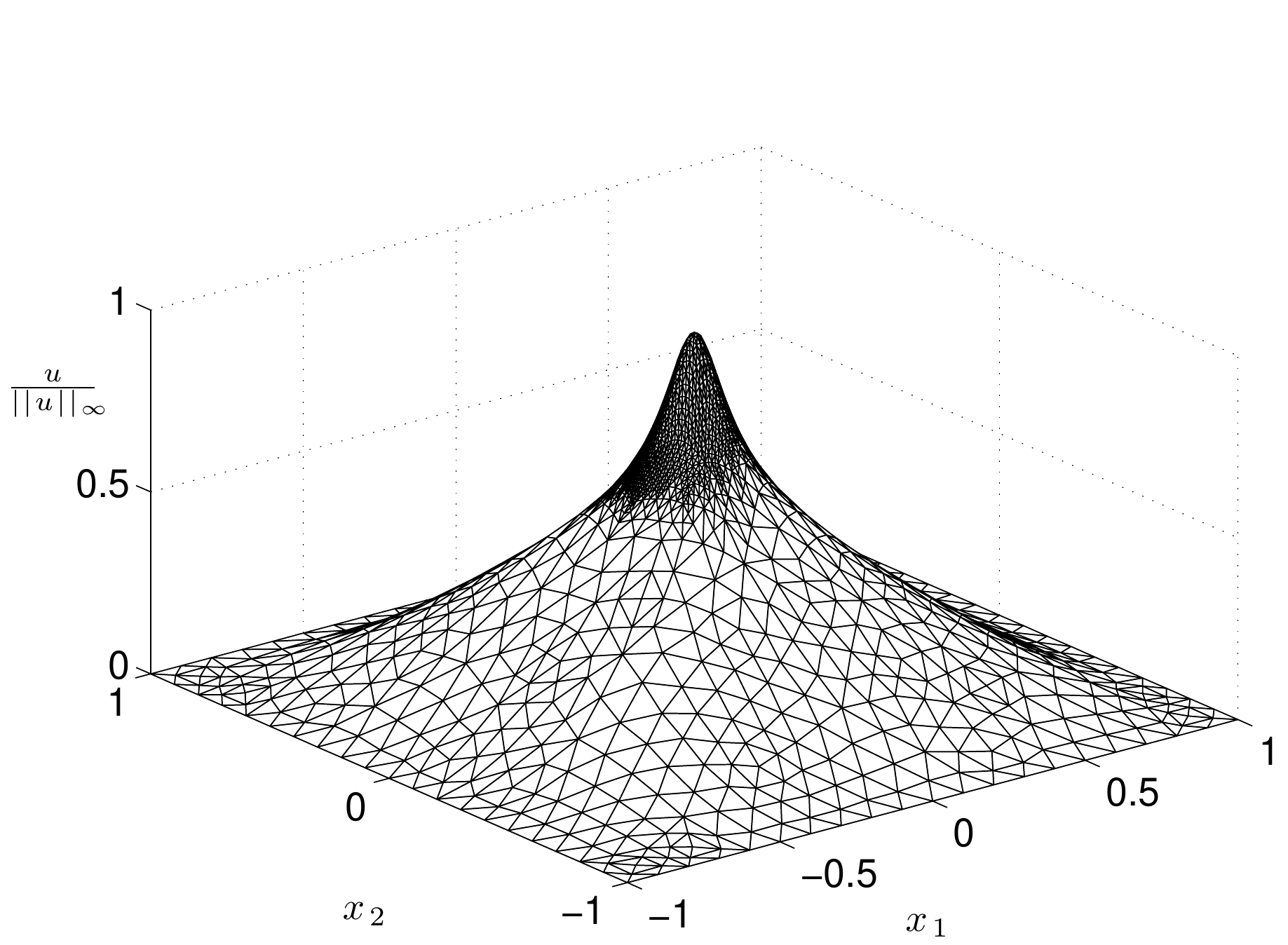}\label{fig:ex5c}}
\quad
\subfigure[$\eps = 0.1$, Four Singularities.]{\includegraphics[width=0.4\textwidth]{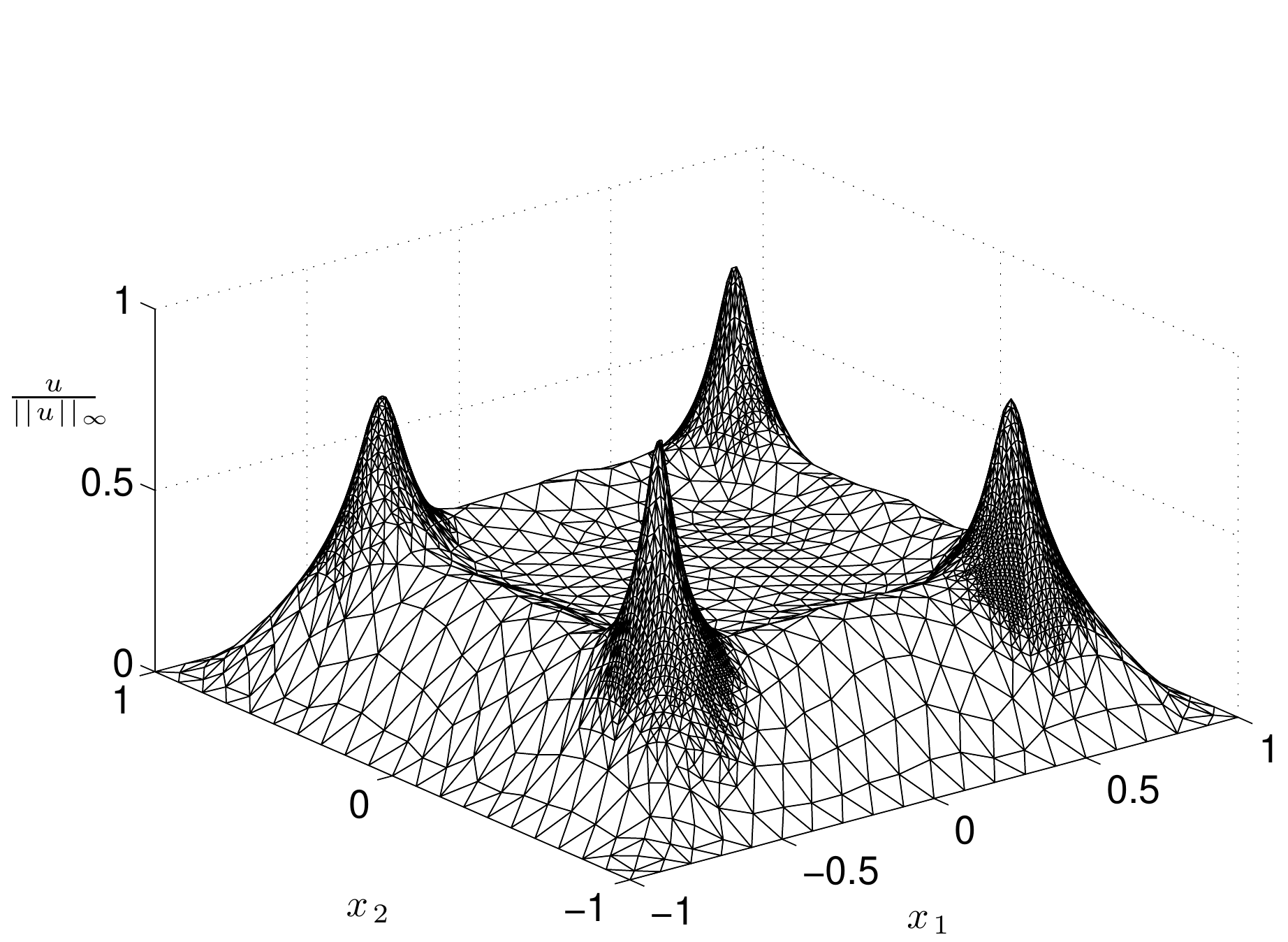}\label{fig:ex5d}}
\parbox{5in}{\caption{Numerical simulations of \eqref{intro_2b} for the square geometry and nonlinearity $f(u) = e^u$.  In panels (a) and (b), solution profiles of \eqref{intro_2b} for $\eps=0.2$ and $\eps=0.1$ are displayed after integration to $\|u\|_{\infty}=10$ and $\al{T_{\eps}}-t =\mathcal{O}(10^{-7})$. For these parameter values, the solution is observed to develop singularities at four discrete points of $\mathcal{S}_{\Omega}$ for $\eps=0.1$ and at the origin only for $\eps=0.2$. \label{fig:ex5}}}
\end{figure}
In Fig.~\ref{fig:ex51_b}, we see qualitative agreement between the numerical and asymptotic predictions and good quantitative agreement when $\eps$ is small. However, the asymptotic theory does not predict $\eps_c$, the threshold between single and multiple blow-up, with high accuracy. This is not surprising since the assumptions which underpin the asymptotic theory: a uniform central region coupled to a propagating boundary effect, do not hold when $\eps\approx\eps_c$.

\al{In the presence of a small random perturbation to the initial data, the four point blow up configuration is not generically stable. Indeed, for simulations of \eqref{intro_2b} initialized with small random noise, the short time solution will develop four peaks, but small discrepancies in their amplitude will result in one being selected for blow-up by the PDE dynamics, before the remaining peaks are able able to fully develop. Over many realizations, one recovers that each of the four possible blow-up locations is selected with uniform probability $1/4$.}

\subsection{\al{Example: Rectangle}}\label{sec:ex_rectangle}

\al{In this section, the blow-up set for the rectangular region $\Omega = [-1,1]\times[0,1]$ is considered. This example shows two things; first, that the determination of the blow-up set can depend on considerations beyond the point symmetries of the domain. Second, that the multiplicity of singularities can change more than once for $\eps\in(0,\eps_c)$, in contrast to behavior seen in the previous 1D and square example. For the second order \alt{problem} \eqref{intro_2a}, the asymptotic theory predicts that singularities should form at the origin in the absence of any noise.}

\begin{figure}[htbp]
\centering
\subfigure[Domain with partial skeleton]{\includegraphics[width=0.475\textwidth]{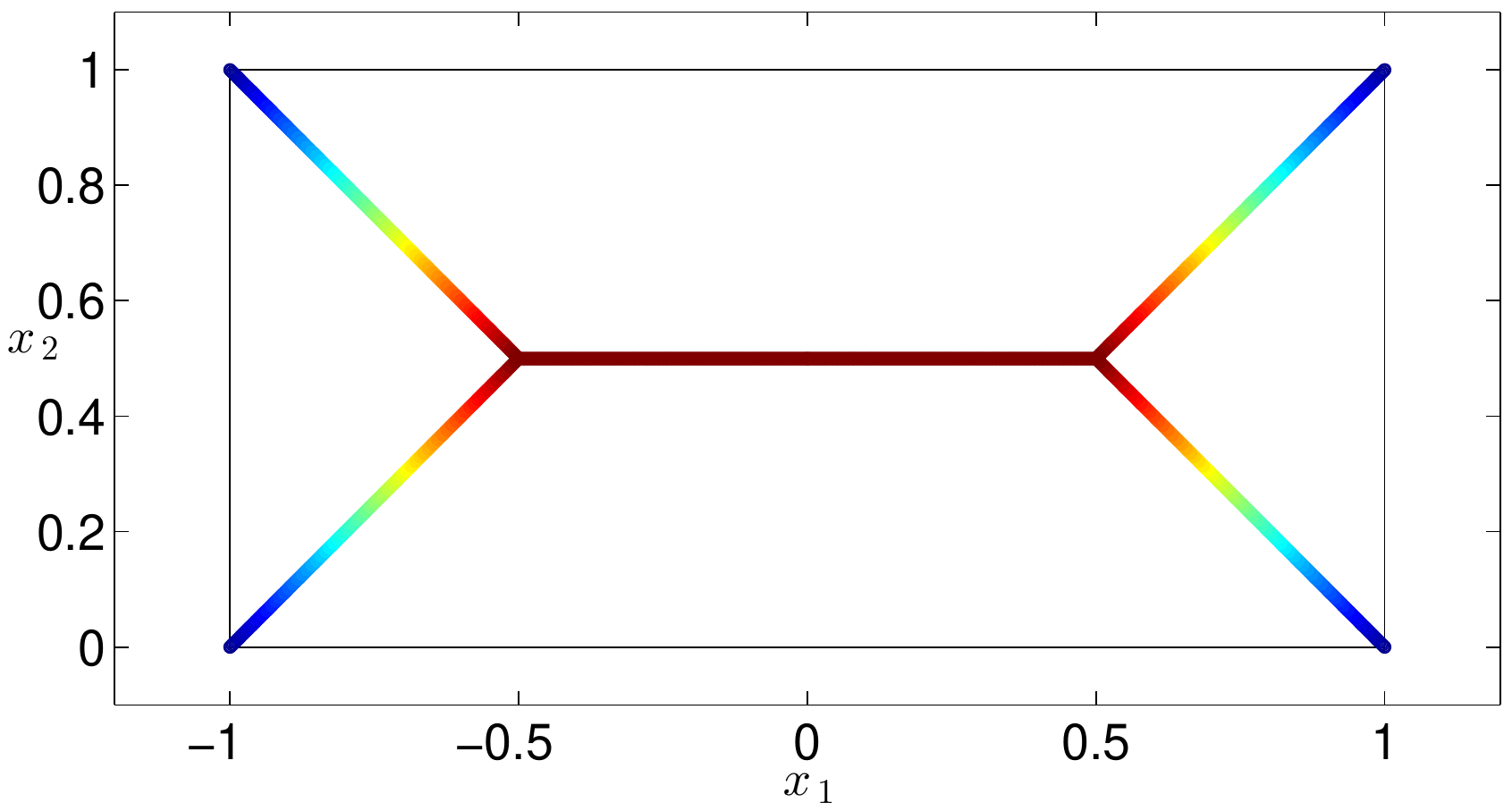}\label{fig:rectangle_1}}\qquad
\subfigure[$\eps=0.05$]{\includegraphics[width=0.45\textwidth]{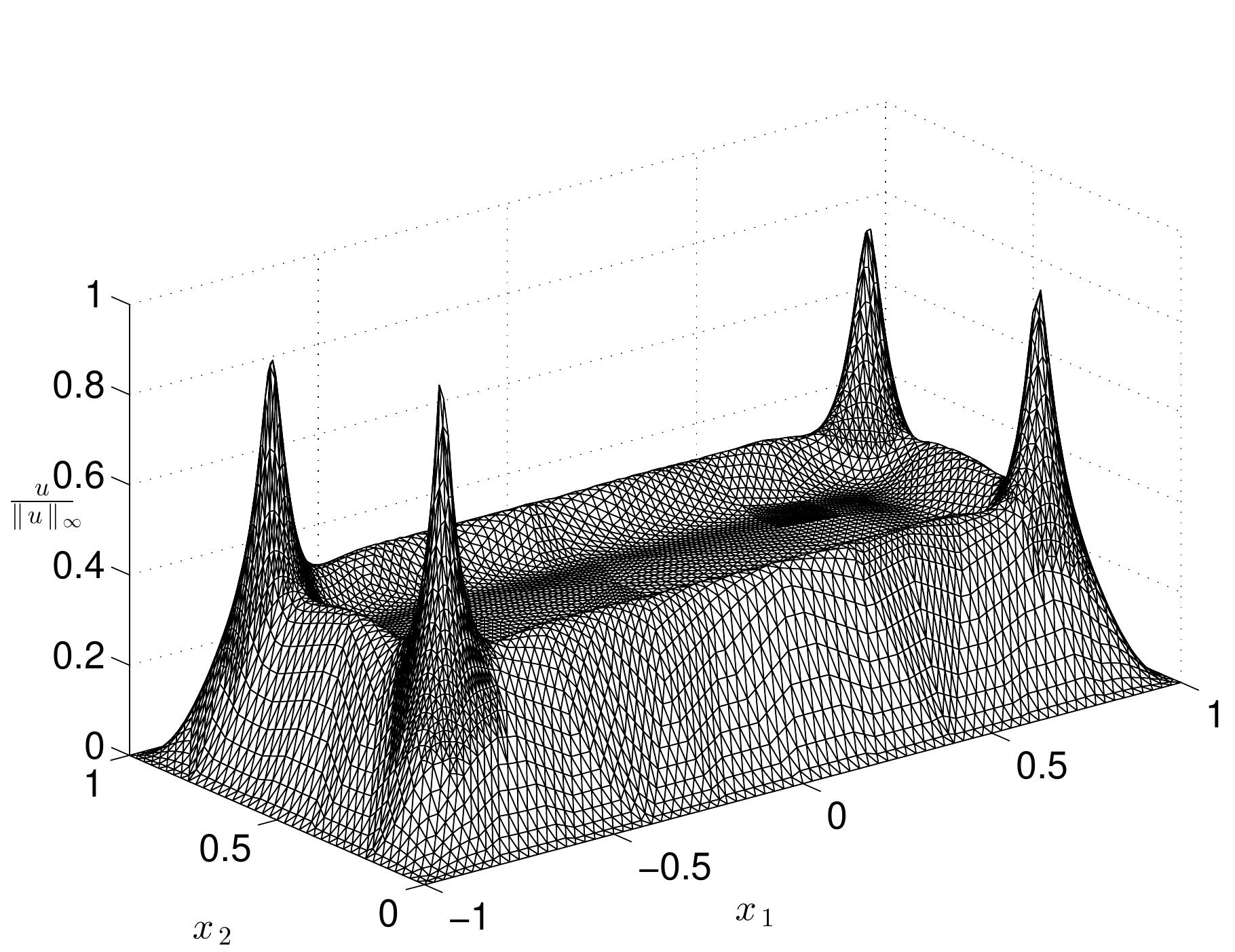}\label{fig:rectangle_2}}\\
\subfigure[$\eps=0.1$]{\includegraphics[width=0.45\textwidth]{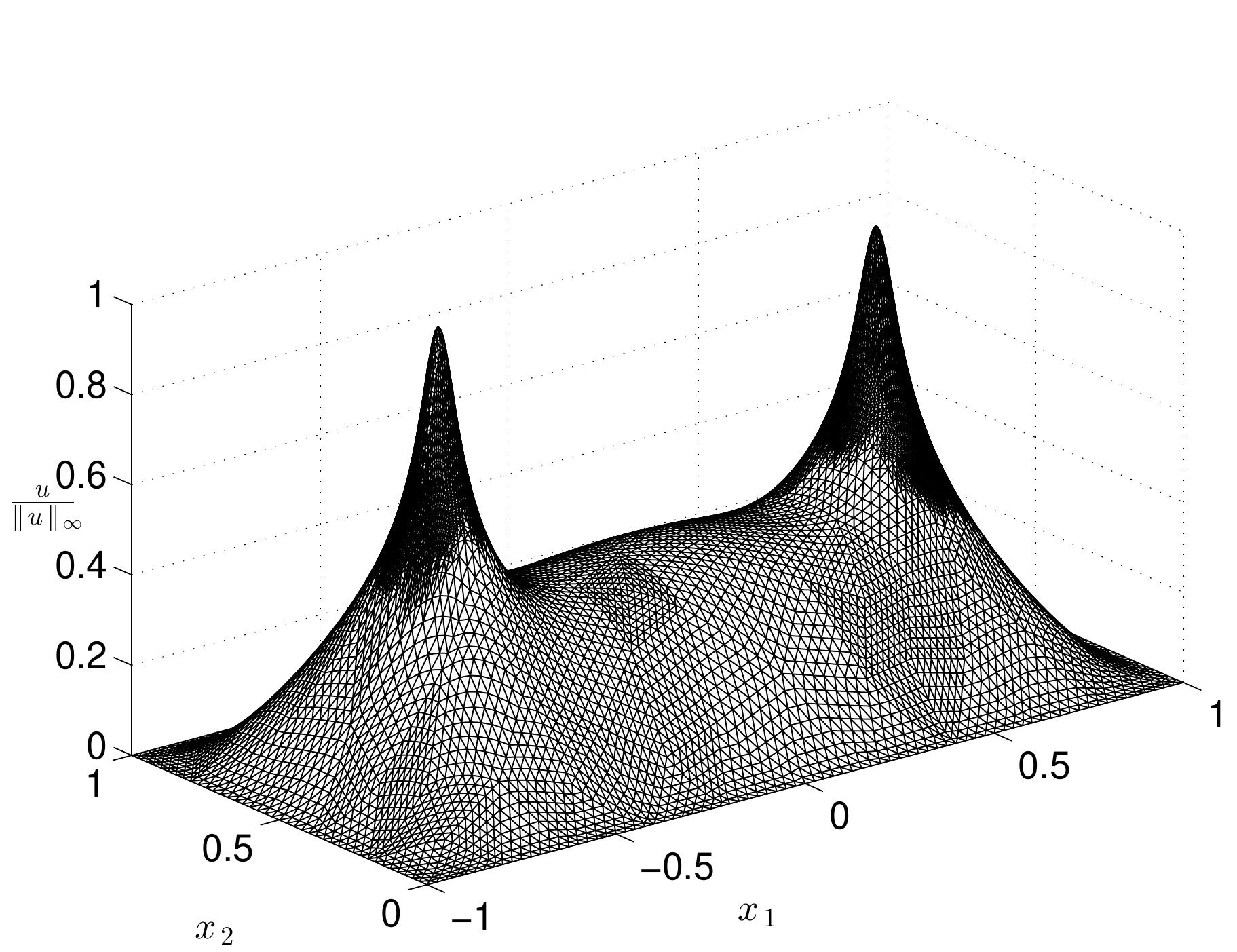}\label{fig:rectangle_3}}\qquad
\subfigure[$\eps=0.2$]{\includegraphics[width=0.45\textwidth]{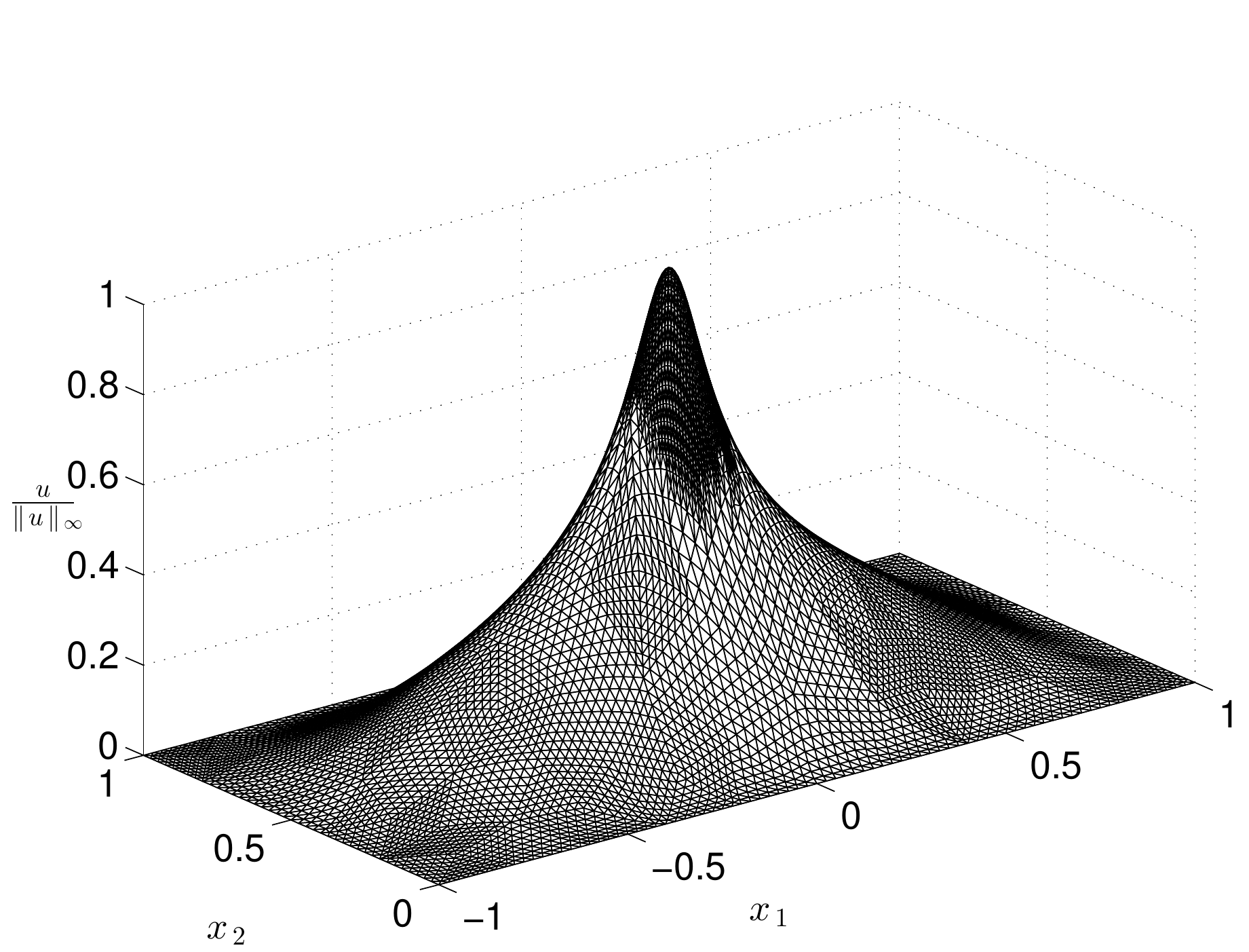}\label{fig:rectangle_4}}
\parbox{5in}{\caption{\al{Panel (a) shows the partial skeleton for the rectangular domain. Panels (b-d) show profiles of solutions to \eqref{intro_2b} for the nonlinearity $f(u) = e^u$ obtained from integration till $\|u\|_{\infty} = 10$.} \label{fig:rectangle}}}
\end{figure}

\al{In the fourth order case, we begin by considering the partial skeleton of the domain which takes the appearance of an \lq\lq envelope\rq\rq\ (cf. Fig.~\ref{fig:rectangle_1}). From the asymptotic theory, we therefore predict that in the absence of noise, \alt{problem} \eqref{intro_2b} should develop four singularities along each of the skeleton segments emanating from the corners of the rectangle. As $\eps$ increases, the location of the singularities moves along the segments until the left and right segments meet each other at the points $[\pm0.5,0]$. As these segments of $\mathcal{S}_{\Omega}$ do not meet directly at the origin as in the square case, the multiplicity of singularities does not change directly from four to one at this critical value of $\eps$. Instead, the multiplicity of singularities decreases to two, with this pair of singularities then being located along the central segment of the skeleton, $[-0.5,0.5]$.}

\al{The asymptotic theory consequently predicts the existence of two critical values $\eps_1<\eps_2<\eps_c$ such that four singularities occur for $\eps\in(0,\eps_1)$, two singularities occur when $\eps\in(\eps_1,\eps_2)$ and one singularity occurs when $\eps\in(\eps_2,\eps_c)$. The numerical experiments displayed in Fig.~\ref{fig:rectangle} illustrate these three outcomes for \alt{problem} \eqref{intro_2b} with the nonlinearity $f(u) = e^u$.}

\al{For the rectangular domain, we therefore have that the blow up set of \eqref{intro_2b} consists of either one, two or four points depending on the value of $\eps$. Moreover, the skeleton theory is able to capture possible singularity configurations of \eqref{intro_2b} that symmetry considerations alone would not. Indeed, \alt{problem} \eqref{intro_2b} may have a complex set of possible blow up location on domains whose skeletons consist of multiple branched segments.}

\subsection{Example: Domain with no particular symmetry}\label{sec:ex_potato}

In this section, the leading order asymptotic theory is applied to \alt{problems} \eqref{intro_2} on the region enclosed by the boundary
\begin{equation}\label{ex:potato1}
\partial \Omega = \{ (r(\theta)\cos\theta ,r(\theta)\sin\theta) \ \left\vert \ \right. 0\leq\theta\leq 2\pi \}, \qquad r(\theta) = 1 + 0.3\,( \cos\theta - \sin 3\theta ).
\end{equation}
The application of the theory to the second order problem \eqref{intro_2a} is relatively simple as the leading order prediction is simply $x^{\ast}_c = \max_{x\in\Omega} \mathrm{d}(x,\partial\Omega)$, independent of $\eps$. This point is calculated numerically to be $x^{\ast}_c = (0.3070,-0.0345)$ for this particular region. To compare this prediction with numerical simulations, the power nonlinearity $f(u) = (1+u)^2$ is chosen and $\eqref{intro_2a}$ is integrated until $\|u\|_{\infty}= 1\times10^{4}$. In Table.~\ref{table:example_potato} the $L^{\infty}$ error between the asymptotic and numerical blow up point predictions is shown to be very small, indicating the asymptotic and numerical predictions are in good agreement. The numerical simulations find the maximum value to occur at the same numerical node point for each value of $\eps$ - hence the errors are identical for each value of $\eps$. This is exactly as predicted by the leading order asymptotic prediction \eqref{eqn:2d:9}, which is independent of $\eps$.

\begin{table}[htbp]
\centering
\begin{tabular}{c|c|c|c}
$\eps$ & $0.1$& $0.15$ & $0.2$ \\[2pt]
\hline
$\|x_c-x_c^{\ast}\|_2$ & $0.0014$  &$0.0014$ & \al{$0.0014$} \\
\end{tabular}
\vspace{0.2in}
\parbox{5in}{\caption{Accuracy of blow-up point predictions for \alt{problem} \eqref{intro_2a} for region \eqref{ex:potato1} with $f(u) = (1+u)^2$ and various $\eps$. The asymptotic prediction $x_c^{\ast} = (0.3070,-0.0345)$ is determined from \eqref{eqn:2d:9} and the estimate $x_c$ is obtained from the maximum of the numerical solution profile at $\|u\|_{\infty}= 1\times10^{4}$.  \label{table:example_potato}}}
\end{table}

\begin{figure}[H]
\centering
\subfigure[Partial Skeleton]{\includegraphics[width=0.38\textwidth]{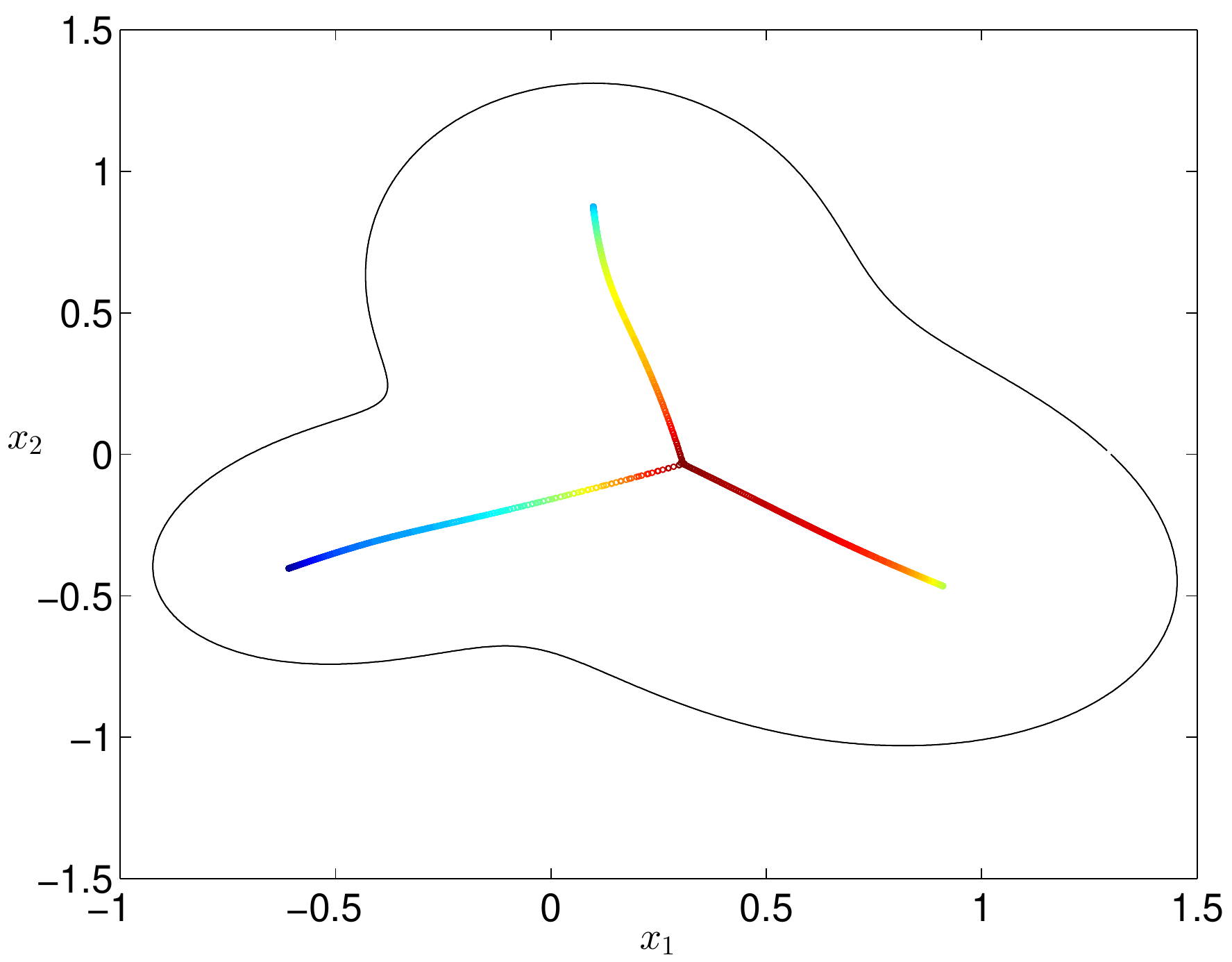} \label{fig:ex_potato2_a} }
\qquad
\subfigure[Blow up Set]{\includegraphics[width=0.4\textwidth]{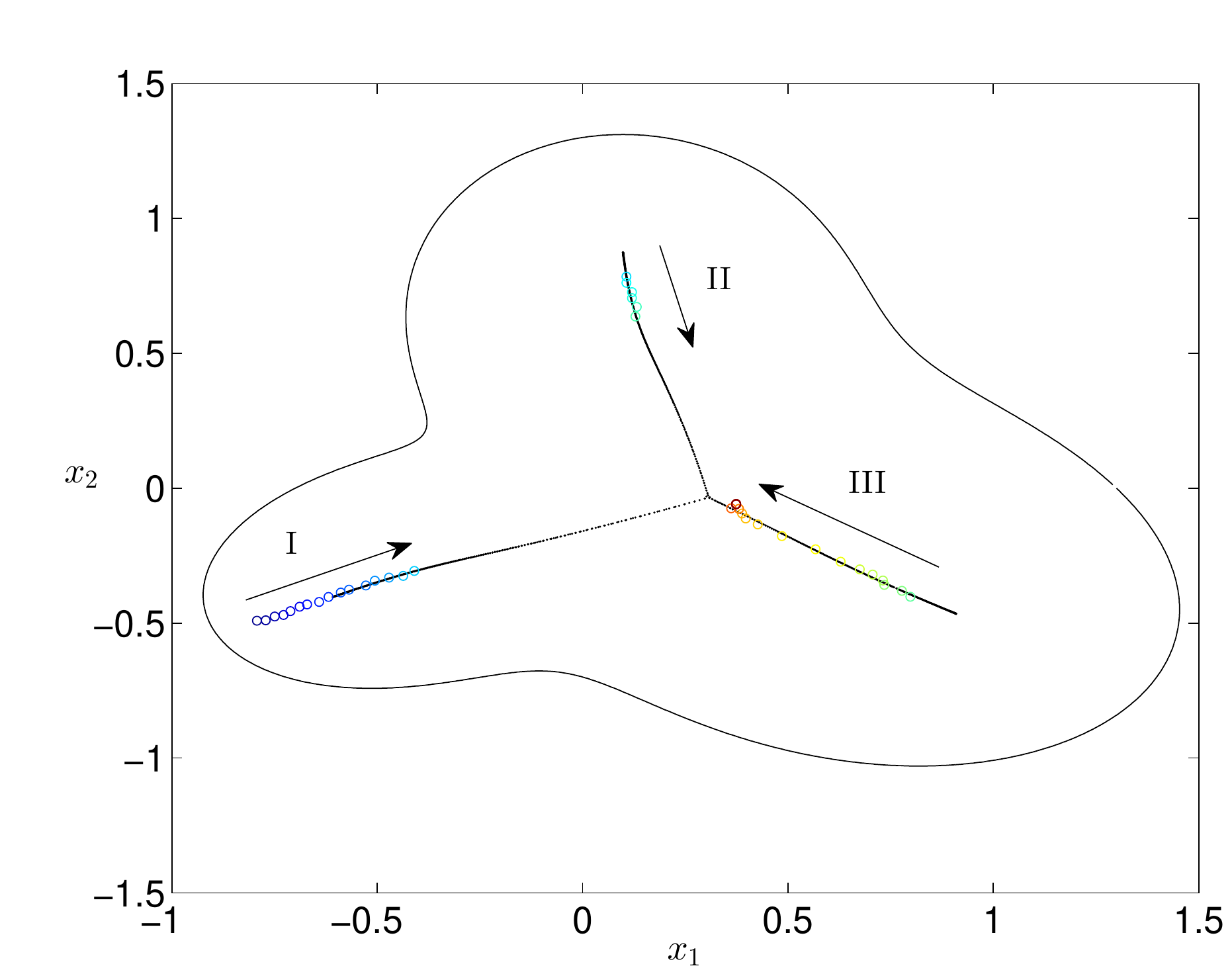} \label{fig:ex_potato2_b} }
\parbox{5in}{\caption{ Left Panel: The domain specified by \eqref{ex:potato1} together with the portion of the skeleton $\mathcal{S}_{\Omega}$, for which $s_{\Omega}(x) <0.735$. The shading along $\mathcal{S}_{\Omega}$ indicates the value of $s_{\Omega}$. Right Panel: Blow-up points overlaid on the skeleton with arrows indicating the direction of increasing $\eps$. As $\eps$ increases, blow-up occurs sequentially on branches $\mathrm{I}$, $\mathrm{II}$ and $\mathrm{III}$ of $\mathcal{S}_{\Omega}$.\label{fig:ex_potato2}}}
\end{figure}

As seen in the previous examples, the application of the leading order asymptotic theory to the fourth order problem \eqref{intro_2b} is considerably more delicate. The first step is to numerically calculate the skeleton $\mathcal{S}_{\Omega}$ for the region, which is displayed in Fig.~\ref{fig:ex_potato2_a}. An inspection of $\mathcal{S}_{\Omega}$ indicates that $T_{\mathcal{S}}>0$, \emph{i.e.} the skeleton arrival time is positive. Therefore, for $\al{T_{\eps}}<T_{\mathcal{S}}$ the leading order asymptotic theory predicts blow-up on discrete points of $\omega(\al{T_{\eps}})$ selected by curvature effects. For $\al{T_{\eps}}\geq T_{\mathcal{S}}$, the theory predicts blow-up at point(s) $x\in \mathcal{S}_{\Omega}$ for which $s_{\Omega}(x) = \eta_0 \phi(\al{T_{\eps}};\eps)$. 

\begin{figure}[htbp]
\centering
\subfigure[$\eps=0.05$]{\includegraphics[width=0.45\textwidth,clip]{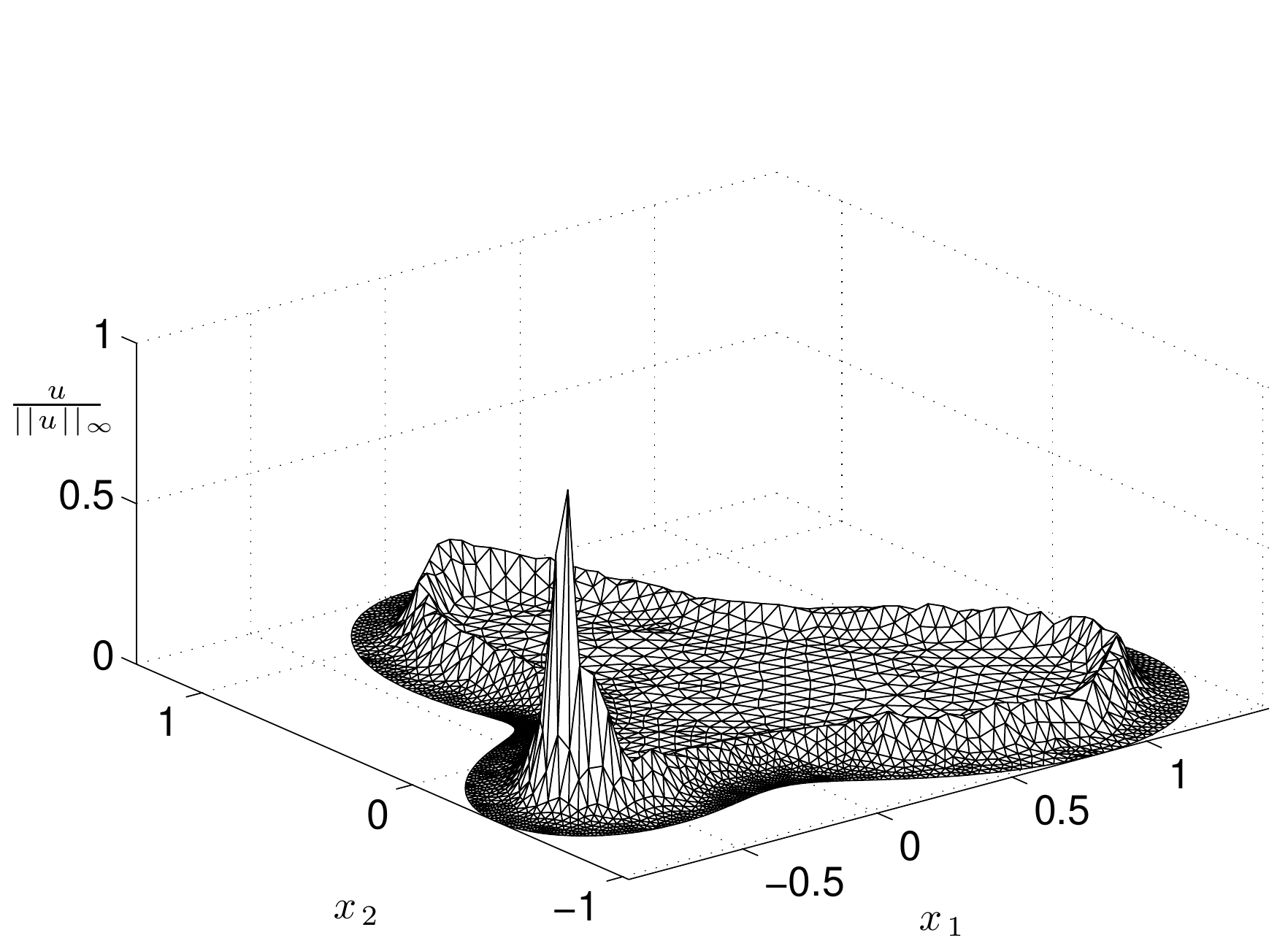} \label{fig:potato2_profile_a} }
\qquad
\subfigure[$\eps=0.1$]{\includegraphics[width=0.45\textwidth,clip]{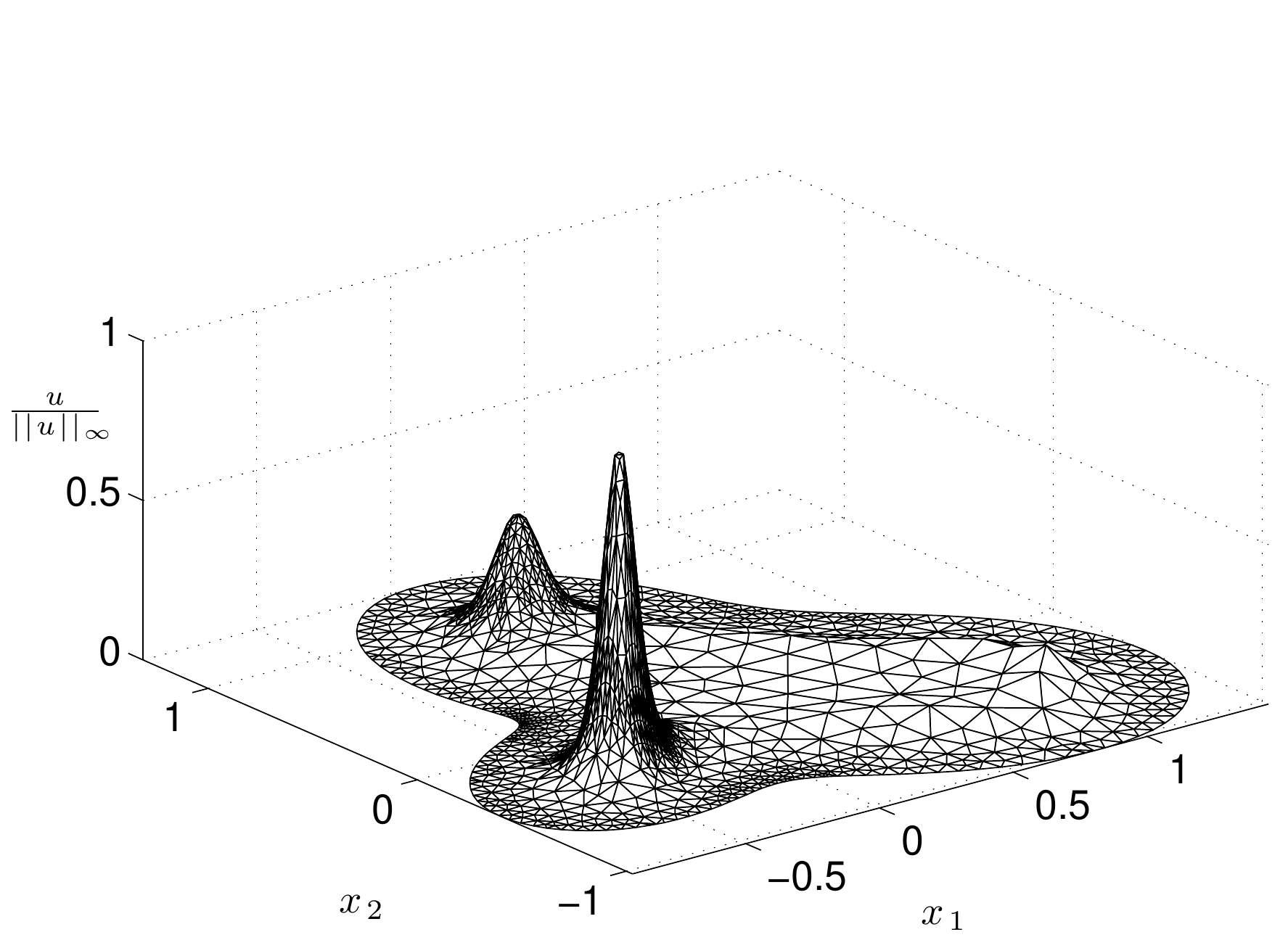} \label{fig:potato2_profile_b} }\\[-5pt]
\subfigure[$\eps=0.102875$]{\includegraphics[width=0.45\textwidth,clip]{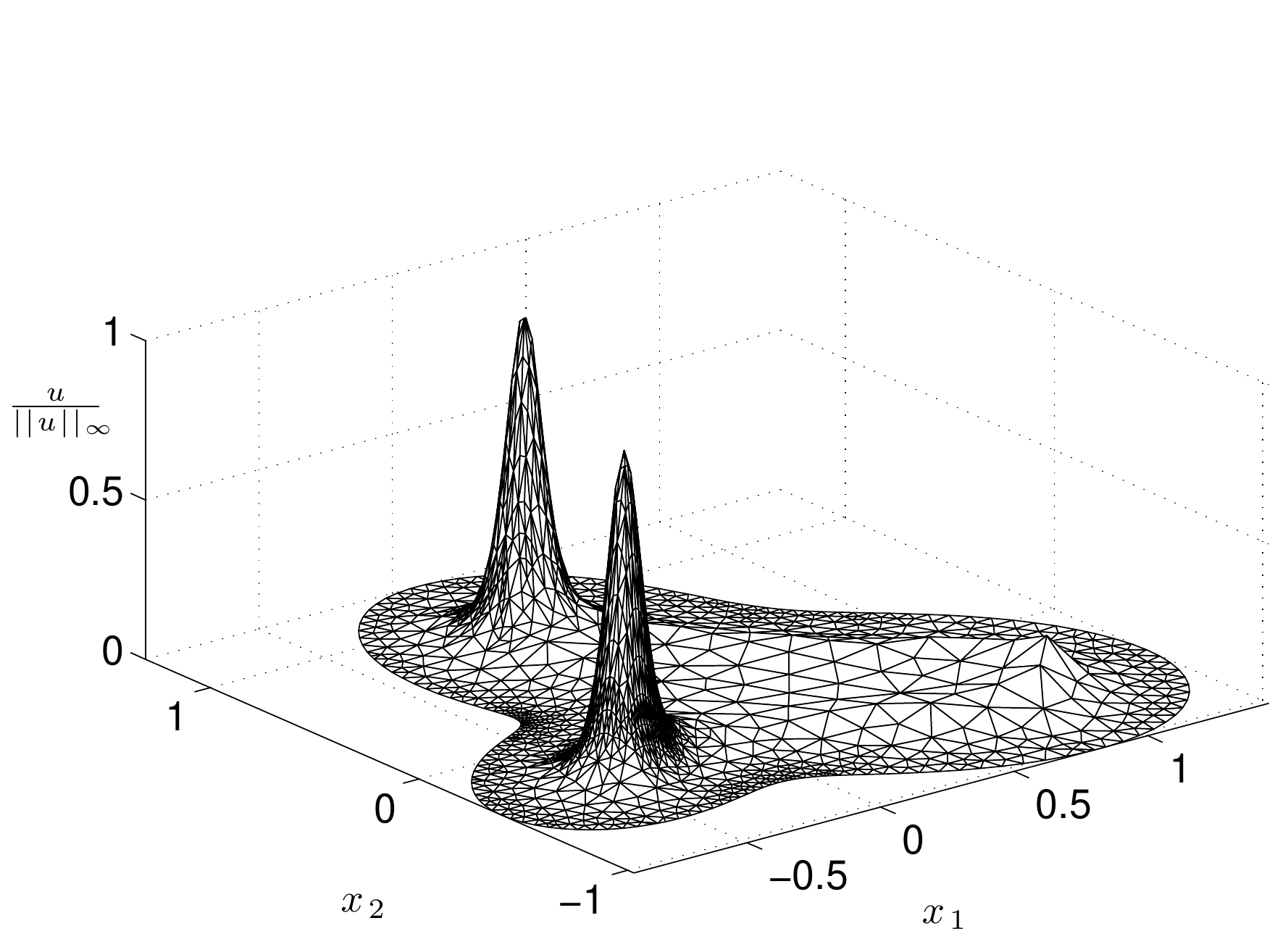} \label{fig:potato2_profile_c} }
\qquad
\subfigure[$\eps=0.13$]{\includegraphics[width=0.45\textwidth,clip]{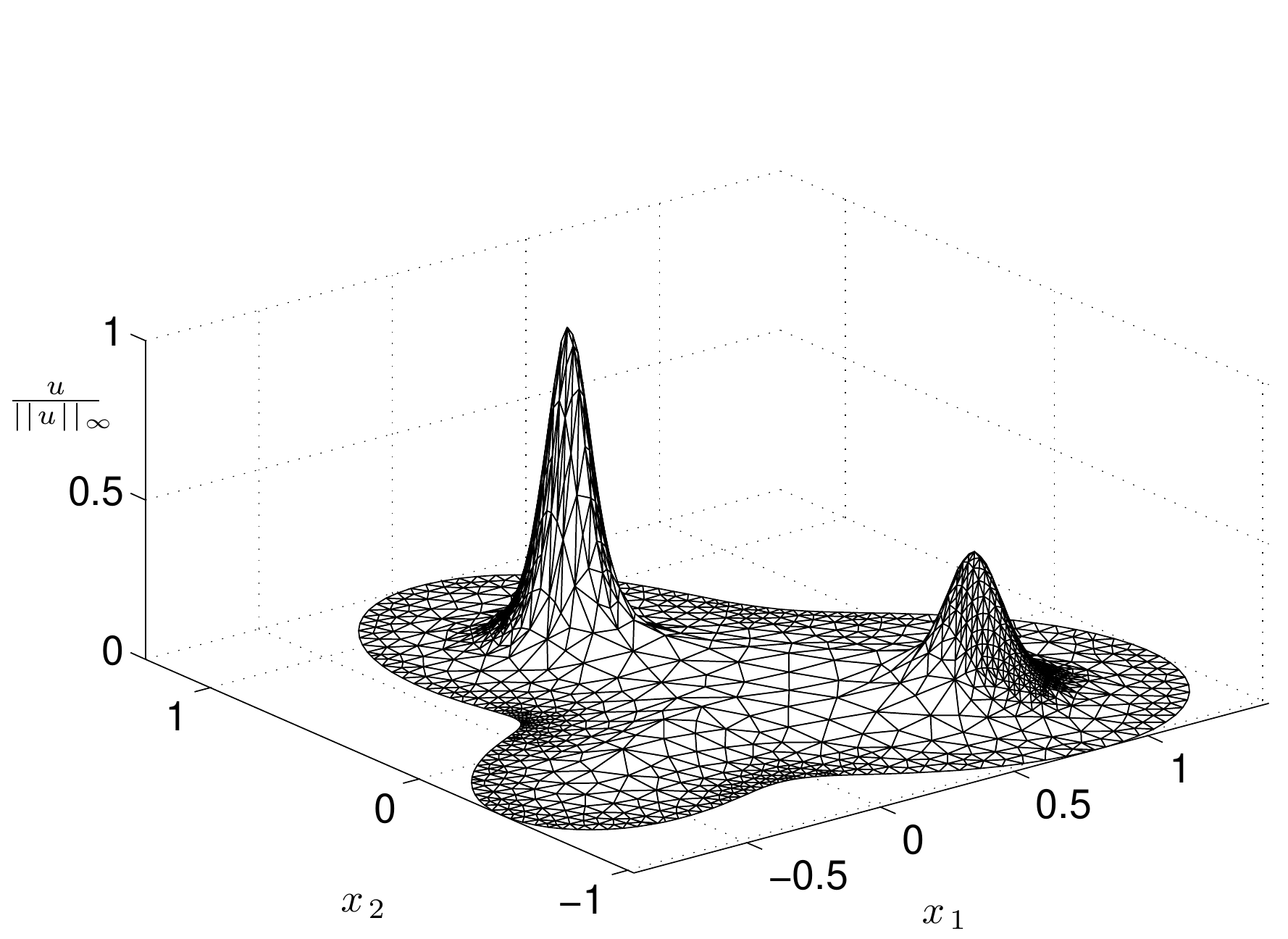} \label{fig:potato2_profile_d} }\\[-5pt]
\subfigure[$\eps=0.1343$]{\includegraphics[width=0.45\textwidth,clip]{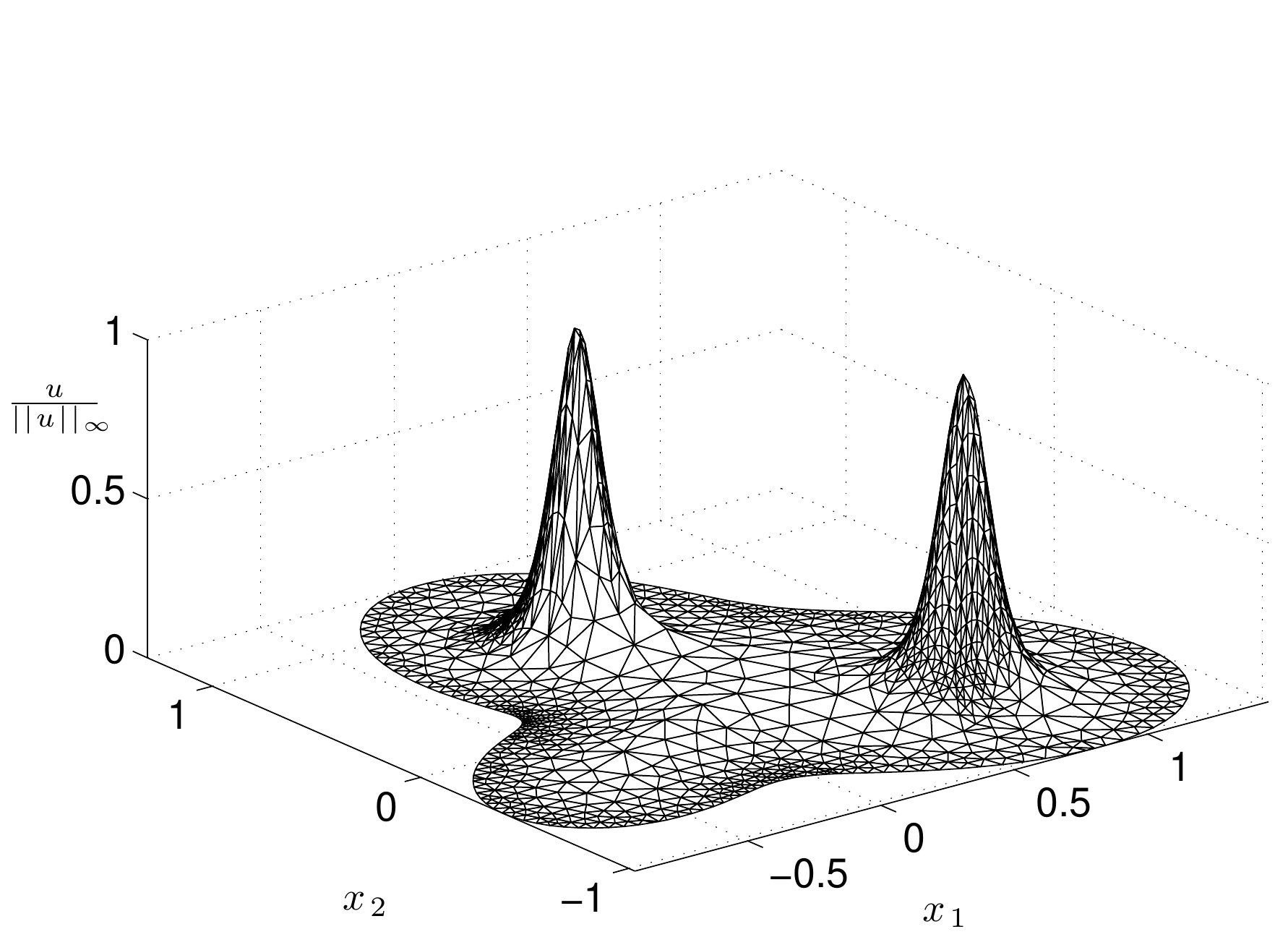} \label{fig:potato2_profile_e} }
\qquad
\subfigure[$\eps=0.14$]{\includegraphics[width=0.45\textwidth]{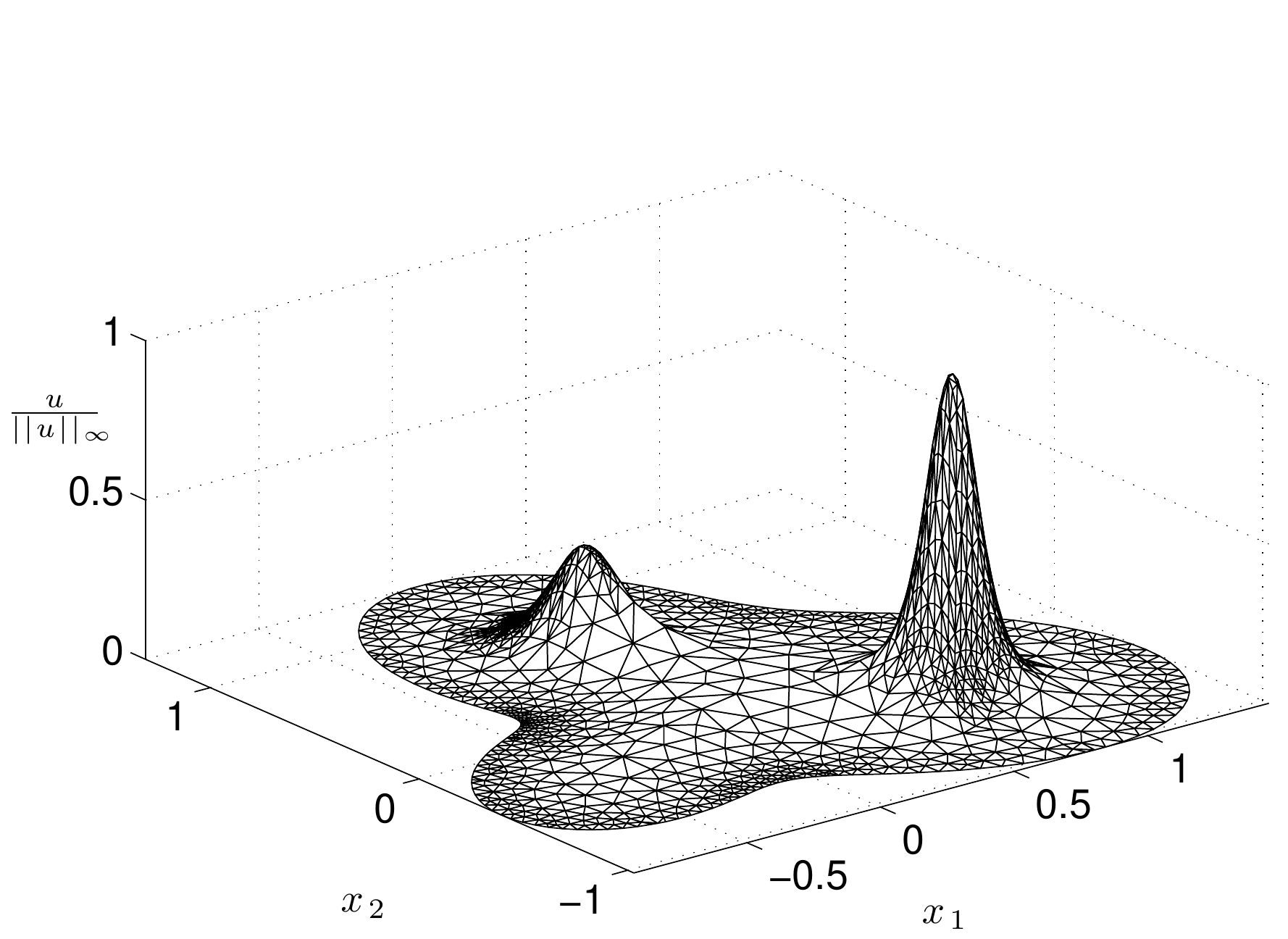} \label{fig:potato2_profile_f} }
\parbox{5in}{\caption{ Solution profiles of \eqref{intro_2b} on the region \eqref{ex:potato1} for $f(u) = (1+u)^2$, $\|u\|_{\infty}= 1\times 10^{3}$ and a range of $\eps$. In panel (a), the solution profile is shown for the case where $\al{T_{\eps}}<T_{\mathcal{S}}$ and blow-up occurs on $\omega(\al{T_{\eps}})$. The structure of $\omega(t)$ is visible as well as the undulations along it due to boundary curvature effects. In panels (b,e,f), solution profiles corresponding to $\eps$ values on segments $\mathrm{I}$, $\mathrm{II}$ and $\mathrm{III}$ of Fig.~\ref{fig:ex_potato2_b} are shown. In panels (c,e), the two peak blow-up solutions very close to the critical $\eps$ values corresponding to the transition from $\mathrm{I}$ to $\mathrm{II}$ and $\mathrm{II}$ to $\mathrm{III}$, are shown.   \label{fig:potato2_profile}}}
\end{figure}

In Fig.~\ref{fig:ex_potato2_b}, the numerically obtained blow-up points for different values of $\eps$ are shown overlaid on the skeleton $\mathcal{S}_{\Omega}$. The direction of the arrows indicate increasing values of $\eps$. Note that for this example, the relevant portion of $\mathcal{S}_{\Omega}$ is composed of three distinct branches, labelled $\mathrm{I}$, $\mathrm{II}$ and $\mathrm{III}$. The blow-up points corresponding to the smallest values of $\eps$ are to the left of \al{segment} $\mathrm{I}$. Indeed, $\al{T_{\eps}}<T_{\mathcal{S}}$ for these values of $\eps$ and the asymptotic theory predicts that blow-up occurs on discrete points of $\omega(\al{T_{\eps}})$ selected by the boundary points of largest curvature (cf. Fig.~\ref{fig:potato2_profile_a}).

As the value of $\eps$ increases, $\al{T_{\eps}}$ exceeds $T_{\mathcal{S}}$, and blow-up occurs initially on segment $\mathrm{I}$ of $\mathcal{S}_{\Omega}$, in agreement with the asymptotic theory (cf. Fig.~\ref{fig:potato2_profile_b}). In addition, the asymptotic theory predicts that blow-up should occur at $x\in\mathcal{S}_{\Omega}$ such that $s_{\Omega}(x) = \eta_0 \phi(\al{T_{\eps}};\eps)$. As the shading in Fig.~\ref{fig:ex_potato2_a} represents values of $s_{\Omega}$, we see that the condition $s_{\Omega}(x)= \eta_0 \phi(\al{T_{\eps}};\eps)$ can be satisfied by multiple $x\in\mathcal{S}_{\Omega}$ so that the leading order theory predicts multiple simultaneous blow up over a range of $\eps$. However, the leading order terms in \eqref{eqn:2d:biharm:3} only predict these points to be local maxima while higher order terms in the expansion \eqref{eqn:2d:biharm:3}, relating to the boundary curvature, will further increase the magnitude of the solution at certain points. It therefore follows that those discrete points with the largest value will be selected by the dynamics of the PDE for blow-up. Consequently, the interaction of the leading and first order terms can result in the blow-up location switching between the three segments of $\mathcal{S}_{\Omega}$ as $\eps$ is increased.

However, as the blow-up point jumps from segment $\mathrm{I}$ to $\mathrm{II}$, there is a finite value of $\eps$ for which multiple blow-up occurs and similarly with the transition from $\mathrm{II}$ to $\mathrm{III}$ (cf. Fig.~\ref{fig:potato2_profile_c}-\ref{fig:potato2_profile_e}). From  Fig.~\ref{fig:potato2_profile_b}, \ref{fig:potato2_profile_d}, we see that when the blow-up occurs on segment $\mathrm{I}$ and $\mathrm{II}$ of $\mathcal{S}_{\Omega}$ respectively, the peak on the subsequent segment of $\mathcal{S}_{\Omega}$ is developing and will eventually form the blow-up point for larger $\eps$.

\section{Three Dimensions}\label{sec:3d}

The analysis developed so far explains the multiple singularity phenomenon of \eqref{intro_2b} by means of a propagating non-monotone boundary layer. Propagating boundary effects from distal segments of $\partial\Omega$ combine to raise the solution value at certain points in $\Omega$ which are in turn selected by the dynamics of the PDE for singularity. This understanding naturally extends to three and higher dimensions, which we now investigate briefly for \eqref{intro_2b} on the cubic region $\Omega = [-1,1]^3$. Extending the analogy of one and two dimensions, we expect there exists an $\eps_c$ such that as $\eps$ increases through $\eps_c$, the multiplicity of singularities increases. As the cube has eight corners, the multiplicity can be expected to go from eight to one at this threshold. Again, we assume $\eps$ is small enough so that global solutions are not present.  
  \begin{figure}[htbp]
  \centering
  \subfigure[$\eps=0.14$]{\includegraphics[width=0.45\textwidth]{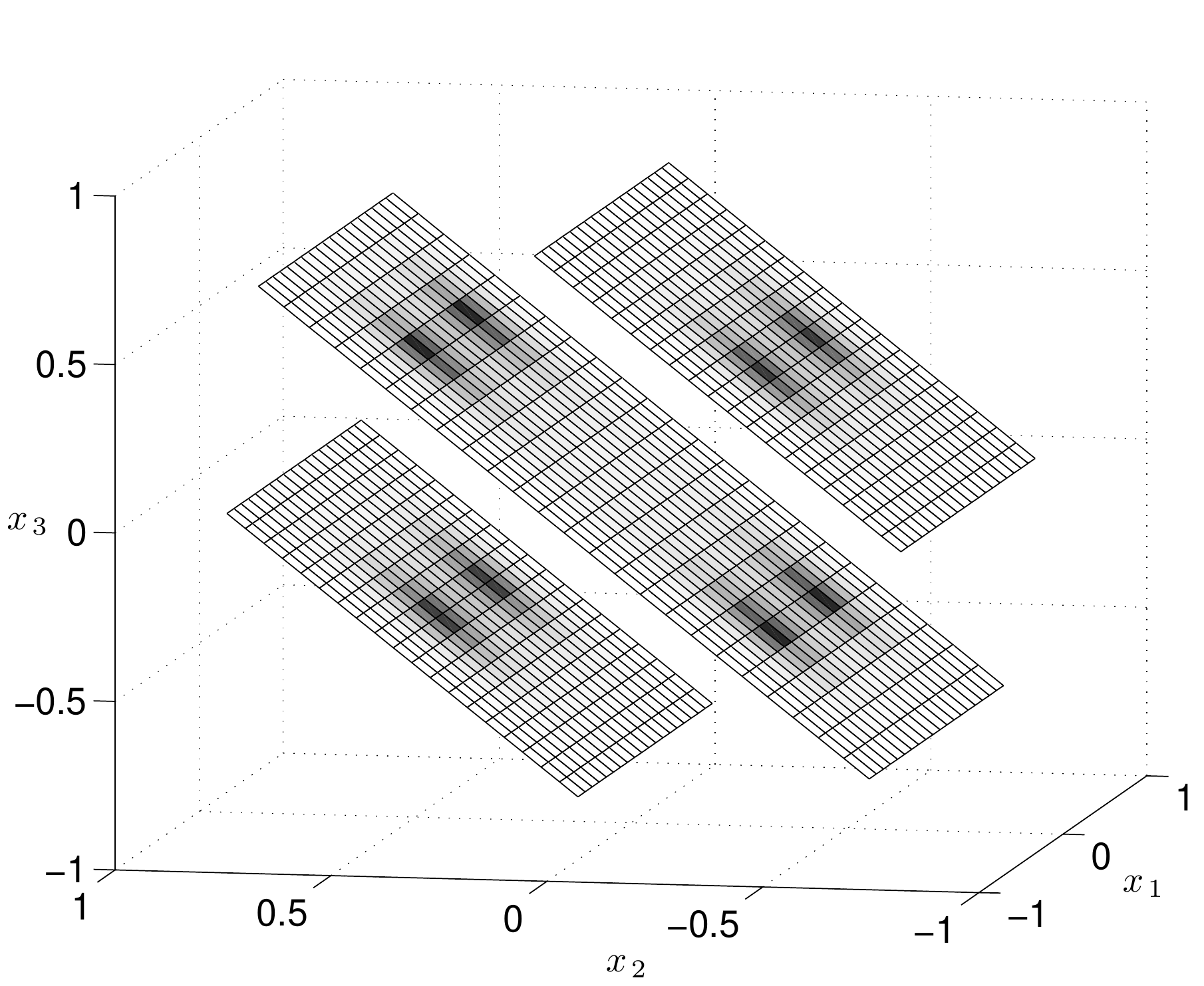}}
\quad     
\subfigure[$\eps=0.2$]{\includegraphics[width=0.45\textwidth]{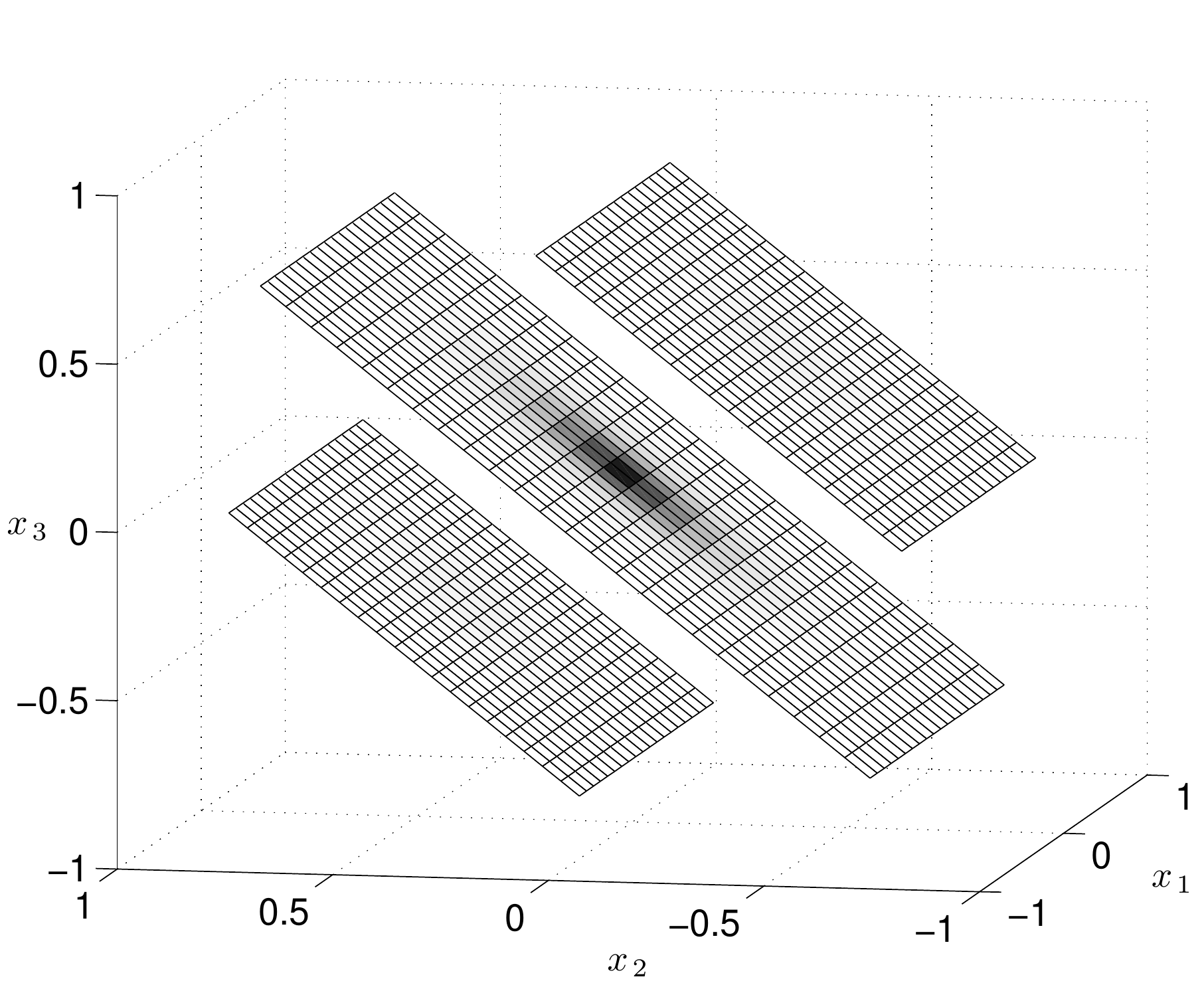}}
\parbox{5in}{\caption{Simulations of \eqref{intro_2b} for the cube $\Omega = [-1,1]^3$ with $f(u) = (1+u)^2$ and $\|u\|_{\infty}=5\times10^2$.  The solution is visualized along three parallel planes intersecting the volume. In panel (a) the profile is shown for $\eps=0.14$ in which case, we have blow-up at eight points. In panel (b), the profile is shown for $\eps=0.2$ an blow-up occurs at the origin.\label{fig:3d}}}
  \end{figure}
  
To reduce \eqref{intro_2b} to a discrete problem in the cubic region $\Omega = [-1,1]^3$, a finite difference method is applied with uniform grid spacing $h=0.05$. This relatively coarse grid cannot be expected to give accurate quantitative agreement and is mainly useful for observing the emergence of the multiple singularities. In Fig.~\ref{fig:3d}, the numerical solution is visualized on three parallel planes with normal vectors $(0,1,1)$. In the left panel, the solution for $\eps=0.14$ is shown for $\|u\|_{\infty}= 5\times10^2$ and can be seen to concentrate on eight distinct points.

It seems natural to conjecture the following from the observed blow-up behavior in one, two and three dimensions: There exists an $L_c(n)$ such that for $\Omega=[-L,L]^n$ and $\psi=0$, \alt{problem} \eqref{intro_1b} exhibits blow-up at $2^n$ distinct points whenever $L>L_c(n)$ while for $L\leq L_c(n)$, blow-up occurs uniquely at the origin. It would be interesting to develop the corresponding analogies of $\omega(t)$ and $\mathcal{S}_{\Omega}$ for bounded dimensional regions to resolve this problem and describe the possible singularity sets for \eqref{intro_1b} under a variety of bounded higher dimensional regions.

 \section{Conclusions}\label{sec:conclusion}

This paper has focussed on the exhibition and explanation of a new and interesting multiple blow-up phenomenon in fourth order parabolic equations. Through a small amplitude asymptotic analysis of the \alt{problem}
  \begin{equation}\label{outro_1}
 \left\{\begin{array}{ll} u_t = -\eps^4 \Delta^2 u + f(u), & \al{ (x,t)\in\Omega_{\al{T_{\eps}}} }; \\[5pt]
 u =\partial_n u = 0, & \al{(x,t)\in\partial\Omega_{\al{T_{\eps}}} };\\[5pt]
 u =  0, & \al{(x,t)\in\Omega_0,} \end{array}\right. 
  \end{equation}
  in spatial dimensions one and two, we have demonstrated and explained how singularities can form simultaneously at multiple points in the domain, \al{in the absence of noise}. The essence of the phenomenon is a non-monotone profile in a stretching boundary layer (cf. Fig.~\ref{fig:two profiles}) which acts to concentrate the solution on certain discrete points in $\Omega$. In the formulation \eqref{outro_1}, the parameter $\eps = L^{-1}$ acts as a length scale for the domain. A consequence of the analysis, is a geometric framework for predicting the singularity set of \eqref{outro_1} for general regions $\Omega$ in one and two spatial dimensions. As seen in the example of \S\ref{sec:ex_potato}, the singularity set has a delicate dependence of the geometry of $\Omega$ and the parameter $\eps$ which can be understood with the theory developed in the present work. \al{For domains with symmetries, we generally observe that multiple singularities are possible for a range of $\eps$ values, while for asymmetric domains (cf. \S \ref{sec:ex_potato}), multiple singularities tend to occur at fixed values of $\eps$ only.}
  
In addition, the asymptotic theory developed here is applicable to the classical second order parabolic semi-linear \alt{problem}
  \begin{equation}\label{outro_2}
 \left\{\begin{array}{ll} u_t = \eps^2 \Delta u + f(u), &\al{ (x,t)\in\Omega_{\al{T_{\eps}}} }; \\[5pt]
 u = 0, & \al{(x,t)\in\partial\Omega_{\al{T_{\eps}}} };\\[5pt]
 u =  0, &\al{(x,t)\in\Omega_0.}\end{array}\right.
  \end{equation}
It accounts for the absence of the multiple singularity phenomena in \eqref{outro_2} through the solution profile in a boundary layer near $\partial\Omega$. In the second order problem, the solution profile is monotone increasing towards its limiting value, whilst in the fourth order case there is an oscillatory approach and therefore a global max. There cannot be an overshoot in the second order problem as there is in the fourth order case, because the former has a maximum principle.
     
In addition, the leading order asymptotic analysis predicts $x_c = \max_{x\in\Omega} \mathrm{d}(x,\partial\Omega)$ to be the point(s) \al{favored} by the dynamics of the PDE for singularity. Underpinning this prediction is the intuition that if $u=0$ on $\partial\Omega$, and $u$ is monotone increasing away from $\partial\Omega$, then points furthest from $\partial\Omega$ will have larger value and are consequently more likely to be selected for singularity by the dynamics of the PDE. To the author's knowledge, such results have previously been established for radially symmetric regions only \cite{AFMC85}. 

The phenomenon described in the present work fits into a family of very interesting and unexpected solution \al{behaviors} associated with higher order PDEs. This \al{behavior} is particularly apparent when contrasted against the well understood second order case. 

There are many interesting avenues of future work which can potentially emanate from this work. First, a rigorous justification of the phenomena described herein is highly desirable - it is possible that the \emph{order-preserving majorizing equation} developed in \cite{GALAK2002} can provide a starting point for the analysis. 

Second, in some degenerate examples, for example a stadium composed of a rectangle with semi-circular end pieces (cf. \cite{ALS2013}), the predictive power of the asymptotic theory can be reduced by the jump in the curvature along $\partial\Omega$. It would therefore be very interesting to study the blow-up set of \eqref{intro_1b} as some regularity conditions on the curvature $\kappa$ of $\partial\Omega$ are relaxed. In addition, how robust is the predictive accuracy of the asymptotic theory in such degenerate cases?

\section*{Acknowledgments}
The author acknowledges many useful discussions with J. Lega and assistance with the implementation of numerical routines from FJ. Sayas. Financial support from the Carnegie Trust for the Universities of Scotland is greatly appreciated.

\end{document}